%% file: main.tex
\documentclass[10pt,CJK]{article}
\usepackage{amsfonts,amsmath,amssymb,amsthm}
\usepackage{mathrsfs}
\usepackage{latexsym}
\usepackage{graphicx}
\usepackage{fancyhdr,cite}
\usepackage{indentfirst}
\usepackage{color}

\theoremstyle{plain}                       
\newtheorem{lemma}{Lemma}[section]
\newtheorem{thm}[lemma]{Theorem}
\newtheorem{cor}[lemma]{Corollary}
\newtheorem{remark}[lemma]{Remark}
\newtheorem{definition}[lemma]{Definition}

\theoremstyle{remark}
\newtheorem*{pf}{Proof}

\numberwithin{equation}{section}

\def\Xint#1{\mathchoice
  {\XXint\displaystyle\textstyle{#1}}%
  {\XXint\textstyle\scriptstyle{#1}}%
  {\XXint\scriptstyle\scriptscriptstyle{#1}}%
  {\XXint\scriptscriptstyle\scriptscriptstyle{#1}}%
  \!\int}
\def\XXint#1#2#3{{\setbox0=\hbox{$#1{#2#3}{\int}$}
  \vcenter{\hbox{$#2#3$}}\kern-.5\wd0}}

\def\dashint{\Xint-}

\newcommand{\myu}[1]{\tilde{#1}_0}
\newcommand{\myl}[2]{\mathcal{#1}_{#2}}
\newcommand{\mys}[1]{#1_\varepsilon}

\begin{document}

\include{article_3}


\clearpage
\end{document}

%% file: article_3.tex
\allowdisplaybreaks
\pagestyle{myheadings}\markboth{$~$ \hfill {\rm Q. Xu,} \hfill $~$} {$~$ \hfill {\rm  } \hfill$~$}
\author{Qiang Xu
\thanks{Email: xuqiang09@lzu.edu.cn.}
\thanks{This work was supported by the National Natural Science Foundation of China (Grant No. 11471147).}
\\
School of Mathematics and Statistics, Lanzhou University, \\
Lanzhou, Gansu 730000, PR China.}

\title{\textbf{Convergence Rates for General Elliptic Homogenization Problems in Lipschitz Domains} }
\maketitle
\begin{abstract}
The paper extends the results obtained by C. Kenig, F. Lin and Z. Shen in \cite{SZW2}
to more general elliptic homogenization problems in two perspectives: lower order terms in the operator and no smoothness on the coefficients.
We do not repeat their arguments.
Instead we find the new weighted-type estimates for the smoothing operator at scale $\varepsilon$,
and combining some techniques developed by Z. Shen in \cite{SZW12} leads to our main results.
In addition, we also obtain sharp $O(\varepsilon)$ convergence rates in $L^{p}$ with $p=2d/(d-1)$,
which were originally established by Z. Shen for elasticity systems in \cite{SZW12}.
Also, this work may be regarded as the extension of \cite{TS,TS2} developed by T. Suslina
concerned with the bounded Lipschitz domain.
\end{abstract}

\section{Introduction and main results}
In this paper, we study convergence rates in periodic homogenization theory for general linear elliptic systems with Dirichlet or Neumann  boundary conditions in a bounded Lipschitz domain.
More precisely, we consider the following operators depending on a parameter $\varepsilon > 0$,
\begin{eqnarray*}
\mathcal{L}_{\varepsilon} =
-\text{div}\big[A(x/\varepsilon)\nabla +V(x/\varepsilon)\big] + B(x/\varepsilon)\nabla +c(x/\varepsilon) +\lambda I
\end{eqnarray*}
where $\lambda\geq 0$ is a constant, and $I=(e^{\alpha\beta})$ is an identity matrix.
Let $d\geq 3$, $m\geq 1$, and $1 \leq i,j \leq d$ and $1\leq \alpha,\beta\leq m$.
Suppose that $A = (a_{ij}^{\alpha\beta})$, $V=(V_i^{\alpha\beta})$, $B=(B_i^{\alpha\beta})$, $c=(c^{\alpha\beta})$ are real measurable functions,
satisfying the following conditions:
\begin{itemize}
\item the uniform ellipticity condition
\begin{equation}\label{a:1}
 \mu |\xi|^2 \leq a_{ij}^{\alpha\beta}(y)\xi_i^\alpha\xi_j^\beta\leq \mu^{-1} |\xi|^2,
 \quad \text{for}~y\in\mathbb{R}^d,~\text{and}~\xi=(\xi_i^\alpha)\in \mathbb{R}^{md},~\text{where}~ \mu>0;
\end{equation}
 (The summation convention for repeated indices is used throughout.)
\item the periodicity condition
\begin{equation}\label{a:2}
A(y+z) = A(y),~~ V(y+z) = V(y),~~ B(y+z) = B(y),~~ c(y+z) = c(y),
~~\text{for}~y\in \mathbb{R}^d ~\text{and}~ z\in \mathbb{Z}^d;
\end{equation}
\item the boundedness condition
\begin{equation}\label{a:3}
 \max\big\{\|V\|_{L^{\infty}(\mathbb{R}^d)},~\|B\|_{L^{\infty}(\mathbb{R}^d)},~\|c\|_{L^{\infty}(\mathbb{R}^d)}\big\}
 \leq \kappa,\qquad \text{where}~\kappa>0.
\end{equation}
\end{itemize}

Throughout this paper, we always assume $\Omega$ is a bounded Lipschitz domain, and $r_0$ denotes the diameter of $\Omega$, unless otherwise stated.
Let $\mathcal{L}_0$ be the homogenized operator associated with $\mathcal{L}_\varepsilon$,
which is expressed by
\begin{equation*}
 \mathcal{L}_0 = -\text{div}(\widehat{A}\nabla+ \widehat{V}) + \widehat{B}\nabla + \widehat{c} + \lambda I,
\end{equation*}
where $\widehat{A}$, $\widehat{V}$, $\widehat{B}$ and $\widehat{c}$
are the constant coefficients, formulated in $\eqref{f:2.1}$.
We assume that $u_\varepsilon,u_0\in H^1(\Omega;\mathbb{R}^m)$ are the weak solutions to
the Dirichlet problems $(\mathbf{DH_\varepsilon})$ and
$(\mathbf{DH_0})$.
If the coefficients of $\mathcal{L}_\varepsilon$ satisfy $\eqref{a:1}-\eqref{a:3}$, then it is well known
that $u_\varepsilon \rightharpoonup u_0$ weakly in $H^1(\Omega;\mathbb{R}^m)$
and strongly in $L^2(\Omega;\mathbb{R}^m)$. Note that
the problem $(\mathbf{DH_0})$ is usually referred to as the homogenized one
of $\mathbf{(DH_\varepsilon)}$, and the related literatures could be found in \cite{ABJLGP,VSO,QXS}.
The primary purpose of this paper is to investigate
the rate of convergence of $\|u_\varepsilon-u_0\|_{L^2(\Omega)}$, as $\varepsilon\to 0$,
in a bounded Lipschitz domain $\Omega\subset \mathbb{R}^d$. As a consequence,
we find an universal way to handle convergence rates
for elliptic operators with rapidly oscillating periodic coefficients under Dirichlet
or Neumann boundary conditions.

\subsection{Main results}

\begin{thm}[Dirichlet condition]\label{thm:1.1}
Suppose that the coefficients of $\mathcal{L}_\varepsilon$ satisfy $\eqref{a:1}-\eqref{a:3}$ and $A=A^*$.
Let $u_\varepsilon$ and $u_0$ be the weak solutions of the Dirichlet problems:
\begin{equation}\label{pde:1.2}
(\mathbf{DH_\varepsilon})\left\{
\begin{aligned}
\mathcal{L}_\varepsilon(u_\varepsilon) &= F &\quad &\text{in}~~\Omega, \\
 u_\varepsilon &= g &\quad&\text{on} ~\partial\Omega,
\end{aligned}\right.
\qquad
(\mathbf{DH_0})\left\{
\begin{aligned}
\mathcal{L}_0(u_0) &= F &\quad &\text{in}~~\Omega, \\
 u_0 &= g &\quad&\text{on} ~\partial\Omega,
\end{aligned}\right.
\end{equation}
with $F\in L^2(\Omega;\mathbb{R}^m)$ and $g\in H^1(\partial\Omega;\mathbb{R}^m)$.
Then we have
\begin{equation}\label{pri:1.1}
\big\|u_\varepsilon - u_0\big\|_{L^2(\Omega)}
\leq C\varepsilon\ln(r_0/\varepsilon)\Big\{\|F\|_{L^2(\Omega)}+\|g\|_{H^1(\partial\Omega)}\Big\}.
\end{equation}
Moreover, if $u_0\in H^2(\Omega;\mathbb{R}^m)$, then for $p=\frac{2d}{d-1}$ we have
\begin{equation}\label{pri:1.2}
\big\|u_\varepsilon - u_0\big\|_{L^p(\Omega)}\leq C\varepsilon \|u_0\|_{H^2(\Omega)},
\end{equation}
where $C$ depends only on $\mu,\kappa,m,d$ and $\Omega$.
\end{thm}

Before stating the Neumann boundary problem, we denote the conormal derivative operator
with respect to $\mathcal{L}_\varepsilon$ on $\partial\Omega$ by
\begin{equation}\label{op:2}
\mathcal{B}_\varepsilon
 = n\cdot\big[A(x/\varepsilon)\nabla + V(x/\varepsilon)\big],
\end{equation}
where $n=(n_1,\cdots,n_d)$ is the outward unit normal vector to $\partial\Omega$.
Its homogenized operator is $\mathcal{B}_0 = n\cdot[\widehat{A}\nabla +\widehat{V}]$,
and the details could be found in \cite[pp.7]{QXS1}.

\begin{thm}[Neumann condition]\label{thm:1.2}
Suppose that the coefficients of $\mathcal{L}_\varepsilon$ satisfy $\eqref{a:1}-\eqref{a:3}$ and $A=A^*$.
Let $u_\varepsilon$ and $u_0$ be the weak solutions of the Neumann problems:
\begin{equation}\label{pde:1.3}
(\mathbf{NH_\varepsilon})\left\{
\begin{aligned}
\mathcal{L}_\varepsilon(u_\varepsilon) &= F &\quad &\emph{in}~~\Omega, \\
 \mathcal{B}_\varepsilon(u_\varepsilon) &= h &\quad&\emph{on} ~\partial\Omega,
\end{aligned}\right.
\qquad
(\mathbf{NH_0})\left\{
\begin{aligned}
\mathcal{L}_0(u_0) &= F &\quad &\emph{in}~~\Omega, \\
 \mathcal{B}_0(u_0) &= h &\quad&\emph{on} ~\partial\Omega,
\end{aligned}\right.
\end{equation}
where $F\in L^2(\Omega;\mathbb{R}^m)$ and $h\in L^2(\partial\Omega;\mathbb{R}^m)$ satisfy the
compatibility condition $\eqref{C:2.1}$, then we have
\begin{equation}\label{pri:1.3}
\big\|u_\varepsilon - u_0\big\|_{L^2(\Omega)}
\leq C\varepsilon\ln(r_0/\varepsilon)\Big\{\|F\|_{L^2(\Omega)}+\|h\|_{L^2(\partial\Omega)}\Big\}.
\end{equation}
Moreover, if $u_0\in H^2(\Omega;\mathbb{R}^m)$, then for $p=\frac{2d}{d-1}$ we have
\begin{equation}\label{pri:1.4}
\big\|u_\varepsilon - u_0\big\|_{L^p(\Omega)}\leq C\varepsilon \|u_0\|_{H^2(\Omega)},
\end{equation}
where $C$ depends only on $\mu,\kappa,m,d$ and $\Omega$.
\end{thm}

Our main results are quite similar to those obtained by C. Kenig, F. Lin and Z. Shen
for $L_\varepsilon$ with Dirichlet or Neumann boundary conditions (see \cite[Theorem 1.1, Theorem 1.2]{SZW2}),
where $L_\varepsilon = -\text{div}[A(x/\varepsilon)\nabla]$.
Compared to theirs,
an obvious progress is that the estimates $\eqref{pri:1.1}$ and $\eqref{pri:1.3}$ do not rely on any smoothness assumption
on the coefficients of $\mathcal{L}_\varepsilon$.
Besides, the operator $\mathcal{L}_\varepsilon$ investigated here is more complicated,
which requires the use of one more corrector $\chi_0$ produced by the coefficient $V$.
Although it would not bring the essential difficulty in most cases (see \cite{QXS,QXS1} for the recent works),
some subtle tools such as the radial maximal operator
are still necessary here to control the behavior of $u_0$ near $\partial\Omega$ or far away from $\partial\Omega$
at scale $\varepsilon$.
Since the nontangential maximal function estimates for $\mathcal{L}_\varepsilon$
in Lipschitz domains have not been established yet,
we can not count on the methods developed in \cite{SZW2} to derive the estimates $\eqref{pri:1.1}$ and $\eqref{pri:1.3}$.
Fortunately, some new findings permit us to transfer all the estimates from $\mathcal{L}_\varepsilon$ to $\mathcal{L}_0$,
and obviously the regularity theories related to $\mathcal{L}_0$ are good enough to be employed.
So, before giving the formal proofs of Theorems $\ref{thm:1.1}$, $\ref{thm:1.2}$,
it is instructive to sketch the main ideas.
For convenience, we take the Dirichlet problem as an example.

Before proceeding further, it is necessary to introduce some notation to ease the later statements.
\subsection{Notation in the paper}\label{sec:1.1}
\begin{itemize}
\item distance function $\delta(x)=\text{dist}(x,\partial\Omega)$, where $x\in\Omega$.
If $x\in\mathbb{R}^d\setminus\Omega$, then we set $\delta(x) = 0$;
\item boundary layer $\Omega\setminus\Sigma_r$,
where $\Sigma_{r}=\{x\in\Omega:\text{dist}(x,\partial\Omega)>r\}$ with $r>0$;
\item cut-off function $\psi_r$ (associated with $\Sigma_r$), satisfying
\begin{equation}\label{def:2.5}
\psi_{r} =1 \quad\text{in}\quad \Sigma_{2r},
\qquad\psi_{r} = 0 \quad\text{outside}\quad \Sigma_{r}, \qquad\text{and}\quad
|\nabla\psi_{r}|\leq C/r;
\end{equation}
\item  level set $S_r = \big\{x\in\Omega:\text{dist}(x,\partial\Omega) = r\big\}$;
\item internal diameter $r_{00} = \max\{r>0:B(x,r)\subset\Omega\},\forall x\in S_r$,
and $c_0 = r_{00}/10$ is referred to as the layer constant.
\item the weighted-type norms:
\begin{equation}
 \|f\|_{L^2(\Sigma_r;\delta)} = \Big(\int_{\Sigma_r}|f(x)|^2\delta(x)dx\Big)^{1/2},
 \qquad
 \|f\|_{L^2(\Sigma_r;\delta^{-1})} = \Big(\int_{\Sigma_r}|f(x)|^2\delta^{-1}(x)dx\Big)^{1/2}.
\end{equation}
\end{itemize}
We denote
$\|f\|_{H^k(\Sigma_r,\delta)}$ by $\sum_{i=0}^k\|\nabla^kf\|_{L^2(\Sigma_r;\delta)}$,
where $k$ is a positive integer, and
$\|\nabla^0f\|_{L^2(\Sigma_r;\delta)} = \|f\|_{L^2(\Sigma_r;\delta)}$.
We take a similar way to define $\|f\|_{H^k(\Sigma_r,\delta^{-1})}$.

\subsection{Outline of the proof of $\eqref{pri:1.2}$}

We first introduce the recent progress on convergence rates made by Z. Shen in \cite{SZW12} and by T. Suslina in \cite{TS,TS2}
through outlining the proof for the estimate $\eqref{pri:1.2}$.
To obtain the estimates $\|u_\varepsilon-u_0\|_{L^p(\Omega)}=O(\varepsilon)$ with $p=2,\text{or}~2d/(d-1)$,
it is sufficient to prove
$\|w_\varepsilon\|_{L^p(\Omega)} = O(\varepsilon)$, 
where $w_\varepsilon$ is the first order approximating of $u_\varepsilon$, defined by
\begin{equation*}
 w_\varepsilon^\beta = u_\varepsilon^\beta - u_0^\beta - \varepsilon\sum_{k=0}^d\chi_{k}^{\beta\gamma}(x/\varepsilon)\varphi_k^\gamma.
\end{equation*}
Note that $\chi_{k}$ with $k=0,\cdots,d$ are correctors,
and $\varphi=(\varphi_k^\gamma)\in H^1(\mathbb{R}^d;\mathbb{R}^{m(d+1)})$
can be later fixed by the concrete target. Before estimating
$\|w_\varepsilon\|_{L^p(\Omega)}$,
we need to calculate the quantity $\|w_\varepsilon\|_{H^1(\Omega)}$.
For this purpose, it is natural to consider
what equation $w_\varepsilon$ satisfies.
According to the fact that $\mathcal{L}_\varepsilon(u_\varepsilon)=\mathcal{L}_0(u_0)$
in $\Omega$ and $u_\varepsilon = u_0$ on $\partial\Omega$,
it is not hard to check that $w_\varepsilon$ satisfies the following Dirichlet boundary value problem
\begin{equation*}
 \left\{\begin{aligned}
  \mathcal{L}_\varepsilon(w_\varepsilon) & =  -\text{div}(\tilde{f}) + \tilde{F}  & \quad &\text{in}~~\Omega,\\
  w_\varepsilon& = \varepsilon\small\sum_{k=0}^{d}\chi_{k,\varepsilon}\varphi_k &\quad &\text{on} ~\partial\Omega,
 \end{aligned}\right.
\end{equation*}
where $\tilde{f}$ and $\tilde{F}$ are really complicated and we will show them later.
Thus the quantity $\|w_\varepsilon\|_{H^1(\Omega)}$ determined by $\tilde{f}$ and $\tilde{F}$
follows from the $H^1$ estimates, and then one may estimate
$\|w_\varepsilon\|_{L^p(\Omega)}$ by a duality argument (see \cite{TS,TS2}).

Precisely speaking, we need to consider the related dual problems as follows.
For any $\Phi\in L^q(\Omega;\mathbb{R}^m)$ with $q=2$ or $2d/(d+1)$,
we say  $\phi_\varepsilon$ and $\phi_0$ are the solutions to the ``adjoint Dirichlet problems''
associated with $\eqref{pde:1.2}$,
if $\phi_\varepsilon$ and $\phi_0$ respectively solve
\begin{equation}\label{pde:2.6}
(\mathbf{DH_\varepsilon^*})\left\{
\begin{aligned}
\mathcal{L}_\varepsilon^*(\phi_\varepsilon) &= \Phi &\quad &\text{in}~~\Omega, \\
 \phi_\varepsilon &= 0 &\quad&\text{on} ~\partial\Omega,
\end{aligned}\right.
\qquad
(\mathbf{DH_0^*})\left\{
\begin{aligned}
\mathcal{L}_0^*(\phi_0) &= \Phi &\quad &\text{in}~~\Omega, \\
 \phi_0 &= 0 &\quad&\text{on} ~\partial\Omega,
\end{aligned}\right.
\end{equation}
where $\mathcal{L}_\varepsilon^*$ is
the adjoint operator associated with $\mathcal{L}_\varepsilon$.
We are now in the position to show the formula
\begin{equation*}
 \int_\Omega w_\varepsilon\Phi dx
 = \int_\Omega \tilde{f}\cdot\nabla\phi_\varepsilon dx
 + \int_\Omega \tilde{F}\phi_\varepsilon dx,
\end{equation*}
provided $\varphi_k$ with $k=0,\cdots,d$ are supported in $\Omega$. We mention that this equality
follows from the second Green's formula associated with
$\mathcal{L}_\varepsilon$ and $\mathcal{L}_\varepsilon^*$,
and then implies the estimate of $\|w_\varepsilon\|_{L^p(\Omega)}$ in terms of $\tilde{f}$ and $\tilde{F}$.

Until now, the key issue concerning the estimates of $\|w_\varepsilon\|_{H^1(\Omega)}$ and $\|w_\varepsilon\|_{L^p(\Omega)}$
is all related to $\tilde{f}$ and $\tilde{F}$.
We present their exact formulas in $\eqref{eq:2.10}$ and $\eqref{eq:2.11}$
since they are too lengthy to be shown here.
The techniques to handle $\tilde{f}$ and $\tilde{F}$ have already been in \cite{SZW12,TS,TS2},
and we directly use the notation $\mathcal{K}$, $\mathcal{I}$ in $\eqref{eq:2.11}$ to
indicate where the related techniques are applied.
Although a little more new notation appears,
if the reader is experienced, one may easily discover
$\mathcal{K}$, $\mathcal{I}$ and $\nabla_j u_0 - \varphi_j$ are the most difficult terms
in $\eqref{eq:2.10}$, where $\nabla_j=\partial/\partial x_j$ with $j=1,\cdots,d$.
To reach our goal, a little trick is to set $\varphi_0 = S^2_\varepsilon(\psi_{4\varepsilon}u_0)$
and $\varphi_j = S^2_\varepsilon(\psi_{4\varepsilon}\nabla_j u_0)$,
where $S_\varepsilon$ is a smoothing operator at scale $\varepsilon$,
and $\psi_{4\varepsilon}$
is a cut-off function supported in
$\Sigma_{4\varepsilon}$.
Here we use cut-off function to avoid analyzing the behavior of
$w_\varepsilon$ on $\partial\Omega$, while the cost paid is that we have to estimate the quantity
$\|u_0\|_{H^1(\Omega\setminus\Sigma_{4\varepsilon})}$.

The first result is $\|w_\varepsilon\|_{H^1(\Omega)} = O(\varepsilon^{\frac{1}{2}})$.
In the proof, we borrow the methods from T. Suslina \cite{TS,TS2} to handle the term $\mathcal{K}$,
and the techniques from Z. Shen \cite{SZW12}
to deal with the terms $\mathcal{I}$ and $\nabla_j u_0 -\varphi_j$.

Furthermore, inspired by T. Suslina \cite{TS,TS2} we employ the duality argument to accelerate the convergence rate.
Resetting $w_\varepsilon$ we then have (see Lemma $\ref{lemma:2.11}$)
\begin{equation}\label{pri:1.5}
\Big|\int_\Omega w_\varepsilon\Phi dx\Big|
\leq C\Big\{\|u_0\|_{H^1(\Omega\setminus\Sigma_{8\varepsilon})}\|\phi_\varepsilon\|_{H^1(\Omega\setminus\Sigma_{9\varepsilon})}
+ \varepsilon\|u_0\|_{H^2(\Omega)}\|\phi_\varepsilon\|_{H^1(\Omega)}\Big\}.
\end{equation}
Obviously, the next thing is to show the quantities
$\|u_0\|_{H^1(\Omega\setminus\Sigma_{8\varepsilon})}$
and $\|\phi_\varepsilon\|_{H^1(\Omega\setminus\Sigma_{9\varepsilon})}$.
The first one is done with the aid of the nontangential maximal function coupled with radial maximal function to
control the behavior of $\nabla u_0$ and $u_0$ on $\Omega\setminus\Sigma_{8\varepsilon}$, respectively.
The original idea belongs to Z. Shen in \cite{SZW12}.
For the quantity $\|\phi_\varepsilon\|_{H^1(\Omega\setminus\Sigma_{9\varepsilon})}$,
we consider the auxiliary function
\begin{equation*}
\Theta_\varepsilon = \phi_\varepsilon - \phi_0 - \varepsilon \chi_{0,\varepsilon}^* S_\varepsilon^2(\psi_{20\varepsilon}\phi_0)
- \varepsilon \chi_{k,\varepsilon}^* S_\varepsilon^2(\psi_{20\varepsilon}\nabla_k\phi_0),
\end{equation*}
and obtain $\|\Theta_\varepsilon\|_{H^1(\Omega)}=O(\varepsilon^{\frac{1}{2}})$ due to the previous result,
where $\chi^*_{k}$ with $k=0,\cdots,d$ are corresponding first order correctors of $\mathcal{L}_\varepsilon^*$.
This argument is developed by T.Suslina in \cite{TS,TS2}. Then it is not hard to reach
$\|\phi_\varepsilon\|_{H^1(\Omega\setminus\Sigma_{9\varepsilon})}\leq \|\Theta_\varepsilon\|_{H^1(\Omega)}
+ \|\phi_0\|_{H^1(\Omega\setminus\Sigma_{9\varepsilon})}$ by the observation that $S_\varepsilon^2(\psi_{20\varepsilon}\phi_0)$
and $S_\varepsilon^2(\psi_{20\varepsilon}\nabla_k\phi_0)$ are supported in $\Sigma_{18\varepsilon}$
and $\Sigma_{18\varepsilon}\cap(\Omega\setminus\Sigma_{9\varepsilon})=\emptyset$.
Consequently,
a routine computation gives $\|w_\varepsilon\|_{L^p(\Omega)} = O(\varepsilon)$
under the condition of $u_0\in H^2(\Omega)$. This implies our result $\eqref{pri:1.2}$.

In fact, when $\Omega$ is merely a Lipschitz domain, the assumption $u_0\in H^2(\Omega;\mathbb{R}^m)$
is not very natural, since $\nabla u_0$ only belongs to $H^{1/2}(\Omega;\mathbb{R}^m)$
under the boundary condition $g\in H^1(\partial\Omega;\mathbb{R}^m)$ (see \cite[Lemma 4.3]{SZW2}).
For this reason,
we manage to get rid of the assumption of $u_0\in H^2(\Omega;\mathbb{R}^m)$
by using the given data in the systems $\eqref{pde:1.2}$.
Compared with the proof of $\eqref{pri:1.2}$,
the improvement of the methods in this paper focuses on proving the estimate $\eqref{pri:1.1}$.

\subsection{Sketch of the proof of $\eqref{pri:1.1}$}
In the following,
we want to explain why we consider the smoothing operator in weighted-type norms.
Perhaps the reader will be disappointed because the motivation is rooted in the tedious computation and
we owe the success to the author's luck.

Let us go back to the estimate $\eqref{pri:1.5}$. In fact, before obtaining it
the quantity $\int_\Omega w_\varepsilon\Phi dx$ is controlled by
\begin{equation*}
C\left\{
 \|u_0\|_{H^1(\Omega\setminus\Sigma_{9\varepsilon})}
 \|\phi_\varepsilon\|_{H^1(\Omega\setminus\Sigma_{9\varepsilon})}
 + \|u_0\|_{H^1(\Omega\setminus\Sigma_{9\varepsilon})}\|\phi_\varepsilon\|_{H^1(\Sigma_{4\varepsilon})}
 + \varepsilon \|u_0\|_{H^2(\Sigma_{4\varepsilon})}\|\phi_\varepsilon\|_{H^1(\Sigma_{4\varepsilon})}
 \right\}.
\end{equation*}
Due to the auxiliary function $\Theta_\varepsilon$,
the estimates of $\phi_\varepsilon$ are transformed into the corresponding ones of $\phi_0$, and
the above expression turns into
\begin{equation}\label{f:1.1}
C\left\{
 \|u_0\|_{H^1(\Omega\setminus\Sigma_{9\varepsilon})}
 \|\phi_0\|_{H^1(\Omega\setminus\Sigma_{9\varepsilon})}
 + \|u_0\|_{H^1(\Omega\setminus\Sigma_{9\varepsilon})}\|\phi_0\|_{H^1(\Sigma_{4\varepsilon})}
 + \varepsilon \|u_0\|_{H^2(\Sigma_{4\varepsilon})}\|\phi_0\|_{H^1(\Sigma_{4\varepsilon})}
 + \text{good terms}
 \right\},
\end{equation}
where ``good terms'' means the terms which do not bring essential difficulties.
An obvious advantage of this argument is that we avoid using the nontangential maximal function of $\phi_\varepsilon$ to
control its behavior near the boundary,
which opens up an opportunity to handle the operator with lower order terms
even though we have not yet established the nontangential maximal function estimates
for $(\mathbf{DH_\varepsilon})$ in Lipschitz domains.
Another important advantage is that it can accelerate the convergence rates
since $\|\Theta_\varepsilon\|_{H^1(\Omega)} = O(\varepsilon^{\frac{1}{2}})$.
Let us observe the expression $\eqref{f:1.1}$ again, and we claim that
the first term $\|u_0\|_{H^1(\Omega\setminus\Sigma_{9\varepsilon})}
 \|\phi_0\|_{H^1(\Omega\setminus\Sigma_{9\varepsilon})}$ is also good,
 since we can prove the following estimates
\begin{equation}\label{f:1.2}
 \|u_0\|_{H^1(\Omega\setminus\Sigma_{9\varepsilon})} = O(\varepsilon^{\frac{1}{2}}),
 \qquad
 \|u_0\|_{H^2(\Sigma_{4\varepsilon})} = O(\varepsilon^{-\frac{1}{2}}),
\end{equation}
in Lemma $\ref{lemma:3.5}$.
The above estimates were first obtained by Z. Shen for elasticity systems in \cite{SZW12}.
Meanwhile the estimates $\eqref{f:1.2}$ suggest that
the second and third terms of $\eqref{f:1.1}$ are the tricky ones, which
exactly arouse the inspiration for the weighted-type norms.
Inspired by the proof in \cite{SZW12},
it is not very hard to derive the following estimates,
\begin{equation*}
\|u_0\|_{H^1(\Omega\setminus\Sigma_{9\varepsilon};\delta)} = O(\varepsilon),
\qquad
\max\big\{\|\nabla^2u_0\|_{L^2(\Sigma_{4\varepsilon};\delta)},
\|\nabla^2\phi_0\|_{L^2(\Sigma_{4\varepsilon};\delta)}\big\} = O\big([\ln(r_0/\varepsilon)]^{\frac{1}{2}}\big),
\end{equation*}
as well as
\begin{equation*}
\|\phi_0\|_{H^1(\Sigma_{4\varepsilon};\delta^{-1})}  = O\big([\ln(r_0/\varepsilon)]^{\frac{1}{2}}\big).
\end{equation*}
Compared to the estimates $\eqref{f:1.2}$,
the above quantities show a noticeable improvement in the sense of the order of $\varepsilon$.
Since $\Omega$ is just a Lipschitz domain,
we can not expect to apply $H^2$ estimates to $u_0$ near the boundary.
In general, the worse the smoothness of the boundary is,
the higher the level of the technicalities employed will be.
The core idea in the proof is to use the radial maximal function,
which was employed by C. Kenig, F. Lin and Z. Shen in \cite{SZW2},
and nontangential maximal function
to control the boundary behavior of the solutions to the corresponding homogenized problems.
In addition, we find $\delta$ as the weighted function highly effective
in eliminating some singularity produced by $\nabla^2 u_0$ near the boundary of $\Omega$.
This is another crucial aspect to explain why we use the above weighted-type norms, and
we refer the reader to the proof of Lemma $\ref{lemma:3.6}$ for precise details.

Observe that if the second and third terms in $\eqref{f:1.1}$ are replaced by
\begin{equation*}
 \|u_0\|_{H^1(\Omega\setminus\Sigma_{9\varepsilon};\delta)}\|\phi_0\|_{H^1(\Sigma_{4\varepsilon};\delta^{-1})}
 \quad \text{and} \quad
 \varepsilon \|u_0\|_{H^2(\Sigma_{4\varepsilon};\delta)}\|\phi_0\|_{H^1(\Sigma_{4\varepsilon};\delta^{-1})},
\end{equation*}
respectively,
then it is not far from reaching $\|w_\varepsilon\|_{L^2(\Omega)} = O(\varepsilon\ln(r_0/\varepsilon))$
due to the duality method. This enlightens us to find the following
weighted-type inequality, 
\begin{equation}\label{pri:1.6}
\begin{aligned}
\Big|\int_\Omega w_\varepsilon\Phi dx\Big|
&\leq C\|u_0\|_{H^1(\Omega\setminus\Sigma_{9\varepsilon})}\|\phi_0\|_{H^1(\Omega\setminus\Sigma_{9\varepsilon})}
 + C\varepsilon\|u_0\|_{H^1(\Omega)}\|\phi_\varepsilon\|_{H^1(\Omega)} \\
& \quad + C\Big\{\| u_0\|_{H^1(\Omega\setminus\Sigma_{8\varepsilon};\delta)}
+\varepsilon \|\nabla u_0\|_{L^2(\Sigma_{4\varepsilon};\delta)}
+\varepsilon \|\nabla^2 u_0\|_{L^2(\Sigma_{4\varepsilon};\delta)}\Big\}\|\phi_0\|_{H^1(\Sigma_{4\varepsilon};\delta^{-1})} \\
& \quad\quad\quad+ C\Big\{\| u_0\|_{H^1(\Omega\setminus\Sigma_{8\varepsilon})}
+\varepsilon \|\nabla u_0\|_{L^2(\Sigma_{4\varepsilon})}
+\varepsilon \|\nabla^2 u_0\|_{L^2(\Sigma_{4\varepsilon})}\Big\}\\
&\quad\quad\quad\quad\quad\cdot\Big\{\|\xi_\varepsilon\|_{H^1(\Omega)} + \|\phi_0\|_{H^1(\Omega\setminus\Sigma_{20\varepsilon})}
+ \varepsilon\|\phi_0\|_{H^1(\Omega)}+ \varepsilon\|\nabla^2\phi_0\|_{L^2(\Sigma_{10\varepsilon})}\Big\},
\end{aligned}
\end{equation}
where $\xi_\varepsilon = \phi_\varepsilon
- \phi_0 -\varepsilon \chi_{0,\varepsilon}^* S_\varepsilon(\psi_{10\varepsilon}\phi_0)
- \varepsilon \chi_{k,\varepsilon}^* S_\varepsilon(\psi_{10\varepsilon}\nabla_k\phi_0)$.
It follows from Theorem $\ref{thm:3.1}$
that $\|\xi_\varepsilon\|_{H^1(\Omega)} = O(\varepsilon^{\frac{1}{2}})$.
Noting the right-hand side of $\eqref{pri:1.6}$,
the first, second and fourth terms produce the factor $O(\varepsilon)$,
while the third term contribute the factor $O(\varepsilon\ln(r_0/\varepsilon))$.
Combining them and applying the duality method finally leads to the desired estimate $\eqref{pri:1.1}$.

Now a big question is how to derive the estimate $\eqref{pri:1.6}$.
Although the computations
are quite lengthy, the central aim is to figure out the factors
$O(\varepsilon)$ or $O(\varepsilon\ln(r_0/\varepsilon))$. To do so,
the crucial ingredients are the following estimates of weighted-type,
\begin{equation}\label{f:1.3}
\|g_\varepsilon S_\varepsilon(f)\|_{L^2(\Sigma_{2\varepsilon};\delta)}
 \leq C\|g\|_{L^2(Y)}\|f\|_{L^2(\Sigma_{2\varepsilon};\delta)},
 \qquad
\|g_\varepsilon S_\varepsilon(f)\|_{L^2(\Sigma_{2\varepsilon};\delta^{-1})}
\leq C\|g\|_{L^2(Y)}\|f\|_{L^2(\Sigma_{2\varepsilon};\delta^{-1})},
\end{equation}
where $f\in L^2(\Omega)$ is supported in $\Sigma_{2\varepsilon}$,
and $g_\varepsilon(x) = g(x/\varepsilon)$ is square integrable and periodic at scale $\varepsilon$. Also,
\begin{equation}\label{f:1.4}
\|f- S_\varepsilon(f)\|_{L^2(\Sigma_{2\varepsilon};\delta)}
 \leq C\varepsilon\|\nabla f\|_{L^2(\Sigma_{\varepsilon};\delta)},
\end{equation}
where $f\in H^1(\Omega)$ is supported in $\Sigma_{\varepsilon}$.
Briefly speaking, the smoothing operator in the weighted-type norms satisfies
the similar properties as in \cite{SZW14,TS2,TS}.
We refer the reader to Lemmas $\ref{lemma:2.7}$ and $\ref{lemma:2.8}$ for the accurate statement.
We end this paragraph
with a comment that the weighted-type norms were discovered to accelerate the convergence rates at first,
and then inspire us to check the estimates $\eqref{f:1.3}$ and $\eqref{f:1.4}$.

For the Neumann problem,
the proof of $\eqref{pri:1.4}$ is quite similar to that given for $\eqref{pri:1.2}$.
A slight difference occurs when we proceed with the duality argument. Here we construct the following ``adjoint Neumann problems''
associated with $\eqref{pde:1.3}$
\begin{equation}\label{pde:2.7}
(\mathbf{NH_\varepsilon^*})\left\{
\begin{aligned}
\mathcal{L}_\varepsilon^*(\phi_\varepsilon) &= \Phi &\quad &\text{in}~~\Omega, \\
 \mathcal{B}_\varepsilon^*(\phi_\varepsilon) &= 0 &\quad&\text{on} ~\partial\Omega,
\end{aligned}\right.
\qquad
(\mathbf{NH_0^*})\left\{
\begin{aligned}
\mathcal{L}_0^*(\phi_0) &= \Phi &\quad &\text{in}~~\Omega, \\
 \mathcal{B}_0^*(\phi_0) &= 0 &\quad&\text{on} ~\partial\Omega,
\end{aligned}\right.
\end{equation}
where $\mathcal{B}_\varepsilon^*$ is the adjoint operator of $\mathcal{B}_\varepsilon$,
and $\mathcal{B}_0^*=n\cdot[\widehat{A}\nabla + \widehat{B}]$
is the homogenized operator of $\mathcal{B}_\varepsilon^*$.
The remainder of the arguments is analogous to that in Theorem $\ref{thm:1.1}$ and we omit them here.

\subsection{Source of the ideas}

Undoubtedly, the ideas used here are inspired from several sources. For example,
some weighted-type inequalities  applied to studying convergence rates have already been in \cite{SZW2},
and the smoothing operator $S_\varepsilon$ at the scale $\varepsilon$ as the improvement of
the Steklov smoothing operator was originally shown in \cite{SZW12},
in which Z. Shen has obtained the estimates like $\eqref{pri:1.2}$ and $\eqref{pri:1.4}$ for $L_\varepsilon$.
The duality argument employed here is motivated by T. Suslina in \cite{TS,TS2}, and the Steklov smoothing operator was originally applied to
homogenization problem by V.V. Zhikov in \cite{ZVVPSE}.
However the estimates $\eqref{pri:1.6}$, $\eqref{f:1.3}$ and $\eqref{f:1.4}$ presented here
seem to be new. We do not offer lengthy heuristics or motivation, but as compensation have tried
to present all the technicalities of the proofs in the later sections.

We mention that in the case of $d=2$, the correctors $\chi_k$ with $k=0,\cdots,d$ are bounded even without the smoothness assumption on
$A$ (see \cite[Section 4.4]{MGLM}). There is a simple way to derive the corresponding results of this paper,
and we will investigate this case in another place.
Finally, we remark that the convergence rates are active topics in homogenization theory, and without attempting to be exhaustive,
we refer the reader to \cite{SAT,SZ,SC,MAFHL,ABJLGP,BMSHSTA,BMSHSTA2,CDJ,ACPD,GG1,GG2,G,VSO,ODVB,SZW14,SZW13,SZW20,SAT2,QXS,QXS1,ZVVPSE}
and the references therein for more results.

This paper is organized as follows. Some definitions and remarks as well as known lemmas and theorems
are introduced in Section 2. We will prove the weighted-type inequalities for the smoothing operator in Section 3,
in which we also introduce two lemmas related to the duality methods (see Lemmas $\ref{lemma:2.10},\ref{lemma:2.11}$).
Section 4 deals with the Dirichlet problem, while Section 5 handles the Neumann problem.

\section{Preliminaries}
Define the correctors $\chi_k = (\chi_{k}^{\alpha\beta})$ associated with $\mathcal{L}_\varepsilon$ as follows:
\begin{equation}\label{pde:2.1}
\left\{ \begin{aligned}
 &L_1(\chi_0) = \text{div}(V)  \quad \text{in}~ \mathbb{R}^d, \\
 &\chi_0\in H^1_{per}(Y;\mathbb{R}^{m^2})~~\text{and}~\int_Y\chi_0 dy = 0
\end{aligned}
\right.
\end{equation}
and
\begin{equation}\label{pde:2.2}
 \left\{ \begin{aligned}
  &L_1(\chi_k^\beta + P_k^\beta) = 0 \quad \text{in}~ \mathbb{R}^d, \\
  &\chi_k^\beta \in H^1_{per}(Y;\mathbb{R}^m)~~\text{and}~\int_Y\chi_k^\beta dy = 0
 \end{aligned}
 \right.
\end{equation}
for $1\leq k\leq d$, where $Y = (0,1]^d \cong \mathbb{R}^d/\mathbb{Z}^d$, and $H^1_{per}(Y;\mathbb{R}^m)$ denotes the closure
of $C^\infty_{per}(Y;\mathbb{R}^m)$ in $H^1(Y;\mathbb{R}^m)$.
Note that $C^\infty_{per}(Y;\mathbb{R}^m)$ is the subset of $C^\infty(Y;\mathbb{R}^m)$, which collects all $Y$-periodic vector-valued functions
(see \cite[pp.56]{ACPD}). By asymptotic expansion arguments (see \cite[pp.103]{ABJLGP} or \cite[pp.31]{VSO}),
we obtain the homogenized operator $\mathcal{L}_0$, and its coefficients
$\widehat{A} = (\hat{a}_{ij}^{\alpha\beta})$, $\widehat{V}=(\hat{V}_i^{\alpha\beta})$,
$\widehat{B} = (\hat{B}_i^{\alpha\beta})$ and $\widehat{c}= (\hat{c}^{\alpha\beta})$ are given by
\begin{equation}\label{f:2.1}
\begin{aligned}
\hat{a}_{ij}^{\alpha\beta} = \int_Y \big[a_{ij}^{\alpha\beta} + a_{ik}^{\alpha\gamma}\frac{\partial\chi_j^{\gamma\beta}}{\partial y_k}\big] dy, \qquad
\hat{V}_i^{\alpha\beta} = \int_Y \big[V_i^{\alpha\beta} + a_{ij}^{\alpha\gamma}\frac{\partial\chi_0^{\gamma\beta}}{\partial y_j}\big] dy, \\
\hat{B}_i^{\alpha\beta} = \int_Y \big[B_i^{\alpha\beta} + B_j^{\alpha\gamma}\frac{\partial\chi_i^{\gamma\beta}}{\partial y_j}\big] dy, \qquad
\hat{c}^{\alpha\beta} = \int_Y \big[c^{\alpha\beta} + B_i^{\alpha\gamma}\frac{\partial\chi_0^{\gamma\beta}}{\partial y_i}\big] dy.
\end{aligned}
\end{equation}

\begin{remark}\label{re:2.1}
\emph{For simplicity of presentation, if $f$ is a periodic function, we will denote $f(x/\varepsilon)$ by $f_{\varepsilon}(x)$.
For example, we usually write
$A_\varepsilon(x) = A(x/\varepsilon)$ and $\chi_{k,\varepsilon}(x) = \chi_k(x/\varepsilon)$,
and their components follow the same simplified way as well. \textbf{Warning:}
the reader do not confuse the two types of notation:
one type is the abridged notation, the other type is the common notation such as
the solution $u_\varepsilon$, the smoothing operator $S_\varepsilon$ and the set $\Sigma_{\varepsilon}$.}
\end{remark}

\begin{definition}[Bilinear form]\label{def:2.1}
\emph{We define the bilinear form associated with $\mathcal{L}_\varepsilon$ as
\begin{equation}\label{pde:2.3}
 \mathrm{B}_\varepsilon[u,\phi]
 = \int_\Omega \Big\{a_{ij,\varepsilon}^{\alpha\beta}\frac{\partial u^\beta}{\partial x_j}
 + V_{i,\varepsilon}^{\alpha\beta}u^\beta\Big\}\frac{\partial \phi^\alpha}{\partial x_i} dx
 + \int_\Omega \Big\{B_{i,\varepsilon}^{\alpha\beta}\frac{\partial u^\beta}{\partial x_i}
  + c_\varepsilon^{\alpha\beta} u^\beta + \lambda u^\alpha \Big\}\phi^\alpha dx
\end{equation}
for any $u,\phi\in H^1(\Omega;\mathbb{R}^m)$.}
\end{definition}

Consider the following Dirichlet and Neumann boundary value problems
\begin{equation}\label{pde:1.1}
 (\mathbf{D}) \left\{\begin{aligned}
  \mathcal{L}_\varepsilon(u_\varepsilon) & =  -\text{div}(f) + F  & \quad &\text{in}~~\Omega,\\
  u_\varepsilon& = g &\quad &\text{on} ~\partial\Omega,
 \end{aligned}\right.
 \qquad (\mathbf{N})
 \left\{\begin{aligned}
  \mathcal{L}_\varepsilon(u_\varepsilon) & =  \text{div}(f) + F  & \quad &\text{in}~~\Omega,\\
  \mathcal{B}_\varepsilon(u_\varepsilon)& = h - n\cdot f &\quad &\text{on} ~\partial\Omega.
 \end{aligned}\right.
\end{equation}

\begin{definition}[Weak solution]\label{def:2.2}
\emph{Let $f=(f_i^\alpha)\in L^2(\Omega;\mathbb{R}^{md})$, $F\in L^{\frac{2d}{d+2}}(\Omega;\mathbb{R}^m)$,
$g=0$ and $h\in H^{-1/2,2}(\partial\Omega;\mathbb{R}^m)$ in $\eqref{pde:1.1}$. Then
\begin{itemize}
 \item[(i)] we say $u_\varepsilon\in H^1_0(\Omega;\mathbb{R}^m)$ is a weak solution to $(\mathbf{D})$, if $u_\varepsilon$ satisfies
\begin{eqnarray*}
 \mathrm{B}_\varepsilon[u_\varepsilon,\phi] = -\int_\Omega f^\alpha\cdot\nabla\phi^\alpha dx + \int_\Omega F^\alpha\phi^{\alpha} dx
 \qquad\forall~\phi\in H^1_0(\Omega;\mathbb{R}^m);
\end{eqnarray*}
 \item[(ii)] we say $u_\varepsilon\in H^1(\Omega;\mathbb{R}^m)$ is a weak solution to $(\mathbf{N})$, if $u_\varepsilon$ satisfies
\begin{eqnarray*}
 \mathrm{B}_\varepsilon[u_\varepsilon,\phi] = -\int_\Omega f^\alpha\cdot\nabla\phi^\alpha dx + \int_\Omega F^\alpha\phi^{\alpha} dx
 + <h,\phi>_{H^{-1/2,2}(\partial\Omega)\times H^{1/2,2}(\partial\Omega)}
 \quad\forall~\phi\in H^1(\Omega;\mathbb{R}^m).
\end{eqnarray*}
\end{itemize}}
\end{definition}

\begin{remark}\label{re:2.2}
\emph{Choose $\phi^\alpha = 1$ in (ii) of Definition $\ref{def:2.1}$, and then we have the compatibility condition
\begin{equation}\label{C:2.1}
 \int_\Omega \Big(B_i^{\alpha\beta}(x/\varepsilon)\frac{\partial u_\varepsilon^\beta}{\partial x_i}
 + c^{\alpha\beta}(x/\varepsilon) u_\varepsilon^\beta \Big)dx + \lambda\int_\Omega u_\varepsilon^\alpha dx
 =\int_\Omega F^\alpha dx +<h^\alpha,1>_{H^{-1/2,2}(\partial\Omega)\times H^{1/2,2}(\partial\Omega)}
\end{equation}
for $\alpha = 1,\ldots,m$, which implies the counterpart of $\eqref{C:2.1}$ in \cite{SZW2} since $B=0, c=0$ and $\lambda=0$ there.}
\end{remark}

\begin{lemma}\label{lemma:2.1}
Suppose that the coefficients of $\mathcal{L}_\varepsilon$ satisfy $\eqref{a:1}$ and $\eqref{a:3}$.
Then the corresponding bilinear form $\mathrm{B}_\varepsilon[\cdot,\cdot]$ satisfies
\begin{itemize}
\item boundedness property
\begin{equation}\label{pri:2.1}
 \big|\mathrm{B}_\varepsilon[u,v]\big| \leq C \|u\|_{H^1(\Omega)}\|v\|_{H^1(\Omega)}
 \qquad \forall~u,v\in H^1(\Omega;\mathbb{R}^m);
\end{equation}
\item coercive property
\begin{equation}\label{pri:2.2}
 \mathbf{c} \| u\|_{H^1(\Omega)}^2 \leq \mathrm{B}_\varepsilon[u,u]  \qquad \forall~u\in H^1(\Omega;\mathbb{R}^m),
\end{equation}
whenever $\lambda\geq\lambda_0$, and $\lambda_0 = \lambda_0(\mu,\kappa,m,d)$.
\end{itemize}
Here, $C$ is dependent on $\mu,\kappa,\lambda,m,d,\Omega$, while $\mathbf{c}$ depends only on $\mu,\kappa,m,d$.
\end{lemma}

\begin{thm}[Dirichlet problem]\label{thm:2.1}
The coefficients of $\mathcal{L}_\varepsilon$ and $\lambda_0$ are given as in Lemma $\ref{lemma:2.1}$.
Suppose $f\in L^2(\mathbb{R}^{md})$,
$F\in L^{\frac{2d}{d+2}}(\Omega;\mathbb{R}^m)$ and $g\in H^{1/2}(\partial\Omega;\mathbb{R}^m)$. Then
$(\mathbf{D})$
has a unique weak solution $u_\varepsilon \in H^1(\Omega;\mathbb{R}^m)$, whenever $\lambda\geq\lambda_0$,
and the solution satisfies the uniform estimate
\begin{equation}\label{pri:2.3}
\|u_\varepsilon\|_{H^1(\Omega)} \leq C \big\{\|f\|_{L^2(\Omega)}+\|F\|_{L^{\frac{2d}{d+2}}(\Omega)}+ \|g\|_{H^{1/2}(\partial\Omega)}\big\},
\end{equation}
where $C$ depends only on $\mu,\kappa,m,d$ and $\Omega$. Moreover, with one more the periodicity condition $\eqref{a:2}$ on the coefficients of
$\mathcal{L}_\varepsilon$,
we then have $u_\varepsilon\rightharpoonup u_0$ weakly in $H^1(\Omega;\mathbb{R}^m)$
and strongly in $L^2(\Omega;\mathbb{R}^m)$ as $\varepsilon\to 0$,
where $u_0$ is the weak solution to the homogenized problem
$\mathcal{L}_0(u_0) = F$ in $\Omega$ and $u_0 = g$ on $\partial\Omega$.
\end{thm}

\begin{thm}[Neumann problem]\label{thm:2.2}
The coefficients of $\mathcal{L}_\varepsilon$ and $\lambda_0$ are shown as in Lemma $\ref{lemma:2.1}$.
Then for any $f\in L^{2}(\Omega;\mathbb{R}^{md})$, $F\in L^{\frac{2d}{d+2}}(\Omega;\mathbb{R}^m)$ and
$h\in H^{-1/2}(\partial\Omega;\mathbb{R}^m)$, there exists a unique weak solution
$u_\varepsilon\in H^1(\Omega;\mathbb{R}^m)$ to $(\mathbf{N})$, whenever $\lambda\geq\lambda_0$.
Furthermore, the solution satisfies the uniform estimate
\begin{equation}\label{pri:2.4}
\|u_\varepsilon\|_{H^1(\Omega)} \leq C\big\{\|f\|_{L^{2}(\Omega)}+\|F\|_{L^{\frac{2d}{d+2}}(\Omega)}
+\|g\|_{H^{-1/2}(\partial\Omega)}\big\},
\end{equation}
where $C$ depends only on $\mu,m,d$ and $\Omega$. If we assume that
$\mathcal{L}_\varepsilon$ additionally satisfies $\eqref{a:2}$, then
the flux converges:
$A_\varepsilon\nabla u_\varepsilon + V_\varepsilon u_\varepsilon \rightharpoonup \widehat{A}\nabla u_0 + \widehat{V}u_0$, and
$B_\varepsilon\nabla u_\varepsilon + (c_\varepsilon+\lambda I)u_\varepsilon \rightharpoonup
\widehat{B}\nabla u_0 + (\widehat{c}+\lambda I)u_0$
weakly in $L^2(\Omega;\mathbb{R}^m)$ as $\varepsilon\to 0$, where $u_0$ satisfies
the corresponding homogenized problem: $\mathcal{L}_0(u_0) = \emph{div}(f)+F$ in $\Omega$ and $\mathcal{B}_0(u_0) = h$ on $\partial\Omega$
with $\mathcal{B}_0 = n\cdot\widehat{V} + n\cdot\widehat{A}\nabla$.
\end{thm}

\begin{remark}
\emph{Theorems $\ref{thm:2.1},\ref{thm:2.2}$ are referred to as the corresponding homogenization theorems.
The reader may find the related proofs in \cite{ABJLGP,VSO}.
Lemma $\ref{lemma:2.1}$ gives the uniqueness and existence of the weak solution to $(\mathbf{D})$ or
$(\mathbf{N})$. We also refer the reader to \cite{QXS,QXS1} and the references therein for more details.}
\end{remark}

\begin{definition}\label{def:2.3}
\emph{Define the adjoint operator $\mathcal{L}_\varepsilon^*$ as
\begin{equation*}
 \mathcal{L}_\varepsilon^*
 = -\text{div}\Big\{A^*(x/\varepsilon)\nabla + B^{*}(x/\varepsilon)\Big\}
 + V^{*}(x/\varepsilon)\nabla
 + c^{*}(x/\varepsilon) + \lambda I,
\end{equation*}
while the corresponding boundary operator becomes
\begin{equation*}
 \mathcal{B}_\varepsilon^* = n\cdot \big[A^{*}(x/\varepsilon)\nabla+B^*(x/\varepsilon)\big].
\end{equation*}
Furthermore, the related bilinear form is given by
\begin{equation*}
 \mathrm{B}_\varepsilon^*[v,\phi]
 = \int_\Omega \Big\{a_{ij,\varepsilon}^{\alpha\beta}\frac{\partial v^\alpha}{\partial x_i}
 + B_{j,\varepsilon}^{\alpha\beta} v^\alpha\Big\} \frac{\partial \phi^\beta}{\partial x_j} dx
 + \int_\Omega \Big\{V_{i,\varepsilon}^{\alpha\beta}\frac{\partial v^\beta}{\partial x_i}
  + c_\varepsilon^{\alpha\beta} v^\alpha + \lambda v^\beta \Big\} \phi^\beta dx
\end{equation*}
for any $v,\phi\in H^1(\Omega;\mathbb{R}^m)$.}
\end{definition}

\begin{remark}\label{re:2.3}
\emph{Let $u_\varepsilon,v_\varepsilon\in H_0^1(\Omega;\mathbb{R}^m)$ be the two weak solutions to the Dirichlet problems
\begin{equation*}
\left\{\begin{aligned}
 \mathcal{L}_\varepsilon(u_\varepsilon) & = \text{div}(f)+F &\quad \text{in}~~~\Omega,\\
 u_\varepsilon & = 0     &\quad \text{on}~\partial\Omega,
\end{aligned}\right.
\qquad\quad
\left\{\begin{aligned}
\mathcal{L}_\varepsilon^*(v_\varepsilon) &= \text{div}(\phi) + \Phi &\quad \text{in}~~~\Omega,\\
 v_\varepsilon &= 0 &\quad \text{on}~\partial\Omega,
\end{aligned}\right.
\end{equation*}
respectively. It follows from $\mathrm{B}_\varepsilon[u_\varepsilon,v_\varepsilon] = \mathrm{B}_\varepsilon^*[v_\varepsilon,u_\varepsilon]$
that
\begin{equation}\label{eq:2.1}
 <\mathcal{L}_\varepsilon(u_\varepsilon),v_\varepsilon>
 = <u_\varepsilon,\mathcal{L}_\varepsilon^*(v_\varepsilon)>.
\end{equation}}
\end{remark}

\begin{remark}\label{re:2.4}
\emph{ If $u_\varepsilon,v_\varepsilon\in H^1(\Omega;\mathbb{R}^m)$ are two weak solutions to the Neumann problems
\begin{equation*}
\left\{\begin{aligned}
 \mathcal{L}_\varepsilon(u_\varepsilon) & = \text{div}(f)+F &\quad \text{in}~~~\Omega,\\
 \mathcal{B}_\varepsilon(u_\varepsilon) & = h-n\cdot f      &\quad \text{on}~\partial\Omega,
\end{aligned}\right.
\qquad\quad
\left\{\begin{aligned}
\mathcal{L}_\varepsilon^*(v_\varepsilon) &= \text{div}(\phi) + \Phi &\quad \text{in}~~~\Omega,\\
\mathcal{B}_\varepsilon^*(v_\varepsilon) &= \eta -n\cdot \phi &\quad \text{on}~\partial\Omega,
\end{aligned}\right.
\end{equation*}
respectively, then we have the second Green's formula
\begin{equation}\label{eq:2.2}
<\mathcal{L}_\varepsilon(u_\varepsilon),v_\varepsilon> - <u_\varepsilon,\mathcal{L}_\varepsilon^*(v_\varepsilon)>
= -<\mathcal{B}_\varepsilon(u_\varepsilon),v_\varepsilon> + <u_\varepsilon,\mathcal{B}_\varepsilon^*(v_\varepsilon)>
\end{equation}
by noting that $\mathrm{B}_\varepsilon[u_\varepsilon,v_\varepsilon] = \mathrm{B}_\varepsilon^*[v_\varepsilon,u_\varepsilon]$.}
\end{remark}

\begin{remark}\label{re:2.5}
\emph{To handle the convergence rates, we define some auxiliary functions via
\begin{equation}\label{eq:2.4}
b_{ik}^{\alpha\gamma}(y) = \hat{a}_{ik}^{\alpha\gamma} - a_{ik}^{\alpha\gamma}(y)
- a_{ij}^{\alpha\beta}(y)\frac{\partial}{\partial y_j}\big\{\chi_{k}^{\beta\gamma}\big\}, \qquad
b_{i0}^{\alpha\gamma}(y) =
\hat{V}_i^{\alpha\gamma} - V_{i}^{\alpha\gamma}(y)-a_{ij}^{\alpha\beta}(y)\frac{\partial}{\partial y_j}\big\{\chi_0^{\beta\gamma}\big\},
\end{equation}
and
\begin{equation}\label{eq:2.5}
\left.\begin{aligned}
\Delta\vartheta_{i}^{\alpha\gamma}
 & = W_i^{\alpha\gamma} := \hat{B}_i^{\alpha\gamma} - B_{i}^{\alpha\gamma}(y)
- B_{j}^{\alpha\beta}(y)\frac{\partial}{\partial y_j}\big\{\chi_i^{\beta\gamma}\big\} & \quad \text{in} ~~\mathbb{R}^d,
&\qquad \int_Y \vartheta_i^{\alpha\beta}(y) dy = 0 ,\\
\Delta \vartheta_0^{\alpha\gamma}
&= W_0^{\alpha\gamma} := \hat{c}^{\alpha\gamma} - c^{\alpha\gamma}(y)
- B_{i}^{\alpha\beta}(y)\frac{\partial}{\partial y_i}\big\{\chi_0^{\beta\gamma}\big\} & \quad \text{in}~~\mathbb{R}^d ,
&\qquad \int_Y \vartheta_0^{\alpha\beta}(y) dy = 0 .
\end{aligned}\right.
\end{equation}
We mention that the existence of $\vartheta_{k}$ is given by \cite[Theorem 4.28]{ACPD}
on account of $\dashint_{Y}\vartheta_k^{\alpha\gamma}(y)dy = 0$ for $k=0,1,\ldots,d$. Furthermore
it is not hard to see that $\vartheta_k^{\alpha\gamma}$ is periodic and belongs to $H^2_{loc}(\mathbb{R}^d)$.}
\end{remark}

\begin{lemma}\label{lemma:2.2}
There exist $E_{jik}^{\alpha\gamma}\in H^1_{per}(Y)$ with $k = 0,1,\ldots,d$, such that
\begin{equation}\label{eq:2.6}
 b_{ik}^{\alpha\gamma} = \frac{\partial}{\partial y_j}\big\{E_{jik}^{\alpha\gamma}\big\}
 \qquad\text{and}\qquad
 E_{jik}^{\alpha\gamma} = - E_{ijk}^{\alpha\gamma},
\end{equation}
where $1\leq i,j\leq d$ and $1\leq\alpha,\gamma\leq m$. Moreover if $\chi_{k}$ is H\"older continuous,
then $E_{jik}^{\alpha\gamma}\in L^\infty(Y)$.
\end{lemma}

\begin{pf}
See \cite[Lemma 2.8]{QXS1}. \qed
\end{pf}

\begin{definition}\label{def:2.4}
\emph{Fix $\zeta\in C_0^\infty(B(0,1/2))$,  and $\int_{\mathbb{R}^d}\zeta = 1$. Define the smoothing operator
\begin{equation}\label{eq:2.7}
S_\varepsilon(f)(x) = f*\zeta_\varepsilon(x) = \int_{\mathbb{R}^d} f(x-y)\zeta_\varepsilon(y)dy,
\end{equation}
where $\zeta_\varepsilon(x) = \varepsilon^{-d}\zeta(x/\varepsilon)$.}
\end{definition}

\begin{lemma}\label{lemma:2.3}
Let $f\in L^p(\mathbb{R}^d)$ for some $1\leq p<\infty$. Then for any $h\in L^p_{per}(\mathbb{R}^d)$,
\begin{equation}\label{pri:2.5}
\|h(\cdot/\varepsilon)S_\varepsilon(f)\|_{L^p(\mathbb{R}^d)}\leq C \|h\|_{L^p(Y)}\|f\|_{L^p(\mathbb{R}^d)},
\end{equation}
where $C$ depends only on $d$.
\end{lemma}

\begin{pf}
See \cite[Lemma 2.1]{SZW12}.
\end{pf}

\begin{lemma}\label{lemma:2.4}
Let $f\in W^{1,p}(\mathbb{R}^d)$ for some $1<p<\infty$. Then we have
\begin{equation}\label{pri:2.6}
\|S_\varepsilon(f) - f\|_{L^p(\mathbb{R}^d)} \leq C\varepsilon\|\nabla f\|_{L^p(\mathbb{R}^d)},
\end{equation}
and furthermore obtain
\begin{equation}\label{pri:2.7}
\|S_\varepsilon(f)\|_{L^2(\mathbb{R}^d)}\leq C\varepsilon^{-1/2}\|f\|_{L^q(\mathbb{R}^d)}
\qquad \text{and}\qquad
\|S_\varepsilon(f)-f\|_{L^2(\mathbb{R}^d)}\leq C\varepsilon^{1/2}\|\nabla f\|_{L^q(\mathbb{R}^d)},
\end{equation}
where $q = \frac{2d}{d+1}$, and $C$ depends only on $d$.
\end{lemma}

\begin{pf}
See \cite[Lemma 2.2]{SZW12}.
\end{pf}


\begin{remark}\label{re:2.8}
\emph{Throughout the paper, let $B(P,r)$ denote the open ball centered at $P$ of radius r,
and the symbol $r_0$ only represents the diameter of $\Omega$.
We say $\partial\Omega\in C^{0,1}$ (Lipschitz), if there exists $R$ such that
for each point $P\in\partial\Omega$ there is a new coordinate system in $\mathbb{R}^d$ obtained
from the standard Euclidean coordinate system translation and rotation so that
$P=(0,0)$ and
\begin{equation*}
 B(P,R)\cap \Omega
 = B(P,R)\cap\big\{(x^\prime,x_d)\in\mathbb{R}^d:x^\prime\in\mathbb{R}^{d-1}
 ~\text{and}~ x_d>\phi(x^\prime)\big\},
\end{equation*}
where $\phi\in C^{0,1}(\mathbb{R}^{d-1})$ is a boundary function with $\phi(0) = 0$ and
$\|\nabla\phi\|_{L^\infty(\mathbb{R}^d)}\leq M_0$. Note that $M_0$ indicates the boundary character of $\Omega$.
In the paper, saying a constant $C$ depends on $\Omega$ means this constant involves both $M_0$ and $|\Omega|$,
where $|\Omega|$ denotes the volume of $\Omega$.}
\end{remark}

\begin{remark}\label{re:2.10}
\emph{For $0\leq r< c_0$,
we may assume that there exist homeomorphisms $\Lambda_r: \partial\Omega\to \partial\Sigma_r = S_r$ such that $\Lambda_0(Q) = Q$,
$|\Lambda_r(Q) - \Lambda_t(P)| \sim |r-t| + |Q-P|$ and
$|\Lambda_r(Q) - \Lambda_t(Q)|\leq C\text{dist}(\Lambda_r(Q),S_t)$ for any
$r>s$ and $P,Q\in\partial\Omega$ (which are bi-Lipschitz maps, see \cite[pp.1014]{SZW2}). Especially, we may have
$\max_{r\in[0,c_0]}\{\|\nabla\Lambda_r\|_{L^\infty(\partial\Omega)},\|\nabla(\Lambda_r^{-1})\|_{L^\infty(\partial\Omega)}\} \leq C(M_0)$.
For a function $h$,
we define the radial maximal function
$\mathcal{M}(h)$ on $\partial\Omega$ as
\begin{equation}\label{def:2.6}
 \mathcal{M}(h)(Q) = \sup\big\{|h(\Lambda_r(Q))|: 0\leq r\leq  c_0\big\} \quad\qquad \forall ~Q\in\partial\Omega.
\end{equation}
We mention that the radial maximal function will play an important role in
the study of convergence rates for Lipschitz domains
(we refer the reader to\cite{SZW2} for the original idea,
and we also refer the reader to \cite[Theorem 5.1]{CTM} for the existence of such bi-Lipschitz maps).}
\end{remark}

\begin{definition}\label{def:2.7}
\emph{The non-tangential maximal function of $u$ is defined by
\begin{equation}\label{def:2.8}
(u)^*(Q) = \sup\big\{ |u(x)|:x\in \Gamma_{N_0}(Q)\big\} \qquad\quad \forall~ Q\in\partial\Omega,
\end{equation}
where $\Gamma_{N_0}(Q) = \{x\in\Omega:|x-Q|\leq N_0\delta(x)\}$ is the cone with vertex $Q$ and aperture $N_0$,
and $N_0>1$ is sufficiently large.}
\end{definition}

\begin{remark}\label{re:2.11}
\emph{Let $h\in L^p(\Omega)$ with $1\leq p<\infty$. For any $r\in(0,c_0)$,
$\Lambda_r$ is given in Remark $\ref{re:2.10}$, and
we can show the estimate of $\|h\|_{L^p(\Omega\setminus\Sigma_r)}$. By $\eqref{def:2.8}$, we note that
$h(\Lambda_r(x))\leq \mathcal{M}(h)(x)$ a.e. $x\in\partial\Omega$ for all $r\in (0,c_0)$.
Then
\begin{equation}\label{pri:2.8}
\begin{aligned}
 \int_{\Omega\setminus\Sigma_r} |h|^p dx
 &= \int_0^r\int_{S_t=\Lambda_{t}(\partial\Omega)} |h(y)|^p dS_t(y)dt \\
 & = \int_0^r\int_{\partial\Omega} |h(\Lambda_t(z))|^p |\nabla\Lambda_{t}|dS(z)dt \leq Cr\int_{\partial\Omega} |\mathcal{M}(h)|^p dS
 \leq  Cr\int_{\partial\Omega} |(h)^{*}|^p dS,
\end{aligned}
\end{equation}
where $C$ depends only on $p$ and the boundary character. 
We note that the first equality is based on the so-called co-area formula $\eqref{eq:2.9}$,
and we use the change of variable in the second one. Besides,
the first inequality follows from Remark $\ref{re:2.10}$. In the last one,
it is not hard to see $\mathcal{M}(h)(Q)\leq (h)^*(Q)$ by comparing Definition $\ref{def:2.7}$ with $\eqref{def:2.6}$.}

\emph{We now explain the co-area formula used here.
Let $Z(0;r)=\{x\in\Omega:0<\delta(x)\leq r\}$,
then $Z(0;r) = \Omega\setminus\Sigma_r$.
Here we point out $|\nabla\delta(x)| =1$ a.e. $x\in\Omega$ without the proof (see \cite[pp.142]{LCE1}).
In view of co-area formula (see \cite[Theorem 3.13]{LCE1}), we have
\begin{equation}\label{eq:2.9}
 \int_{\Omega\setminus\Sigma_r} |h|^p dx = \int_{Z(0;r)} |h|^p dx
 = \int_{0}^{r}\int_{\{x\in\Omega:\delta(x)=t\}}\frac{|h|^p}{|\nabla\delta|}d\mathcal{H}^{d-1}dt
 =\int_{0}^{r}\int_{S_t}|h|^p dS_tdt,
\end{equation}
where $S_r=\{x\in\Omega:\delta(x) = t\}$, $d\mathcal{H}^{d-1}$ is the ($d-1$)-dimensional Hausdorff measure,
and $dS_t=d\mathcal{H}^{d-1}(S_t)$ denotes the surface measure of $S_t$.}
\end{remark}

\begin{lemma}\label{lemma:2.5}
Let $\Omega$ be a Lipschitz domain, and $\mathcal{M}$ associated with $c_0$ is defined in Remark $\ref{re:2.10}$.
Then for any $h\in H^1(\Omega)$, we have the following estimate
\begin{equation}\label{pri:2.9}
 \|\mathcal{M}(h)\|_{L^2(\partial\Omega)}
 \leq C\|h\|_{H^1(\Omega\setminus\Sigma_{c_0})},
\end{equation}
where $C$ depends only on $d,c_0$ and the character of $\Omega$.
\end{lemma}

\begin{pf}
See \cite[Lemma 2.24]{QXS1}.
\end{pf}

\begin{lemma}\label{lemma:2.6}
 Suppose that $u_\varepsilon,u_0\in H^1(\Omega;\mathbb{R}^m)$ satisfy $\mathcal{L}_\varepsilon(u_\varepsilon) = \mathcal{L}_0(u_0)$ in
 $\Omega$. Let
 \begin{equation}\label{eq:2.12}
  w_\varepsilon^\beta = u_\varepsilon^\beta - u_0^\beta - \varepsilon\sum_{k=0}^d\chi_{k}^{\beta\gamma}(x/\varepsilon)\varphi_k^\gamma,
 \end{equation}
 where $\varphi=(\varphi_k^\gamma)\in H^1(\mathbb{R}^d;\mathbb{R}^{m(d+1)})$. Then
 \begin{itemize}
  \item[\emph{(i)}] if $u_\varepsilon = u_0$ on $\partial\Omega$, we have
 \begin{equation}\label{pde:2.4}
  \mathcal{L}_\varepsilon(w_\varepsilon) = -\emph{div}(\tilde{f}) + \tilde{F} \quad \text{in}~\Omega,
  \qquad w_\varepsilon = \varepsilon \sum_{k=0}^d\chi_{k,\varepsilon}\varphi_k \quad\text{on}~\partial\Omega;
\end{equation}
 \item [\emph{(ii)}] if $\mathcal{B}_\varepsilon(u_\varepsilon) = \mathcal{B}_0(u_0)$ on $\partial\Omega$, we have
  \begin{equation}\label{pde:2.5}
  \mathcal{L}_\varepsilon(w_\varepsilon) = -\emph{div}(\tilde{f}) + \tilde{F} \quad \text{in}~\Omega,
  \qquad \mathcal{B}_\varepsilon(w_\varepsilon) = n\cdot\tilde{f} + \varepsilon n\cdot\mathcal{J}\quad \text{on}~\partial\Omega,
\end{equation}
\end{itemize}
where $n$ is the outward unit normal vector to $\partial\Omega$.
Note that $\tilde{f} = (\tilde{f}^\alpha_i)$ and $\tilde{F} = (\tilde{F}^\alpha)$ are given in $\eqref{eq:2.10}$,
and $\mathcal{J} = (\mathcal{J}_i^\alpha)$ is shown in $\eqref{eq:2.11}$.
\end{lemma}

\begin{pf}
 See \cite[Lemma 5.1]{QXS1}.
 \qed
\end{pf}

\begin{remark}\label{re:2.12}
\emph{For simplicity of presentation in Lemma $\ref{lemma:2.6}$, we set
\begin{equation}\label{eq:2.10}
\begin{aligned}
& \tilde{f}^\alpha_i
 = \mathcal{K}_i^\alpha
+\big[\hat{a}_{ij}^{\alpha\beta}-a_{ij,\varepsilon}^{\alpha\beta}\big]
\Big[\frac{\partial u_0^\beta}{\partial x_j}-\varphi_j^\beta\Big]
+\big[\hat{V}_i^{\alpha\beta}-V_{i,\varepsilon}^{\alpha\beta}\big]\big[u_0^\beta
-\varphi_0^\beta\big] - \varepsilon\big(\mathcal{I}_{i}^\alpha+\mathcal{J}_i^\alpha\big) \\
& \tilde{F}^\alpha
= [\hat{B}_i^{\alpha\beta}-B_{i,\varepsilon}^{\alpha\beta}]
\big[\frac{\partial u_0^\beta}{\partial x_i}-\varphi^\beta_i\big]
+\big[\hat{c}^{\alpha\beta}-c_{\varepsilon}^{\alpha\beta}\big]\big[u_0^\beta-\varphi^\beta_0\big]
-\varepsilon\big(\mathcal{M}^\alpha + \mathcal{N}^\alpha\big)
\end{aligned}
\end{equation}
where
\begin{equation}\label{eq:2.11}
\begin{aligned}
& \mathcal{I}_i^\alpha
 = a_{ij,\varepsilon}^{\alpha\beta}\sum_{k=0}^d\chi_{k,\varepsilon}^{\beta\gamma}
 \frac{\partial}{\partial x_j}\big\{\varphi^\gamma_k\big\}
+V_{i,\varepsilon}^{\alpha\beta}\sum_{k=0}^d\chi_{k,\varepsilon}^{\beta\gamma}
\varphi^\gamma_k,
\qquad
\mathcal{J}_{i}^\alpha = \sum_{k=0}^d\Big(\frac{\partial \vartheta_k^{\alpha\gamma}}{\partial y_i}\Big)_\varepsilon\varphi^\gamma_k,
\qquad
\mathcal{K}_{i}^\alpha =  \sum_{j=0}^db_{ij,\varepsilon}^{\alpha\gamma}\varphi^\gamma_j, \\
& \mathcal{M}^\alpha =
\sum_{k=0}^d\Big[\Big(\frac{\partial \vartheta_k^{\alpha\gamma}}{\partial y_i}\Big)_\varepsilon
+ B_{i,\varepsilon}^{\alpha\beta}\chi_{k,\varepsilon}^{\beta\gamma}\Big]
\frac{\partial}{\partial x_i}\big\{\varphi_k^\gamma\big\},
\qquad
\mathcal{N}^\alpha =
\big[c_{\varepsilon}^{\alpha\beta}+\lambda e^{\alpha\beta}\big]
 \sum_{k=0}^d\chi_{k,\varepsilon}^{\beta\gamma}\varphi^\gamma_k,
 \qquad y =x/\varepsilon.
\end{aligned}
\end{equation} }
\end{remark}

\begin{lemma}\label{lemma:3.1}
Suppose that the coefficients of $\mathcal{L}_\varepsilon$ satisfy $\eqref{a:1}$, $\eqref{a:2}$ and $\eqref{a:3}$.
Let $w_\varepsilon$ be defined by $\eqref{eq:2.12}$, and $u_\varepsilon,u_0$ be
the weak solutions to $\eqref{pde:1.2}$.
Then we have
\begin{equation}\label{pri:3.1}
\begin{aligned}
 \|w_\varepsilon\|_{H^1(\Omega)}
 \leq C\Big\{\|\varpi_\varepsilon\vec{\phi}\|_{L^2(\Omega\setminus\Sigma_{2\varepsilon})}
 &+\|\nabla u_0 - \nabla\vec{\varphi}\|_{L^2(\Omega)} \\
 &+\|u_0 - \varphi_0\|_{L^2(\Omega)}
 +\varepsilon\|\varpi_\varepsilon\nabla\vec{\phi}\|_{L^2(\Omega)}
 +\varepsilon\|\varpi_\varepsilon\vec{\phi}\|_{L^2(\Omega)} \Big\},
\end{aligned}
\end{equation}
where $\vec{\varphi}=(\varphi_1,\cdots,\varphi_d)$, $\vec{\phi}=(\varphi_0,\vec{\varphi})$,
and $\varpi_\varepsilon$ denotes the periodic functions
(partially or fully) depending on the coefficients of $\mathcal{L}_\varepsilon$,
the correctors $\{\chi_k\}_{k=0}^d$,
and auxiliary functions $\{b_{ik},E_{jik},\nabla\vartheta_k\}_{k=0}^d$.
\end{lemma}

\begin{remark}\label{re:2.6}
\emph{As we proceed to prove the above lemma and others,
we are often confronted with the periodic functions such as the coefficients of $\mathcal{L}_\varepsilon$,
the correctors $\{\chi_k\}_{k=0}^d$, as well as some auxiliary functions $\{b_{ik},E_{jik},\nabla\vartheta_k\}_{k=0}^d$ in the calculations.
These periodic functions are actually the known quantities.
The algebra combination of these periodic functions
is always lengthy to write, so we denote it by $\varpi_\varepsilon$ for short. If
ignoring the form of the different combinations,
then the notation $\varpi_\varepsilon$ will play a similar role as the constant $C$ does in the estimates,
which becomes an universal periodic function determined by $\mu,\kappa,m,d$. }
\end{remark}

\begin{pf}
In view of (i) in Lemma $\ref{lemma:2.6}$, we have
\begin{equation*}
 \mathcal{L}_\varepsilon(w_\varepsilon) = \text{div}(\tilde{f}) + \tilde{F} \quad \text{in}~\Omega,
 \qquad w_\varepsilon = \varepsilon\sum_{k=0}^d\chi_{k,\varepsilon}\varphi_k\quad \text{on}~\partial\Omega,
\end{equation*}
where $\tilde{f}$ and $\tilde{F}$ are given in Remark $\ref{re:2.12}$. It is reasonable to consider dividing
$w_\varepsilon$ into $w_{\varepsilon,1}$ and $w_{\varepsilon,2}$, and they satisfy
\begin{equation*}
(1)\left\{
\begin{aligned}
\mathcal{L}_\varepsilon(w_{\varepsilon,1}) &= -\text{div}(\tilde{f} - \mathcal{K})+\tilde{F} &\quad &\text{in}~~\Omega, \\
 w_{\varepsilon,1} &= \varepsilon\sum_{k=0}^d\chi_{k,\varepsilon}\varphi_k &\quad&\text{on} ~\partial\Omega,
\end{aligned}\right.
\qquad
(2)\left\{
\begin{aligned}
\mathcal{L}_\varepsilon(w_{\varepsilon,2}) &= -\text{div}(\mathcal{K}) &\quad &\text{in}~~\Omega, \\
 w_{\varepsilon, 2} &= 0 &\quad&\text{on} ~\partial\Omega,
\end{aligned}\right.
\end{equation*}
respectively.

For (1), it follows from $\eqref{pri:2.3}$ that
\begin{equation}\label{f:3.1}
\begin{aligned}
\|w_{\varepsilon,1}\|_{H^1(\Omega)}
\leq C\Big\{\|\nabla u_0 - \nabla\vec{\varphi}\|_{L^2(\Omega)}
 &+\|u_0 - \varphi_0\|_{L^2(\Omega)} \\
 &+\varepsilon\|\varpi_\varepsilon\nabla\vec{\phi}\|_{L^2(\Omega)}
 +\varepsilon\|\varpi_\varepsilon\vec{\phi}\|_{L^2(\Omega)}\Big\}
+ \varepsilon\|\chi_{k,\varepsilon}\varphi_k\|_{H^{1/2}(\partial\Omega)}
\end{aligned}
\end{equation}
with summation convention applied to $k$ from $0$ to $d$.
We now handle the term of $\|\chi_{k,\varepsilon}\varphi_k\|_{H^{1/2}(\partial\Omega)}$, and then
\begin{equation}\label{f:3.2}
\begin{aligned}
\|\chi_{k,\varepsilon}\varphi_k\|_{H^{1/2}(\partial\Omega)}
&\leq C\|(1-\psi_\varepsilon)\chi_{k,\varepsilon}\varphi_k\|_{H^{1}(\Omega)} \\
&\leq C\|(1-\psi_\varepsilon)\chi_{k,\varepsilon}\varphi_k\|_{L^{2}(\Omega)}
+ C\|\nabla[(1-\psi_\varepsilon)\chi_{k,\varepsilon}\varphi_k]\|_{L^{2}(\Omega)} \\
&\leq C\|\varpi_\varepsilon\vec{\phi}\|_{L^2(\Omega\setminus\Sigma_{2\varepsilon})}
+ C\varepsilon^{-1}\|\varpi_\varepsilon\vec{\phi}\|_{L^2(\Omega\setminus\Sigma_{2\varepsilon})}
+ C\|\varpi_\varepsilon\nabla\vec{\phi}\|_{L^2(\Omega)},
\end{aligned}
\end{equation}
where $\varpi_\varepsilon$ is explained in Remark $\ref{re:2.6}$, here depending on $\chi_{k,\varepsilon}$ or $\nabla\chi_{k,\varepsilon}$,
and $\vec{\phi} = (\varphi_0,\cdots,\varphi_d)$.
Hence, plugging $\eqref{f:3.2}$ back into $\eqref{f:3.1}$, we obtain
\begin{equation}\label{f:3.6}
\begin{aligned}
\|w_{\varepsilon,1}\|_{H^1(\Omega)}
\leq C\Big\{\|\varpi_\varepsilon\vec{\phi}\|_{L^2(\Omega\setminus\Sigma_{2\varepsilon})}
&+ \|\nabla u_0 - \nabla\vec{\varphi}\|_{L^2(\Omega)} \\
&+\|u_0 - \varphi_0\|_{L^2(\Omega)}
 +\varepsilon\|\varpi_\varepsilon\nabla\vec{\phi}\|_{L^2(\Omega)}
 +\varepsilon\|\varpi_\varepsilon\vec{\phi}\|_{L^2(\Omega)}\Big\}.
\end{aligned}
\end{equation}

For (2), in view of (i) in Definition $\ref{def:2.2}$, we have
\begin{equation}\label{f:3.3}
\mathrm{B}_\varepsilon[w_{\varepsilon,2},v] = \int_\Omega \mathcal{K}\cdot\nabla v dx =: R(v)
\end{equation}
for any $v\in H^1_0(\Omega;\mathbb{R}^m)$.
According to Lemma $\ref{lemma:2.2}$, $R(v)$ in $\eqref{f:3.3}$ satisfies
\begin{eqnarray*}
R(v)& = & \varepsilon\int_\Omega\sum_{k=0}^d\Big\{\frac{\partial}{\partial x_j}\big[{E}_{jik,\varepsilon}^{\alpha\gamma}\big]
\varphi^\gamma_k\Big\}\frac{\partial v^\alpha}{\partial x_i}dx\\
&=&\varepsilon\int_\Omega
\frac{\partial}{\partial x_j}\Big\{\sum_{k=0}^d{E}_{jik,\varepsilon}^{\alpha\gamma}\varphi^\gamma_k\Big\}
\frac{\partial v^\alpha}{\partial x_i} dx
-\varepsilon\int_\Omega
\sum_{k=0}^d\Big\{{E}_{jik,\varepsilon}^{\alpha\gamma}
\frac{\partial}{\partial x_j}\big[\varphi^\gamma_k\big]\Big\}
\frac{\partial v^\alpha}{\partial x_i} dx\\
&=:& R_1(v) - R_2(v).
\end{eqnarray*}
Note that due to the antisymmetry of $E_{jik}$ with respect to $i,j$, we obtain
\begin{eqnarray*}
R_1(v) &=& \varepsilon\int_\Omega\frac{\partial}{\partial x_j}
\Big\{\big[\psi_\varepsilon+(1-\psi_\varepsilon)\big]{E}_{jik,\varepsilon}^{\alpha\gamma}\varphi^\gamma_k
\Big\}\frac{\partial v^\alpha}{\partial x_i} dx\\
&=&\varepsilon\int_\Omega\frac{\partial}{\partial x_j}
\Big\{(1-\psi_\varepsilon){E}_{jik,\varepsilon}^{\alpha\gamma}\varphi^\gamma_k\Big\}\frac{\partial v^\alpha}{\partial x_i} dx
-~\varepsilon\int_\Omega\psi_\varepsilon{E}_{jik,\varepsilon}^{\alpha\gamma}\varphi^\gamma_k
\frac{\partial^2 v^\alpha}{\partial x_i\partial x_j} dx \\
&=&\varepsilon\int_\Omega\frac{\partial}{\partial x_j}
\Big\{(1-\psi_\varepsilon){E}_{jik,\varepsilon}^{\alpha\gamma}\varphi^\gamma_k
\Big\}\frac{\partial v^\alpha}{\partial x_i} dx,
\end{eqnarray*}
where $\psi_\varepsilon\in C_0^\infty(\Omega)$ satisfies $\eqref{def:2.5}$, and $k=0,1,\ldots,d$. Moreover, we have
\begin{equation*}
\begin{aligned}
R_1(v) & = - \int_\Omega
\frac{\partial\psi_\varepsilon}{\partial x_j}{E}_{jik,\varepsilon}^{\alpha\gamma}\varphi^\gamma_k
\frac{\partial v^\alpha}{\partial x_i} dx
+\int_\Omega (1-\psi_\varepsilon)b_{ik,\varepsilon}^{\alpha\gamma}\varphi^\gamma_k
\frac{\partial v^\alpha}{\partial x_i} dx
 + \varepsilon \int_\Omega(1-\psi_\varepsilon){E}_{jik,\varepsilon}^{\alpha\gamma}
\frac{\partial\varphi^\gamma_k}{\partial x_j}\frac{\partial v^\alpha}{\partial x_i} dx,
\end{aligned}
\end{equation*}
and this indicates
\begin{equation}\label{f:3.4}
\begin{aligned}
\big|R_1(v)\big| & \leq \big\{\|{E}_{\varepsilon}\vec{\phi}\|_{L^2(\Omega\setminus\Sigma_{2\varepsilon})}
+\|b_{\varepsilon}\vec{\phi}\|_{L^2(\Omega\setminus\Sigma_{2\varepsilon})}
+ \varepsilon\|{E}_{\varepsilon}\nabla\vec{\phi}\|_{L^2(\Omega\setminus\Sigma_{2\varepsilon})}\big\}
\|\nabla v\|_{L^2(\Omega\setminus\Sigma_{2\varepsilon})}.
\end{aligned}
\end{equation}
Meanwhile we arrive at
\begin{equation}\label{f:3.5}
 \big|R_2(v)\big| \leq \varepsilon\|{E}_\varepsilon\nabla\vec{\phi}\|_{L^2(\Omega)}\|\nabla v\|_{L^2(\Omega)}.
\end{equation}

Let $v=w_{\varepsilon,2}$. In view of $\eqref{pri:2.2}$, $\eqref{f:3.3}$, $\eqref{f:3.4}$ and $\eqref{f:3.5}$, we have
\begin{equation}\label{f:3.7}
\|w_{\varepsilon,2}\|_{H_0^1(\Omega)}
\leq C\big\{\|{E}_{\varepsilon}\vec{\phi}\|_{L^2(\Omega\setminus\Sigma_{2\varepsilon})}
+\|b_{\varepsilon}\vec{\phi}\|_{L^2(\Omega\setminus\Sigma_{2\varepsilon})}
+\varepsilon\|{E}_{\varepsilon}\nabla\vec{\phi}\|_{L^2(\Omega)}\big\}.
\end{equation}

It is clear to see that the desired estimate $\eqref{pri:3.1}$
follows from $\eqref{f:3.6}$ and $\eqref{f:3.7}$, which completes the proof.
\qed
\end{pf}

\section{Weighted-type inequalities and duality lemmas}
The core techniques of this paper are introduced in this section.
As we mentioned before a crucial reason why we developed the weighted estimates is
that the distance function $\delta$ used to help cancel the singularity of $\nabla^2 u_0$ near the boundary of
$\Omega$. Therefore $\delta$ is chosen to be the weighted function. To achieve our goal,
it is natural to expect the weighted function $\delta$
to pass through the convolution freely in the calculations,
and Lemmas $\ref{lemma:2.9}-\ref{lemma:2.8}$ exactly realize this thinking.

\begin{lemma}\label{lemma:2.9}
Let $\delta(x)$, $\Sigma_{2\varepsilon}$ be defined in Subsection $\ref{sec:1.1}$. Then for any $x\in\Sigma_{2\varepsilon}$ we have
\begin{equation}\label{pri:2.10}
 \big|S_\varepsilon(\delta)(x)\big|\leq 2\delta(x),
 \qquad \big|S_\varepsilon(\delta^{-1})(x)\big|\leq 2\delta^{-1}(x).
\end{equation}
\end{lemma}

\begin{pf}
It is clear to see that
\begin{equation*}
\begin{aligned}
\Big|\delta(x) - \int_{\mathbb{R}^d}\zeta_\varepsilon(x-y)\delta(y)dy\Big|
&\leq \int_{\mathbb{R}^d}\zeta_\varepsilon(x-y)\big|\delta(x)-\delta(y)\big|dy \\
&\leq  \int_{\mathbb{R}^d}\zeta_\varepsilon(x-y)\|\nabla\delta\|_{L^\infty(\Omega)}|x-y|dy
\leq \varepsilon \int_{\mathbb{R}^d}\zeta_\varepsilon(x-y)dy
\leq \varepsilon,
\end{aligned}
\end{equation*}
where we use the fact of $\|\nabla\delta\|_{L^\infty(\Omega)} = 1$ (see \cite[Theorem 3.14]{LCE1}) in the third inequality.
Since $\delta(x)>\varepsilon$ whenever $x\in\Sigma_{2\varepsilon}$, we have
\begin{equation*}
\big|S_\varepsilon(\delta)(x)\big|
\leq \Big|\delta(x) - \int_{\mathbb{R}^d}\zeta_\varepsilon(x-y)\delta(y)dy\Big| + \delta(x)
\leq 2\delta(x).
\end{equation*}

By the same token, we have
\begin{equation*}
\begin{aligned}
\Big|\delta^{-1}(x) - \int_{\mathbb{R}^d}\zeta_\varepsilon(x-y)\delta^{-1}(y)dy\Big|
&\leq  \int_{\mathbb{R}^d}\zeta_\varepsilon(x-y)\|\nabla\delta\|_{L^\infty(\Omega)}|x-y|\delta^{-1}(x)\delta^{-1}(y)dy \\
&\leq \varepsilon \int_{\mathbb{R}^d}\zeta_\varepsilon(x-y)\delta^{-1}(x)\delta^{-1}(y)dy
 \leq \delta^{-1}(x)\int_{\mathbb{R}^d}\zeta_\varepsilon(x-y) dy = \delta^{-1}(x),
\end{aligned}
\end{equation*}
where we point out $y\in\Sigma_{\varepsilon}$ at most, and therefore $\delta(y)>\varepsilon$. Thus we have
$|S_\varepsilon(\delta^{-1})(x)|\leq 2\delta^{-1}(x)$, and complete the proof.
\qed

\end{pf}

\begin{lemma}\label{lemma:2.7}
 Let $f\in L^2(\Omega)$ be supported in $\Sigma_{2\varepsilon}$, and $g\in L_{per}^2(Y)$, then we have
 \begin{equation}\label{pri:2.11}
  \Big(\int_{\Sigma_{2\varepsilon}}\big|g(x/\varepsilon)S_\varepsilon(f)(x)\big|^2 \delta^{-1}(x)dx\Big)^{1/2}
  \leq C\|g\|_{L^2(Y)}\Big(\int_{\Sigma_{2\varepsilon}}\big|f(x)\big|^2\delta^{-1}(x)dx\Big)^{1/2},
 \end{equation}
 and
 \begin{equation}\label{pri:2.12}
 \Big(\int_{\Sigma_{2\varepsilon}}\big|g(x/\varepsilon)S_\varepsilon(f)(x)\big|^2 \delta(x)dx\Big)^{1/2}
  \leq C\|g\|_{L^2(Y)}\Big(\int_{\Sigma_{2\varepsilon}}\big|f(x)\big|^2\delta(x)dx\Big)^{1/2},
 \end{equation}
 where $C$ depends only on $d$ and $\|\zeta\|_{L^\infty(B(0,1/2))}$.
\end{lemma}

\begin{pf}
 Noting that for any $x\in\Sigma_{2\varepsilon}$, we have
 \begin{equation*}
 \begin{aligned}
 \Big|\int_{\mathbb{R}^d}\zeta_\varepsilon(x-y)f(y)dy\Big|^2
 & \leq \int_{\mathbb{R}^d}\zeta_\varepsilon(x-y)|f(y)|^2\delta^{-1}(y)dy\int_{\mathbb{R}^d}\zeta_\varepsilon(x-y)\delta(y)dy \\
 & \leq 2\delta(x) \int_{\mathbb{R}^d}\zeta_\varepsilon(x-y)|f(y)|^2\delta^{-1}(y)dy,
 \end{aligned}
 \end{equation*}
 where we use H\"older's inequality in the first inequality, and the estimate $\eqref{pri:2.10}$ in the last one. Thus
 \begin{equation*}
 \begin{aligned}
 \int_{\Sigma_{2\varepsilon}}\big|g(x/\varepsilon)S_\varepsilon(f)(x)\big|^2 \delta^{-1}(x)dx
&  \leq 2\int_{\Sigma_{2\varepsilon}} |g(x/\varepsilon)|^2\int_{\mathbb{R}^d}\zeta_\varepsilon(x-y)|f(y)|^2\delta^{-1}(y)dy dx \\
&  \leq C\sup_{y\in\mathbb{R}^d}\int_{B(y,1/2)}|g|^2 dz \int_{\Sigma_{2\varepsilon}}|f(y)|^2\delta^{-1}(y)dy \\
&   \leq C\|g\|_{L^2(Y)}^2 \int_{\Sigma_{2\varepsilon}}|f(y)|^2\delta^{-1}(y)dy.
 \end{aligned}
 \end{equation*}
 Taking square root on the both sides, we have the desired estimate $\eqref{pri:2.11}$. By the same token, we have
 \begin{equation*}
 \begin{aligned}
  \Big|\int_{\mathbb{R}^d}\zeta_\varepsilon(x-y)f(y) dy\Big|
  & \leq \int_{\mathbb{R}^d}\zeta_\varepsilon(x-y)|f(y)|^2\delta(y) dy\int_{\mathbb{R}^d}\zeta_\varepsilon(x-y)\delta^{-1}(y) dy  \\
  & \leq 2\delta^{-1}(x) \int_{\mathbb{R}^d}\zeta_\varepsilon(x-y)|f(y)|^2\delta(y) dy,
 \end{aligned}
 \end{equation*}
 and this gives
 \begin{equation*}
 \begin{aligned}
 \Big(\int_{\Sigma_{2\varepsilon}}\big|g(x/\varepsilon)S_\varepsilon(f)(x)\big|^2 \delta(x)dx\Big)^{1/2}
&  \leq \sqrt{2}\Big(\int_{\Sigma_{2\varepsilon}} |g(x/\varepsilon)|^2\int_{\mathbb{R}^d}\zeta_\varepsilon(x-y)|f(y)|^2\delta(y)dy dx\Big)^{1/2} \\
&   \leq C\|g\|_{L^2(Y)}\Big( \int_{\Sigma_{2\varepsilon}}|f(y)|^2\delta(y)dy\Big)^{1/2}.
 \end{aligned}
 \end{equation*}
We have completed the proof.
 \qed
\end{pf}

\begin{lemma}\label{lemma:2.8}
Let $f\in H^1(\Omega)$ be supported in $\Sigma_{\varepsilon}$, then we obtain
\begin{equation}\label{pri:2.13}
\Big(\int_{\Sigma_{2\varepsilon}} \big|f(x) - S_\varepsilon(f)(x)\big|^2\delta(x)dx\Big)^{1/2}
\leq C\varepsilon\Big(\int_{\Sigma_{\varepsilon}}\big|\nabla f(x)\big|^2 \delta(x) dx\Big)^{1/2}.
\end{equation}
where $C$ depends only on $d$.
\end{lemma}

\begin{pf}
Let $|y|\leq 1$, then we first obtain
\begin{equation}\label{f:2.2}
 \int_{\Sigma_{2\varepsilon}} \big|f(x) - f(x-\varepsilon y)\big|^2 \delta(y)dy
 \leq C\varepsilon^2 \int_{\Sigma_\varepsilon}\big|\nabla f(z)\big|^2 \delta(z)dz.
\end{equation}
To see this estimate, we start with
\begin{equation*}
 f(x) - f(x-\varepsilon y) = \varepsilon\int_0^1 \nabla f(x+(t-1)\varepsilon y)\cdot y dt.
\end{equation*}
Then we have
\begin{equation*}
 \big|f(x) - f(x-\varepsilon y)\big|^2\delta(x) \leq  C\varepsilon^2\int_0^1 |\nabla f(x+(t-1)\varepsilon y)|^2 \delta(x) dt.
\end{equation*}
Integrating both sides with respect to $x$ on $\Sigma_{2\varepsilon}$, we arrive at
\begin{equation*}
\begin{aligned}
 \int_{\Sigma_{2\varepsilon}}\big|f(x) - f(x-\varepsilon y)\big|^2\delta(x) dx
& \leq  C\varepsilon^2 \int_{\Sigma_{2\varepsilon}}\int_0^1 |\nabla f(x+(t-1)\varepsilon y)|^2 \delta(x) dt dx\\
& \leq C\varepsilon^2 \int_{\Sigma_{\varepsilon}}\int_0^1 |\nabla f(z)|^2 \delta(z+(1-t)\varepsilon y) dt dz
 \leq C\varepsilon^2 \int_{\Sigma_{\varepsilon}}|\nabla f(z)|^2 \delta(z)dz.
\end{aligned}
\end{equation*}
Note that since $t\in[0,1]$ and $|y|\leq 1$, we have $z = x + (t-1)\varepsilon y\in \Sigma_{\varepsilon}$. In the last inequality,
it follows from the mean value theorem that for any $z\in\Sigma_\varepsilon$,
\begin{equation*}
\delta(z + (1-t)\varepsilon y)\leq \varepsilon(1-t)|y|\|\nabla\delta(z)\|_{L^\infty(\Omega)} + \delta(z)
\leq \varepsilon + \delta(z)\leq 2\delta(z),
\end{equation*}
since $\delta(z)> \varepsilon$ and $\|\nabla\delta\|_{L^\infty(\Omega)} = 1$.

Next step, we will prove Minkowski's inequality of a weighted-type.
\begin{equation*}
\begin{aligned}
\int_{\Sigma_{2\varepsilon}}\big|f(x)-S_\varepsilon(f)(x)\big|^2 \delta(x) dx
&= \int_{\Sigma_{2\varepsilon}}\Big|\int_{\mathbb{R}^d}\zeta_\varepsilon(x-y)\big[f(x)-f(y)\big]dy\Big|^2 \delta(x) dx \\
&\leq \int_{\Sigma_{2\varepsilon}} \int_{\mathbb{R}^d} \zeta_\varepsilon(x-y)\delta^{-1}(y) dy
\int_{\mathbb{R}^d} \zeta_\varepsilon(x-y)\big|f(x)-f(y)\big|^2\delta(y) dy \delta(x) dx \\
&\leq 2\int_{\mathbb{R}^d} \zeta_\varepsilon(x-y)\int_{\Sigma_{2\varepsilon}}\big|f(x)-f(y)\big|^2\delta(y) dydx \\
& = 2\int_{|y|\leq 1} \zeta(y)\int_{\Sigma_{2\varepsilon}}\big|f(x)-f(x-\varepsilon y)\big|^2\delta(x-\varepsilon y)dxdy   \\
& \leq 4\int_{|y|\leq 1} \zeta(y)\int_{\Sigma_{2\varepsilon}}\big|f(x)-f(x-\varepsilon y)\big|^2\delta(x)dx dy\\
&\leq C\varepsilon^2\int_{\Sigma_{\varepsilon}}\big|\nabla f(z)\big|^2\delta(z) dy,
\end{aligned}
\end{equation*}
where we use H\"older's inequality in the first inequality, the estimate $\eqref{pri:2.10}$ in the second one. In the third inequality,
we observe that $\delta(x-\varepsilon y)\leq 2\delta(x)$ whenever $|y|\leq 1$ and $x\in\Sigma_{2\varepsilon}$.
Besides, we use $\eqref{f:2.2}$ in the last one.
This implies the estimate $\eqref{pri:2.13}$, and we have completed the proof.
\qed
\end{pf}

\begin{remark}\label{re:2.13}
\emph{For ease of notations, we write
\begin{equation*}
 \|f\|_{L^2(\Sigma_{r};\delta)} = \Big(\int_{\Sigma_{r}}|f(x)|^2\delta(x) dx\Big)^{1/2},
 \quad
 \|f\|_{L^2(\Sigma_{r};\delta^{-1})} = \Big(\int_{\Sigma_{r}}|f(x)|^2\delta^{-1}(x) dx\Big)^{1/2},
\end{equation*}
and $\|f\|_{H^1(\Sigma_{r};\delta)} = \|f\|_{L^2(\Sigma_{r};\delta)} + \|\nabla f\|_{L^2(\Sigma_{r};\delta)}$.
Then from H\"older's inequality, it is not hard to see
\begin{equation}\label{pri:2.14}
 \|fg\|_{L^1(\Sigma_{r})} \leq \|f\|_{L^2(\Sigma_r;\delta)}\|g\|_{L^2(\Sigma_r;\delta^{-1})}.
\end{equation}
Also, it follows from the definition that
\begin{equation}\label{pri:2.19}
  \|f\|_{H^1(\Omega\setminus\Sigma_r;\delta)}\leq \sqrt{2r}\|f\|_{H^1(\Omega\setminus\Sigma_r)}.
\end{equation}
Moreover, if $\varphi\in L^2(\Omega)$, and $g_\varepsilon,f$ are given in Lemmas $\ref{lemma:2.7}$ and $\ref{lemma:2.8}$, then we obtain
\begin{equation}\label{pri:2.15}
 \|g_\varepsilon S_\varepsilon(f)\varphi\|_{L^1(\Sigma_{2\varepsilon})}
 \leq \|g_\varepsilon S_\varepsilon(f)\|_{L^2(\Sigma_{2\varepsilon};\delta^{-1})}
 \|\varphi\|_{L^2(\Sigma_{2\varepsilon};\delta)}
 \leq C\|g\|_{L^2(Y)}\|f\|_{L^2(\Sigma_{2\varepsilon};\delta^{-1})}\|\varphi\|_{L^2(\Sigma_{2\varepsilon};\delta)},
\end{equation}
and
\begin{equation}\label{pri:2.16}
 \|g_\varepsilon S_\varepsilon(f)\varphi\|_{L^1(\Sigma_{2\varepsilon})}
 \leq \|g_\varepsilon S_\varepsilon(f)\|_{L^2(\Sigma_{2\varepsilon};\delta)}
 \|\varphi\|_{L^2(\Sigma_{2\varepsilon};\delta^{-1})}
 \leq C\|g\|_{L^2(Y)}\|f\|_{L^2(\Sigma_{2\varepsilon};\delta)}\|\varphi\|_{L^2(\Sigma_{2\varepsilon};\delta^{-1})},
\end{equation}
where we use estimates $\eqref{pri:2.14}$, $\eqref{pri:2.11}$ and $\eqref{pri:2.12}$. Similarly we have
\begin{equation}\label{pri:2.17}
 \|[f- S_\varepsilon(f)]\varphi\|_{L^1(\Sigma_{2\varepsilon})}
 \leq \|f- S_\varepsilon(f)\|_{L^2(\Sigma_{2\varepsilon};\delta)}
 \|\varphi\|_{L^2(\Sigma_{2\varepsilon};\delta^{-1})}
 \leq C\varepsilon\|\nabla f\|_{L^2(\Sigma_{\varepsilon};\delta)}\|\varphi\|_{L^2(\Sigma_{2\varepsilon};\delta^{-1})},
\end{equation}
where C depends only on $d$ and $\|\zeta\|_{L^\infty(B(0,1/2))}$.
In the end, we mention that the estimates $\eqref{pri:2.15}-\eqref{pri:2.17}$ are frequently used
in the proof of Lemma $\ref{lemma:2.10}$. }
\end{remark}

It is known that by energy methods we just have the convergence rate $O(\varepsilon^{\frac{1}{2}})$.
To accelerate the convergence rate from $O(\varepsilon^{\frac{1}{2}})$ to $O(\varepsilon)$,
we employ the duality methods which developed by T. Suslina in \cite{TS,TS2}.
Compared to theirs, we make some improvements from a technical standpoint.
There are two lemmas related to duality methods in the paper,
and the first one is more complicated than the second one. We mention that the second one
is based on the assumption of $u_0\in H^2(\Omega)$. Because of this,
we do not worry about how to control
$\|\nabla^2 u_0\|_{L^2(\Omega)}$ by the given data such as the source term and the boundary term.
Of course, if $\Omega$ is a smooth domain, there is no worry about it since we can use the
$H^2$ estimate directly.
However, this paper is concentrating on the boundary value problems concerned with Lipschitz domains.
In this case, as we have explained before the assumption of $u_0\in H^2(\Omega)$ is not very natural.
Observing that the distance function $\delta$ is helpful in cancelling
the singularity of $\nabla^2 u_0$ near the boundary of $\Omega$,
the weighted-type inequalities are expected to be established in the duality argument,
and then the following lemma makes this thinking come true.

\begin{lemma}[Duality lemma I]\label{lemma:2.10}
Let $w_\varepsilon$ be given in $\eqref{eq:2.12}$ by choosing $\varphi_0 = S_\varepsilon(\psi_{4\varepsilon}u_0)$ and
$\varphi_k = S_\varepsilon(\psi_{4\varepsilon}\nabla_k u_0)$,
where the weak solutions $u_\varepsilon$ and $u_0$ satisfy $\eqref{pde:1.2}$ or $\eqref{pde:1.3}$.
For any $\Phi\in L^2(\Omega;\mathbb{R}^m)$,
we let $\phi_\varepsilon$ and $\phi_0$ be the weak solutions to the corresponding adjoint problems $\eqref{pde:2.6}$ or $\eqref{pde:2.7}$.
Then we have
\begin{equation}\label{eq:2.13}
 \int_\Omega w_\varepsilon\Phi dx
 = \int_\Omega \tilde{f}\cdot\nabla\phi_\varepsilon dx
 + \int_\Omega \tilde{F}\phi_\varepsilon dx,
\end{equation}
where $\tilde{f}$ and $\tilde{F}$ are defined in $\eqref{eq:2.10}$. Moreover, if we assume
\begin{equation}\label{eq:2.14}
\xi_\varepsilon = \phi_\varepsilon - \phi_0 - \varepsilon\chi_{0,\varepsilon}^*S_\varepsilon(\psi_{10\varepsilon}\phi_0)
+ \varepsilon\chi_{k,\varepsilon}^*S_\varepsilon(\psi_{10\varepsilon}\nabla_k\phi_0),
\end{equation}
where $\chi_j^*$ with $j=0,\cdots,d$ are the corresponding correctors associated with $\mathcal{L}_\varepsilon^*$.
Then we obtain the estimate
\begin{equation}\label{pri:2.18}
\begin{aligned}
\Big|\int_\Omega w_\varepsilon\Phi dx\Big|
&\leq C\|u_0\|_{H^1(\Omega\setminus\Sigma_{9\varepsilon})}\|\phi_0\|_{H^1(\Omega\setminus\Sigma_{9\varepsilon})}
 + C\varepsilon\|u_0\|_{H^1(\Omega)}\|\phi_\varepsilon\|_{H^1(\Omega)} \\
& \qquad + C\Big\{\| u_0\|_{H^1(\Omega\setminus\Sigma_{8\varepsilon};\delta)}
+\varepsilon \|\nabla u_0\|_{L^2(\Sigma_{4\varepsilon};\delta)}
+\varepsilon \|\nabla^2 u_0\|_{L^2(\Sigma_{4\varepsilon};\delta)}\Big\}\|\phi_0\|_{H^1(\Sigma_{4\varepsilon};\delta^{-1})} \\
& \qquad\qquad\qquad+ C\Big\{\| u_0\|_{H^1(\Omega\setminus\Sigma_{8\varepsilon})}
+\varepsilon \|\nabla u_0\|_{L^2(\Sigma_{4\varepsilon})}
+\varepsilon \|\nabla^2 u_0\|_{L^2(\Sigma_{4\varepsilon})}\Big\}\\
&\qquad\qquad\qquad\qquad\qquad\cdot\Big\{\|\xi_\varepsilon\|_{H^1(\Omega)} + \|\phi_0\|_{H^1(\Omega\setminus\Sigma_{20\varepsilon})}
+ \varepsilon\|\phi_0\|_{H^1(\Omega)}+ \varepsilon\|\nabla^2\phi_0\|_{L^2(\Sigma_{10\varepsilon})}\Big\},
\end{aligned}
\end{equation}
where $C$ depends on $\mu,\kappa,m,d$ and $\Omega$.
\end{lemma}

\begin{pf}
First of all, we prove the equality $\eqref{eq:2.13}$ under different types of boundary conditions.
In the case of the Dirichlet problem, $u_\varepsilon$ and $u_0$ are given in Theorem $\ref{thm:1.1}$,
while $\phi_\varepsilon$ and $\phi_0$ satisfy $\eqref{pde:2.6}$ with $\Phi\in L^2(\Omega;\mathbb{R}^m)$
according to the assumption of this lemma. Hence, in view of $\eqref{eq:2.1}$ and $\eqref{pde:2.4}$
we have
\begin{equation}\label{f:2.3}
\begin{aligned}
\int_\Omega w_\varepsilon \Phi dx
 = <\mathcal{L}_\varepsilon(w_\varepsilon),\phi_\varepsilon>
 = <-\text{div}(\tilde{f})+\tilde{F},\phi_\varepsilon>
  = \int_\Omega \tilde{f}\cdot\nabla \phi_\varepsilon dx + \int_\Omega \tilde{F} \phi_\varepsilon dx,
\end{aligned}
\end{equation}
where we use the facts of $w_\varepsilon = 0$ and $\phi_\varepsilon = 0$ on $\partial\Omega$.
For the Neumann problem, $u_\varepsilon$ and $u_0$ are given in Theorem $\ref{thm:1.2}$,
and $\phi_\varepsilon$ and $\phi_0$ solve $\eqref{pde:2.7}$ for any $\Phi\in L^2(\Omega;\mathbb{R}^m)$.
It follows from $\eqref{eq:2.2}$ that
\begin{equation}\label{f:2.4}
 \int_\Omega w_\varepsilon\Phi dx
 = <-\text{div}(\tilde{f})+\tilde{F},\phi_\varepsilon> + <\mathcal{B}_\varepsilon(w_\varepsilon),\phi_\varepsilon>
 = \int_\Omega \tilde{f}\cdot\nabla \phi_\varepsilon dx + \int_\Omega \tilde{F} \phi_\varepsilon dx,
\end{equation}
where we use the facts of $\mathcal{B}_\varepsilon(w_\varepsilon) = n\cdot\tilde{f}$
and $\mathcal{B}_\varepsilon^*(\phi_\varepsilon) = 0$ on $\partial\Omega$.
Thus the equalities $\eqref{f:2.3}$ and $\eqref{f:2.4}$ show that $\eqref{eq:2.13}$ holds
for $\mathcal{L}_\varepsilon$ with either Dirichlet or Neumann boundary conditions under the assumptions of this lemma.

We now proceed to prove the estimate $\eqref{pri:2.18}$, and
the main task is to estimate the right-hand side of  $\eqref{eq:2.13}$,
and its first term can be expressed by
\begin{equation}\label{f:3.45}
 \begin{aligned}
 \Big|\int_\Omega\tilde{f}\cdot\nabla\phi_\varepsilon dx\Big|
 &\leq \Big|\int_\Omega\mathcal{K}\cdot\nabla\phi_\varepsilon dx\Big|
  + C\int_\Omega \big|\nabla u_0 - S_\varepsilon(\psi_{4\varepsilon}\nabla u_0)\big|\big|\nabla\phi_\varepsilon\big| dx \\
 &  + C\int_\Omega \big|u_0 - S_\varepsilon(\psi_{4\varepsilon}u_0)\big|\big|\nabla\phi_\varepsilon\big| dx
    + \varepsilon\int_\Omega\big|\mathcal{I}\cdot\nabla\phi_\varepsilon\big| dx
    + \varepsilon\int_\Omega\big|\mathcal{J}\cdot\nabla\phi_\varepsilon\big| dx.
 \end{aligned}
\end{equation}
Clearly, we can divide the estimate of $\eqref{f:3.45}$ in five parts.

\textbf{Part 1}: we study the term of $\int_\Omega \mathcal{K}\cdot \nabla \phi_\varepsilon dx$.  In view of $\eqref{f:3.3}$, we have
\begin{equation*}
\int_\Omega \mathcal{K}\cdot \nabla \phi_\varepsilon dx = R(\phi_\varepsilon) = R_1(\phi_\varepsilon) - R_2(\phi_\varepsilon).
\end{equation*}
It follows from the estimate $\eqref{f:3.4}$ that $R_1(\phi_\varepsilon)=0$,
because $S_\varepsilon(\psi_{4\varepsilon}u_0)$ and
$S_\varepsilon(\psi_{4\varepsilon}\nabla u_0)$ are supported in $\Sigma_{3\varepsilon}$.
We only compute the term of  $R_2(\phi_\varepsilon)$. By the definition of $\xi_\varepsilon$ in $\eqref{eq:2.14}$,
\begin{equation}\label{f:3.46}
\begin{aligned}
R_2(\phi_\varepsilon) &= \varepsilon \int_\Omega E_{jik,\varepsilon}^{\alpha\gamma}\nabla_j\big\{
S_\varepsilon(\psi_{4\varepsilon}\nabla_ku_0^\gamma)\big\}\nabla_i\phi_\varepsilon^\alpha dx\\
&= \int_\Omega E_{jik,\varepsilon}^{\alpha\gamma}S_\varepsilon(\nabla_j\psi_{4\varepsilon}\nabla_ku_0^\gamma)\nabla_i\phi_\varepsilon^\alpha dx
+\varepsilon\int_\Omega E_{jik,\varepsilon}^{\alpha\gamma}S_\varepsilon(\psi_{4\varepsilon}\nabla_{jk}^2u_0^\gamma)
\nabla_i\phi_\varepsilon^\alpha dx \\
&\leq C\|u_0\|_{H^1(\Omega\setminus\Sigma_{9\varepsilon})}\|\nabla\phi_\varepsilon\|_{L^2(\Omega\setminus\Sigma_{9\varepsilon})}\\
&+ \varepsilon\int_\Omega E_{jik,\varepsilon}^{\alpha\gamma}S_\varepsilon(\psi_{4\varepsilon}\nabla_{jk}^2u_0^\gamma)
\nabla_i\big[\xi_\varepsilon^\alpha +\phi_0^\alpha
+\varepsilon\chi_0^{\alpha\beta}(x/\varepsilon)S_\varepsilon(\psi_{10\varepsilon}\phi_0^\beta)
+\varepsilon\chi_k^{\alpha\beta}(x/\varepsilon)S_\varepsilon(\psi_{10\varepsilon}\nabla_k\phi_0^\beta)\big]dx.
\end{aligned}
\end{equation}
Note that
\begin{equation}\label{f:3.47}
\begin{aligned}
\Big|\int_\Omega E_{jik,\varepsilon}^{\alpha\gamma}S_\varepsilon(\psi_{4\varepsilon}\nabla_{jk}^2u_0^\gamma)
\nabla_i\xi_\varepsilon^\alpha dx\Big|
\leq \big\|\varpi_\varepsilon S_\varepsilon(\psi_{4\varepsilon}\nabla^2u_0)\big\|_{L^2(\Omega)}\|\nabla\xi_\varepsilon\|_{L^2(\Omega)}
\leq C\|\nabla^2 u_0\|_{L^2(\Sigma_{4\varepsilon})}\|\nabla\xi_\varepsilon\|_{L^2(\Omega)}
\end{aligned}
\end{equation}
where we use $\eqref{pri:2.5}$ in the second inequality, the notation $\varpi_\varepsilon$ represents the
periodic functions determined by $E_{jik,\varepsilon}$.
(In the following proof,
there will be a lot of algebra combinations of the different periodic functions. Ignoring
their concrete form, we often use the notation $\varpi_\varepsilon$
to denote them. We refer the reader to Remark $\ref{re:2.6}$ for more explanations. If we say
$\varpi_\varepsilon$ depends on some periodic functions, that means
$\varpi_\varepsilon$ represents one of their algebra combinations.)
\begin{equation}\label{f:3.48}
\begin{aligned}
\Big|\int_\Omega E_{jik,\varepsilon}^{\alpha\gamma}S_\varepsilon(\psi_{4\varepsilon}\nabla_{jk}^2u_0^\gamma)
\nabla_i\phi_0^\alpha dx\Big|
& \leq \big\|\varpi_\varepsilon S_\varepsilon(\psi_{4\varepsilon}\nabla^2u_0)\big\|_{L^2(\Sigma_{4\varepsilon};\delta)}
\|\nabla\phi_0\|_{L^2(\Sigma_{4\varepsilon};\delta^{-1})} \\
& \leq C\|\nabla^2 u_0\|_{L^2(\Sigma_{4\varepsilon};\delta)}\|\nabla\phi_0\|_{L^2(\Sigma_{4\varepsilon};\delta^{-1})}
\end{aligned}
\end{equation}
where we use $\eqref{pri:2.14}$ in the first inequality and $\eqref{pri:2.12}$ in the last one.
\begin{equation}\label{f:3.49}
\begin{aligned}
&\Big|\int_\Omega E_{jik,\varepsilon}^{\alpha\gamma}S_\varepsilon(\psi_{4\varepsilon}\nabla_{jk}^2u_0^\gamma)
\nabla_i\big[\varepsilon\chi_{0,\varepsilon}^{\alpha\beta}S_\varepsilon(\psi_{10\varepsilon}\phi_0^\beta)\big] dx\Big| \\
&\leq \int_\Omega \Big|E_{jik,\varepsilon}^{\alpha\gamma}S_\varepsilon(\psi_{4\varepsilon}\nabla_{jk}^2u_0^\gamma)
\big[\nabla_i\chi_{0,\varepsilon}^{\alpha\beta}S_\varepsilon(\psi_{10\varepsilon}\phi_0^\beta)
+ \varepsilon\chi_{0,\varepsilon}^{\alpha\beta}\nabla_iS_\varepsilon(\psi_{10\varepsilon}\phi_0^\beta) \big] \Big|dx \\
& \leq \big\|\varpi_\varepsilon S_\varepsilon(\psi_{4\varepsilon}\nabla^2u_0)\big\|_{L^2(\Sigma_{4\varepsilon};\delta)}
\big\|\nabla\chi_{0,\varepsilon}S_\varepsilon(\psi_{10\varepsilon}\phi_0)\big\|_{L^2(\Sigma_{4\varepsilon};\delta^{-1})}
+ \varepsilon\big\|\varpi_\varepsilon S_\varepsilon(\psi_{4\varepsilon}\nabla^2u_0)\big\|_{L^2(\Omega)}
  \big\|\chi_{0,\varepsilon}\nabla S_\varepsilon(\psi_{10\varepsilon}\phi_0)\big\|_{L^2(\Omega)} \\
& \leq C\|\nabla^2 u_0\|_{L^2(\Sigma_{4\varepsilon};\delta)}\|\phi_0\|_{L^2(\Sigma_{10\varepsilon};\delta^{-1})}
+ C\|\nabla^2 u_0\|_{L^2(\Sigma_{4\varepsilon})}\|\phi_0\|_{L^2(\Omega\setminus\Sigma_{20\varepsilon})}
+ C\varepsilon \|\nabla^2 u_0\|_{L^2(\Sigma_{4\varepsilon})}\|\nabla\phi_0\|_{L^2(\Sigma_{10\varepsilon})},
\end{aligned}
\end{equation}
where we use H\"older's inequality in the second inequality,
and the estimates $\eqref{pri:2.5}$, $\eqref{pri:2.11}$ and $\eqref{pri:2.12}$ are employed in the last one. Similarly, we acquire
\begin{equation}\label{f:3.50}
\begin{aligned}
\Big|\int_\Omega E_{jik,\varepsilon}^{\alpha\gamma}S_\varepsilon(\psi_{4\varepsilon}\nabla_{jk}^2u_0^\gamma)
&\nabla_i\big[\varepsilon\chi_{k,\varepsilon}^{\alpha\beta}S_\varepsilon(\psi_{10\varepsilon}\nabla_k\phi_0^\beta)\big] dx\Big|  \\
&\leq C\|\nabla^2 u_0\|_{L^2(\Sigma_{4\varepsilon};\delta)}\|\nabla\phi_0\|_{L^2(\Sigma_{10\varepsilon};\delta^{-1})}  \\
& \qquad\qquad+ C\|\nabla^2 u_0\|_{L^2(\Sigma_{4\varepsilon})}\|\nabla\phi_0\|_{L^2(\Omega\setminus\Sigma_{20\varepsilon})}
+ C\varepsilon \|\nabla^2 u_0\|_{L^2(\Sigma_{4\varepsilon})}\|\nabla^2\phi_0\|_{L^2(\Sigma_{10\varepsilon})}.
\end{aligned}
\end{equation}
Combining $\eqref{f:3.46}-\eqref{f:3.50}$, we have
\begin{equation}\label{f:3.51}
\begin{aligned}
\Big|\int_\Omega \mathcal{K}\cdot\nabla\varphi_\varepsilon dx \Big|
& \leq  \big|R_2(\phi_\varepsilon)\big| \\
&  \leq C\|u_0\|_{H^1(\Omega\setminus\Sigma_{9\varepsilon})}\|\nabla\phi_\varepsilon\|_{L^2(\Omega\setminus\Sigma_{9\varepsilon})}
     + C\varepsilon\|\nabla^2 u_0\|_{L^2(\Sigma_{4\varepsilon};\delta)}\|\phi_0\|_{H^1(\Sigma_{4\varepsilon};\delta^{-1})} \\
&  + C\varepsilon\|\nabla^2 u_0\|_{L^2(\Sigma_{4\varepsilon})}\Big\{\|\nabla\xi_\varepsilon\|_{L^2(\Omega)}
   + \|\phi_0\|_{H^1(\Omega\setminus\Sigma_{20\varepsilon})} + \varepsilon\|\nabla^2\phi_0\|_{L^2(\Sigma_{10\varepsilon})}
   +\varepsilon\|\phi_0\|_{H^1(\Sigma_{4\varepsilon})}\Big\},
\end{aligned}
\end{equation}
where $\|\phi_0\|_{H^1(\Sigma_{4\varepsilon};\delta^{-1})} = \|\nabla\phi_0\|_{L^2(\Sigma_{4\varepsilon};\delta^{-1})}
+ \|\phi_0\|_{L^2(\Sigma_{4\varepsilon};\delta^{-1})}$.

\textbf{Part 2}: we now consider
\begin{equation}\label{f:3.52}
\begin{aligned}
\int_\Omega\big|\nabla u_0
- S_\varepsilon(\psi_{4\varepsilon}\nabla u_0)\big||\nabla\phi_\varepsilon|dx
&\leq \int_\Omega\big|(1-\psi_{4\varepsilon})\nabla u_0\big||\nabla\phi_\varepsilon|dx
+ \int_\Omega\big|\psi_{4\varepsilon}\nabla u_0 - S_\varepsilon(\psi_{4\varepsilon}\nabla u_0)\big||\nabla\phi_\varepsilon|dx\\
&\leq C\|\nabla u_0\|_{L^2(\Omega\setminus\Sigma_{8\varepsilon})}\|\nabla\phi_\varepsilon\|_{L^2(\Omega\setminus\Sigma_{8\varepsilon})}
+ \int_\Omega\big|\psi_{4\varepsilon}\nabla u_0 - S_\varepsilon(\psi_{4\varepsilon}\nabla u_0)\big| \\
&\qquad\qquad\qquad\cdot\big|\nabla\big[\xi_\varepsilon + \phi_0 + \varepsilon\chi_{0,\varepsilon}^*S_\varepsilon(\psi_{10\varepsilon}\phi_0)
+\varepsilon\chi_{k,\varepsilon}^*S_\varepsilon(\psi_{10\varepsilon}\nabla_k\phi_0)\big]\big| dx.
\end{aligned}
\end{equation}
Note that
\begin{equation}\label{f:3.53}
\begin{aligned}
\int_\Omega\big|\psi_{4\varepsilon}\nabla u_0 - S_\varepsilon(\psi_{4\varepsilon}\nabla u_0)\big||\nabla\xi_\varepsilon|dx
&\leq \big\|\psi_{4\varepsilon}\nabla u_0 - S_\varepsilon(\psi_{4\varepsilon}\nabla u_0)\big\|_{L^2(\Omega)}
\|\nabla\xi_\varepsilon\|_{L^2(\Omega)} \\
& \leq C\varepsilon\|\nabla(\psi_{4\varepsilon}\nabla u_0)\|_{L^2(\mathbb{R}^d)}\|\nabla\xi_\varepsilon\|_{L^2(\Omega)}\\
& \leq C\big\{\|\nabla u_0\|_{L^2(\Omega\setminus\Sigma_{8\varepsilon})}
+ \varepsilon\|\nabla^2 u_0\|_{L^2(\Sigma_{4\varepsilon})}\big\}\|\nabla \xi_\varepsilon\|_{L^2(\Omega)},
\end{aligned}
\end{equation}
where we use $\eqref{pri:2.6}$ in the second inequality.
\begin{equation}\label{f:3.54}
\begin{aligned}
\int_\Omega\big|\psi_{4\varepsilon}\nabla u_0 - S_\varepsilon(\psi_{4\varepsilon}\nabla u_0)\big||\nabla\phi_0|dx
&\leq \big\|\psi_{4\varepsilon}\nabla u_0 - S_\varepsilon(\psi_{4\varepsilon}\nabla u_0)\big\|_{L^2(\Sigma_{4\varepsilon};\delta)}
\|\nabla\phi_0\|_{L^2(\Sigma_{4\varepsilon};\delta^{-1})} \\
&\leq C\varepsilon\|\nabla(\psi_{4\varepsilon}\nabla u_0)\|_{L^2(\Sigma_{\varepsilon};\delta)}\|\nabla\phi_0\|_{L^2(\Sigma_{4\varepsilon};\delta^{-1})} \\
&\leq C\big\{\|\nabla u_0\|_{L^2(\Omega\setminus\Sigma_{8\varepsilon};\delta)} + \varepsilon\|\nabla^2 u_0\|_{L^2(\Sigma_{4\varepsilon};\delta)}\big\}
\|\nabla\phi_0\|_{L^2(\Sigma_{4\varepsilon};\delta^{-1})},
\end{aligned}
\end{equation}
where we use $\eqref{pri:2.14}$ in the first inequality, and the estimate $\eqref{pri:2.13}$ in the second one.
\begin{equation}\label{f:3.55}
\begin{aligned}
\int_\Omega\big|\psi_{4\varepsilon}\nabla u_0
&- S_\varepsilon(\psi_{4\varepsilon}\nabla u_0)\big|
\big|\nabla\big[\varepsilon\chi_{0,\varepsilon}^*S_\varepsilon(\psi_{10\varepsilon}\phi_0)\big]\big|dx \\
&\leq \big\|\psi_{4\varepsilon}\nabla u_0 - S_\varepsilon(\psi_{4\varepsilon}\nabla u_0)\big\|_{L^2(\Sigma_{4\varepsilon};\delta)}
\big\|\nabla\chi_{0,\varepsilon}^*S_\varepsilon(\psi_{10\varepsilon}\phi_0)\big\|_{L^2(\Sigma_{4\varepsilon};\delta^{-1})}\\
& +\varepsilon\big\|\psi_{4\varepsilon}\nabla u_0 - S_\varepsilon(\psi_{4\varepsilon}\nabla u_0)\big\|_{L^2(\Omega)}
\big\|\chi_{0,\varepsilon}^*\nabla S_\varepsilon(\psi_{10\varepsilon}\phi_0)\big\|_{L^2(\Omega)} \\
&\leq C\varepsilon\|\nabla(\psi_{4\varepsilon}\nabla u_0)\|_{L^2(\Sigma_{\varepsilon};\delta)}
\|\psi_{10\varepsilon}\phi_0\|_{L^2(\Sigma_{4\varepsilon};\delta^{-1})}
+ C\varepsilon^2\|\nabla(\psi_{4\varepsilon}\nabla u_0)\|_{L^2(\mathbb{R}^d)}\|\nabla(\psi_{10\varepsilon}\phi_0)\|_{L^2(\mathbb{R}^d)} \\
&\leq C\big\{\|\nabla u_0\|_{L^2(\Omega\setminus\Sigma_{8\varepsilon};\delta)}
+\varepsilon \|\nabla^2 u_0\|_{L^2(\Sigma_{4\varepsilon};\delta)}\big\}
\|\phi_0\|_{L^2(\Sigma_{10\varepsilon};\delta^{-1})} \\
& \qquad\qquad\qquad + C\big\{\|\nabla u_0\|_{L^2(\Omega\setminus\Sigma_{8\varepsilon})}
+\varepsilon \|\nabla^2 u_0\|_{L^2(\Sigma_{4\varepsilon})}\big\}\|\phi_0\|_{L^2(\Omega\setminus\Sigma_{20\varepsilon})}\\
& \qquad\qquad\qquad\qquad\qquad
\qquad + C\varepsilon\big\{\|\nabla u_0\|_{L^2(\Omega\setminus\Sigma_{8\varepsilon})}
+\varepsilon \|\nabla^2 u_0\|_{L^2(\Sigma_{4\varepsilon})}\big\}\|\nabla\phi_0\|_{L^2(\Sigma_{10\varepsilon})},
\end{aligned}
\end{equation}
where we use H\"older's inequality in the first inequality, the estimates $\eqref{pri:2.6}$ and $\eqref{pri:2.13}$ in the second one.
By the same token, we also have
\begin{equation}\label{f:3.56}
\begin{aligned}
\int_\Omega\big|\psi_{4\varepsilon}\nabla u_0
&- S_\varepsilon(\psi_{4\varepsilon}\nabla u_0)\big|
\big|\nabla\big[\varepsilon\chi_{k,\varepsilon}^*S_\varepsilon(\psi_{10\varepsilon}\nabla_k\phi_0)\big]\big|dx \\
&\leq C\big\{\|\nabla u_0\|_{L^2(\Omega\setminus\Sigma_{8\varepsilon};\delta)}
+\varepsilon \|\nabla^2 u_0\|_{L^2(\Sigma_{4\varepsilon};\delta)}\big\}
\|\nabla \phi_0\|_{L^2(\Sigma_{10\varepsilon};\delta^{-1})} \\
& \qquad\qquad\qquad\qquad + C\big\{\|\nabla u_0\|_{L^2(\Omega\setminus\Sigma_{8\varepsilon})}
+\varepsilon \|\nabla^2 u_0\|_{L^2(\Sigma_{4\varepsilon})}\big\}\|\nabla\phi_0\|_{L^2(\Omega\setminus\Sigma_{20\varepsilon})}\\
& \qquad\qquad\qquad\qquad\qquad\qquad
\qquad + C\varepsilon\big\{\|\nabla u_0\|_{L^2(\Omega\setminus\Sigma_{8\varepsilon})}
+\varepsilon \|\nabla^2 u_0\|_{L^2(\Sigma_{4\varepsilon})}\big\}\|\nabla^2\phi_0\|_{L^2(\Sigma_{10\varepsilon})}.
\end{aligned}
\end{equation}
Collecting $\eqref{f:3.52}-\eqref{f:3.56}$, we obtain
\begin{equation}\label{f:3.57}
\begin{aligned}
&\int_\Omega\big|\nabla u_0
- S_\varepsilon(\psi_{4\varepsilon}\nabla u_0)\big||\nabla\phi_\varepsilon|dx\\
&\leq C\|\nabla u_0\|_{L^2(\Omega\setminus\Sigma_{8\varepsilon})}\|\nabla\phi_\varepsilon\|_{L^2(\Omega\setminus\Sigma_{8\varepsilon})}
+ C\Big\{\|\nabla u_0\|_{L^2(\Omega\setminus\Sigma_{8\varepsilon};\delta)}
+\varepsilon \|\nabla^2 u_0\|_{L^2(\Sigma_{4\varepsilon};\delta)}\Big\}\|\phi_0\|_{H^1(\Sigma_{4\varepsilon};\delta^{-1})} \\
& + C\Big\{\|\nabla u_0\|_{L^2(\Omega\setminus\Sigma_{8\varepsilon})}
+\varepsilon \|\nabla^2 u_0\|_{L^2(\Sigma_{4\varepsilon})}\Big\}
\Big\{\|\nabla \xi_\varepsilon\|_{L^2(\Omega)} + \|\phi_0\|_{H^1(\Omega\setminus\Sigma_{20\varepsilon})}
+ \varepsilon\|\nabla\phi_0\|_{L^2(\Omega)}+ \varepsilon\|\nabla^2\phi_0\|_{L^2(\Sigma_{10\varepsilon})}\Big\}.
\end{aligned}
\end{equation}

\textbf{Part 3}: an argument similar to the one used in \textbf{Part 2} shows that
\begin{equation}\label{f:3.58}
\begin{aligned}
&\int_\Omega\big| u_0
- S_\varepsilon(\psi_{4\varepsilon} u_0)\big||\nabla\phi_\varepsilon|dx\\
&\leq C\| u_0\|_{L^2(\Omega\setminus\Sigma_{8\varepsilon})}\|\nabla\phi_\varepsilon\|_{L^2(\Omega\setminus\Sigma_{8\varepsilon})}
+ C\Big\{\| u_0\|_{L^2(\Omega\setminus\Sigma_{8\varepsilon};\delta)}
+\varepsilon \|\nabla u_0\|_{L^2(\Sigma_{4\varepsilon};\delta)}\Big\}\|\phi_0\|_{H^1(\Sigma_{4\varepsilon};\delta^{-1})} \\
& + C\Big\{\| u_0\|_{L^2(\Omega\setminus\Sigma_{8\varepsilon})}
+\varepsilon \|\nabla u_0\|_{L^2(\Sigma_{4\varepsilon})}\Big\}
\Big\{\|\nabla \xi_\varepsilon\|_{L^2(\Omega)} + \|\phi_0\|_{H^1(\Omega\setminus\Sigma_{20\varepsilon})}
+ \varepsilon\|\nabla\phi_0\|_{L^2(\Omega)}+ \varepsilon\|\nabla^2\phi_0\|_{L^2(\Sigma_{10\varepsilon})}\Big\}.
\end{aligned}
\end{equation}

\textbf{Part 4}: according to the expression of $\mathcal{I}$ in $\eqref{eq:2.11}$, we first have
\begin{equation}\label{f:3.59}
\begin{aligned}
\int_\Omega\big|\mathcal{I}\cdot\nabla\phi_\varepsilon\big| dx
= \int_\Omega \Big|\big[a_{ij,\varepsilon}^{\alpha\beta}\chi_{0,\varepsilon}^{\beta\gamma}\nabla_jS_\varepsilon(\psi_{4\varepsilon}u_0)
&+a_{ij,\varepsilon}^{\alpha\beta}\chi_{k,\varepsilon}^{\beta\gamma}\nabla_jS_\varepsilon(\psi_{4\varepsilon}\nabla_k u_0)\\
&+V_{i,\varepsilon}^{\alpha\beta}\chi_{0,\varepsilon}^{\beta\gamma}S_\varepsilon(\psi_{4\varepsilon}u_0)
+V_{i,\varepsilon}^{\alpha\beta}\chi_{k,\varepsilon}^{\beta\gamma}S_\varepsilon(\psi_{4\varepsilon}\nabla_ku_0)\big]
\nabla_i\phi_\varepsilon^\alpha\Big| dx.
\end{aligned}
\end{equation}
Then we can show the estimate of $\eqref{f:3.59}$ term by term. The first one is
\begin{equation}\label{f:3.60}
\begin{aligned}
\int_\Omega
\Big|a_{ij,\varepsilon}^{\alpha\beta}\chi_{0,\varepsilon}^{\beta\gamma}\nabla_jS_\varepsilon(\psi_{4\varepsilon}u_0)\Big|
\big|\nabla_i\phi_\varepsilon^\alpha\big|dx
&\leq \varepsilon^{-1}\int_\Omega
\big|\varpi_\varepsilon S_\varepsilon(\nabla\psi_{4\varepsilon}u_0)\big|
|\nabla\phi_\varepsilon|dx
+ \int_\Omega
\big|\varpi_\varepsilon S_\varepsilon(\psi_{4\varepsilon}\nabla u_0)\big|
|\nabla\phi_\varepsilon|dx  \\
&\leq \varepsilon^{-1}\big\|\varpi_\varepsilon S_\varepsilon(\nabla\psi_{4\varepsilon}u_0)\big\|_{L^2(\mathbb{R}^d)}
\|\nabla\phi_\varepsilon\|_{L^2(\Omega\setminus\Sigma_{9\varepsilon})}
+ \big\|\varpi_\varepsilon S_\varepsilon(\psi_{4\varepsilon}\nabla u_0)\big\|_{L^2(\mathbb{R}^d)}
\|\nabla\phi_\varepsilon\|_{L^2(\Omega)}\\
&\leq C\varepsilon^{-1}\|\nabla\psi_{4\varepsilon}u_0\|_{L^2(\mathbb{R}^d)}
\|\nabla\phi_\varepsilon\|_{L^2(\Omega\setminus\Sigma_{9\varepsilon})}
+ \|\psi_{4\varepsilon}\nabla u_0\|_{L^2(\mathbb{R}^d)}
\|\nabla\phi_\varepsilon\|_{L^2(\Omega)}\\
&\leq C\varepsilon^{-1}\|u_0\|_{L^2(\Omega\setminus\Sigma_{8\varepsilon})}
\|\nabla\phi_\varepsilon\|_{L^2(\Omega\setminus\Sigma_{9\varepsilon})}
+\|\nabla u_0\|_{L^2(\Sigma_{4\varepsilon})}
\|\nabla\phi_\varepsilon\|_{L^2(\Omega)},
\end{aligned}
\end{equation}
where $\varpi_\varepsilon$ depends on $A$ and $\chi_{0}$,
and $S_\varepsilon(\nabla\psi_{4\varepsilon}u_0)$ is supported in $\Omega\setminus\Sigma_{9\varepsilon}$ in the first inequality.
The estimate $\eqref{pri:2.5}$ is used in the third one. Then we study
\begin{equation*}
\begin{aligned}
\int_\Omega
\Big|a_{ij,\varepsilon}^{\alpha\beta}\chi_{k,\varepsilon}^{\beta\gamma}\nabla_jS_\varepsilon(\psi_{4\varepsilon}\nabla_ku_0)\Big|
\big|\nabla_i\phi_\varepsilon^\alpha\big|dx
\leq \varepsilon^{-1}\int_\Omega
\big|\varpi_\varepsilon S_\varepsilon(\nabla\psi_{4\varepsilon}\nabla u_0)\big|
\big|\nabla\phi_\varepsilon\big|dx
+ \int_\Omega
\big|\varpi_\varepsilon S_\varepsilon(\psi_{4\varepsilon}\nabla^2 u_0)\big|
|\nabla\phi_\varepsilon|dx.
\end{aligned}
\end{equation*}
Note that $\phi_\varepsilon$ in the last term is replaced by
$\xi_\varepsilon + \phi_0
+\varepsilon\chi_{0,\varepsilon}^{*}S_\varepsilon(\psi_{10\varepsilon}\phi_0)
+ \varepsilon\chi_{k,\varepsilon}^{*}S_\varepsilon(\psi_{10\varepsilon}\nabla_k\phi_0)$, and then
the right-hand side of the above inequality can be controlled by
\begin{equation*}
\begin{aligned}
& \varepsilon^{-1}\big\|\varpi_\varepsilon S_\varepsilon(\nabla\psi_{4\varepsilon}\nabla u_0)\big\|_{L^2(\mathbb{R}^d)}
\|\nabla\phi_\varepsilon\|_{L^2(\Omega\setminus\Sigma_{9\varepsilon})}
+ \big\|\varpi_\varepsilon S_\varepsilon(\psi_{4\varepsilon}\nabla^2 u_0)\big\|_{L^2(\mathbb{R}^d)}
\Big\{\|\nabla\xi_\varepsilon\|_{L^2(\Omega)}
 + \big\|\chi_{0,\varepsilon}^*S_\varepsilon(\nabla\psi_{10\varepsilon}\phi_0)\big\|_{L^2(\mathbb{R}^d)} \\
&+ \varepsilon\big\|\chi_{0,\varepsilon}^*S_\varepsilon(\psi_{10\varepsilon}\nabla\phi_0)\big\|_{L^2(\mathbb{R}^d)}
 + \big\|\chi_{k,\varepsilon}^*S_\varepsilon(\nabla\psi_{10\varepsilon}\nabla_k\phi_0)\big\|_{L^2(\mathbb{R}^d)}
 + \varepsilon\big\|\chi_{k,\varepsilon}^*S_\varepsilon(\psi_{10\varepsilon}\nabla (\nabla_k\phi_0))\big\|_{L^2(\mathbb{R}^d)} \Big\} \\
& + \big\|\varpi_\varepsilon S_\varepsilon(\psi_{4\varepsilon}\nabla^2 u_0)\big\|_{L^2(\Sigma_{3\varepsilon};\delta)}
\Big\{\|\nabla\phi_0\|_{L^2(\Sigma_{3\varepsilon};\delta^{-1})}
+ \|\nabla\chi_{0,\varepsilon}^*S_\varepsilon(\psi_{10\varepsilon}\phi_0)\|_{L^2(\Sigma_{3\varepsilon};\delta^{-1})}
+ \|\nabla\chi_{k,\varepsilon}^*S_\varepsilon(\psi_{10\varepsilon}\nabla_k\phi_0)\|_{L^2(\Sigma_{3\varepsilon};\delta^{-1})}\Big\},
\end{aligned}
\end{equation*}
where we use H\"older's inequality and $\eqref{pri:2.14}$.
We continue to apply the estimates $\eqref{pri:2.5}$, $\eqref{pri:2.11}$ and $\eqref{pri:2.12}$ to the above
expression, and finally obtain
\begin{equation}\label{f:3.61}
\begin{aligned}
&\int_\Omega
\Big|a_{ij,\varepsilon}^{\alpha\beta}\chi_{k,\varepsilon}^{\beta\gamma}\nabla_jS_\varepsilon(\psi_{4\varepsilon}\nabla_ku_0)\Big|
\big|\nabla_i\phi_\varepsilon^\alpha\big|dx\\
&\leq C\varepsilon^{-1}\|\nabla u_0\|_{L^2(\Omega\setminus\Sigma_{8\varepsilon})}
\|\nabla\phi_\varepsilon\|_{L^2(\Omega\setminus\Sigma_{9\varepsilon})}
+ C\|\nabla^2 u_0\|_{L^2(\Sigma_{3\varepsilon};\delta)}
\Big\{\|\nabla\phi_0\|_{L^2(\Sigma_{3\varepsilon};\delta^{-1})}
+ \|\phi_0\|_{H^1(\Sigma_{3\varepsilon};\delta^{-1})}\Big\} \\
&+ C\|\nabla^2 u_0\|_{L^2(\Sigma_{4\varepsilon})}\Big\{
\|\nabla\xi_\varepsilon\|_{L^2(\Omega)}+ \|\phi_0\|_{H^1(\Omega\setminus\Sigma_{20\varepsilon})}
+\varepsilon\|\nabla\phi_0\|_{L^2(\Sigma_{10\varepsilon})}
+\varepsilon\|\nabla^2\phi_0\|_{L^2(\Sigma_{10\varepsilon})}\Big\} ,
\end{aligned}
\end{equation}
Then it is easy to derive
\begin{equation}\label{f:3.62}
\begin{aligned}
\int_\Omega
\Big|\big[V_{i,\varepsilon}^{\alpha\beta}\chi_{0,\varepsilon}^{\beta\gamma}
&S_\varepsilon(\psi_{4\varepsilon}u_0)
+V_{i,\varepsilon}^{\alpha\beta}\chi_{k,\varepsilon}^{\beta\gamma}S_\varepsilon(\psi_{4\varepsilon}\nabla_ku_0)\big]
\nabla_i\phi_\varepsilon^\alpha\Big| dx \\
&\leq \Big\{\big\|\varpi_\varepsilon S_\varepsilon(\psi_{4\varepsilon}u_0)\big\|_{L^2(\mathbb{R}^d)}
+\big\|\varpi_\varepsilon S_\varepsilon(\psi_{4\varepsilon}\nabla u_0)\big\|_{L^2(\mathbb{R}^d)}\Big\}\|\nabla\phi_\varepsilon\|_{L^2(\Omega)}
\leq C\|u_0\|_{H^1(\Omega)}\|\nabla\phi_\varepsilon\|_{L^2(\Omega)},
\end{aligned}
\end{equation}
where $\varpi_\varepsilon$ depends on
$V,\chi_k$ with $k=0,\cdots,d$, and we use $\eqref{pri:2.5}$ in the last inequality.

Combining $\eqref{f:3.59}-\eqref{f:3.62}$, we obtain
\begin{equation}\label{f:3.63}
\begin{aligned}
\int_\Omega\big|\mathcal{I}\cdot\nabla\phi_\varepsilon\big| dx
&\leq C\varepsilon^{-1}\|\nabla u_0\|_{L^2(\Omega\setminus\Sigma_{8\varepsilon})}
\|\nabla\phi_\varepsilon\|_{L^2(\Omega\setminus\Sigma_{9\varepsilon})}
+ \|u_0\|_{H^1(\Omega)}\|\nabla\phi_\varepsilon\|_{L^2(\Omega)} \\
&+ C\|\nabla^2 u_0\|_{L^2(\Sigma_{4\varepsilon})}\Big\{
\|\nabla\xi_\varepsilon\|_{L^2(\Omega)}+ \|\phi_0\|_{H^1(\Omega\setminus\Sigma_{20\varepsilon})}
+\varepsilon\|\nabla\phi_0\|_{L^2(\Sigma_{10\varepsilon})}
+\varepsilon\|\nabla^2\phi_0\|_{L^2(\Sigma_{10\varepsilon})}\Big\} \\
& + C\|\nabla^2 u_0\|_{L^2(\Sigma_{3\varepsilon};\delta)}
\Big\{\|\nabla\phi_0\|_{L^2(\Sigma_{3\varepsilon};\delta^{-1})}
+ \|\phi_0\|_{H^1(\Sigma_{3\varepsilon};\delta^{-1})}\Big\}.
\end{aligned}
\end{equation}

\textbf{Part 5}: we make the procedure as in previous parts, and it is not hard to show
\begin{equation}\label{f:3.64}
\begin{aligned}
\int_\Omega\big|\mathcal{J}\cdot\nabla\phi_\varepsilon\big| dx
&\leq \int_\Omega\Big|\Big(\frac{\partial\vartheta_0^{\alpha\gamma}}{\partial y_i}\Big)_\varepsilon S_\varepsilon(\psi_{4\varepsilon}u_0^\gamma)
+ \Big(\frac{\partial\vartheta_k^{\alpha\gamma}}{\partial y_i}\Big)_\varepsilon S_\varepsilon(\psi_{4\varepsilon}\nabla_k u_0^\gamma)\Big|
|\nabla_i\phi_\varepsilon^\alpha| dx \\
&\leq \Big\{\big\|\varpi_\varepsilon S_\varepsilon(\psi_{4\varepsilon}u_0)\big\|_{L^2(\mathbb{R}^d)}
+ \big\|\varpi_\varepsilon S_\varepsilon(\psi_{4\varepsilon}\nabla u_0)\big\|_{L^2(\mathbb{R}^d)}\Big\}
\|\nabla\phi_\varepsilon\|_{L^2(\Omega)}
\leq C\|u_0\|_{H^1(\Omega)}\|\nabla\phi_\varepsilon\|_{L^2(\Omega)},
\end{aligned}
\end{equation}
where $y = x/\varepsilon$, and $\varpi_\varepsilon$ depends on $\nabla\vartheta_k$ with $k=0,\cdots,d$.
We use the estimate $\eqref{pri:2.5}$ in the last inequality.

Hence, we can summarize the five parts, and it follows from $\eqref{f:3.45}$,
$\eqref{f:3.51}$, $\eqref{f:3.57}$, $\eqref{f:3.58}$, $\eqref{f:3.63}$ and $\eqref{f:3.64}$ that
\begin{equation}\label{f:3.65}
 \begin{aligned}
 \Big|\int_\Omega\tilde{f}\cdot\nabla\phi_\varepsilon dx\Big|
&\leq C\|u_0\|_{H^1(\Omega\setminus\Sigma_{9\varepsilon})}\|\nabla\phi_\varepsilon\|_{L^2(\Omega\setminus\Sigma_{9\varepsilon})}
 + C\varepsilon\|u_0\|_{H^1(\Omega)}\|\nabla\phi_\varepsilon\|_{L^2(\Omega)} \\
& + C\Big\{\| u_0\|_{H^1(\Omega\setminus\Sigma_{8\varepsilon};\delta)}
+\varepsilon \|\nabla u_0\|_{L^2(\Sigma_{4\varepsilon};\delta)}
+\varepsilon \|\nabla^2 u_0\|_{L^2(\Sigma_{4\varepsilon};\delta)}\Big\}\|\phi_0\|_{H^1(\Sigma_{4\varepsilon};\delta^{-1})} \\
& \qquad\qquad+ C\Big\{\| u_0\|_{H^1(\Omega\setminus\Sigma_{8\varepsilon})}
+\varepsilon \|\nabla u_0\|_{L^2(\Sigma_{4\varepsilon})}
+\varepsilon \|\nabla^2 u_0\|_{L^2(\Sigma_{4\varepsilon})}\Big\}\\
&\qquad\qquad\qquad\qquad\cdot\Big\{\|\nabla \xi_\varepsilon\|_{L^2(\Omega)} + \|\phi_0\|_{H^1(\Omega\setminus\Sigma_{20\varepsilon})}
+ \varepsilon\|\nabla\phi_0\|_{L^2(\Omega)}+ \varepsilon\|\nabla^2\phi_0\|_{L^2(\Sigma_{10\varepsilon})}\Big\}.
 \end{aligned}
\end{equation}

Using the same argument as in the proof of $\eqref{f:3.65}$, we can easily carry out
the estimate for the second term
in the right-hand side of $\eqref{eq:2.13}$, that is
\begin{equation}\label{f:3.66}
\Big|\int_\Omega\tilde{F}\phi_\varepsilon dx\Big|
\leq C\int_\Omega\big|\nabla u_0 - S_\varepsilon(\psi_{4\varepsilon}\nabla u_0)\big||\phi_\varepsilon|dx
+ C\int_\Omega\big| u_0 - S_\varepsilon(\psi_{4\varepsilon} u_0)\big||\phi_\varepsilon|dx
+ \varepsilon\int_\Omega|(\mathcal{M}+\mathcal{N})\phi_\varepsilon|dx.
\end{equation}
It is clear to see that the first term in the right-hand side of $\eqref{f:3.66}$ is similar to the proof
in \textbf{Part 2}.
Thus we show the computations without explanations:
\begin{equation}\label{f:3.67}
\begin{aligned}
\int_\Omega\big|\nabla u_0
&- S_\varepsilon(\psi_{4\varepsilon}\nabla u_0)\big||\phi_\varepsilon|dx
\leq \int_\Omega(1-\psi_{4\varepsilon})|\nabla u_0||\phi_\varepsilon| dx
+ \int_\Omega\big|\psi_{4\varepsilon}\nabla u_0 - S_\varepsilon(\psi_{4\varepsilon}\nabla u_0)\big||\phi_\varepsilon|dx \\
&\leq \|\nabla u_0\|_{L^2(\Omega\setminus\Sigma_{8\varepsilon})}\|\phi_\varepsilon\|_{L^2(\Omega\setminus\Sigma_{8\varepsilon})}
+ \big\|\psi_{4\varepsilon}\nabla u_0 - S_\varepsilon(\psi_{4\varepsilon}\nabla u_0)\big\|_{L^2(\Sigma_{3\varepsilon;\delta})}
\|\phi_0\|_{L^2(\Sigma_{3\varepsilon};\delta^{-1})}\\
& + \big\|\psi_{4\varepsilon}\nabla u_0 - S_\varepsilon(\psi_{4\varepsilon}\nabla u_0)\big\|_{L^2(\mathbb{R}^d)}
\Big\{\|\xi_\varepsilon\|_{L^2(\Omega)}
+ \varepsilon\|\chi_{0,\varepsilon}^*S_\varepsilon(\psi_{10\varepsilon}\phi_0)\|_{L^2(\mathbb{R}^d)}
+ \varepsilon\|\chi_{0,\varepsilon}^*S_\varepsilon(\psi_{10\varepsilon}\phi_0)\|_{L^2(\mathbb{R}^d)} \Big\}\\
&\leq \|\nabla u_0\|_{L^2(\Omega\setminus\Sigma_{8\varepsilon})}\|\phi_\varepsilon\|_{L^2(\Omega\setminus\Sigma_{8\varepsilon})}
+ C\Big\{\|\nabla u_0\|_{L^2(\Omega\setminus\Sigma_{8\varepsilon};\delta)}
+ \varepsilon\|\nabla^2 u_0\|_{L^2(\Sigma_{4\varepsilon};\delta)}\Big\}\|\phi_0\|_{L^2(\Sigma_{3\varepsilon};\delta^{-1})}\\
&+ C\Big\{\|\nabla u_0\|_{L^2(\Omega\setminus\Sigma_{8\varepsilon})}
+\varepsilon\|\nabla^2 u_0\|_{L^2(\Sigma_{4\varepsilon})}\Big\}\Big\{\|\xi_\varepsilon\|_{L^2(\Omega)}
+ \varepsilon\|\phi_0\|_{H^1(\Omega)}\Big\}.
\end{aligned}
\end{equation}
By the same token, it is not hard to see that
\begin{equation}\label{f:3.68}
\begin{aligned}
\int_{\Omega}\big|u_0 &-S_\varepsilon(\psi_{4\varepsilon}u_0)\big||\phi_\varepsilon| dx \\
&\leq \|u_0\|_{L^2(\Omega\setminus\Sigma_{8\varepsilon})}\|\phi_\varepsilon\|_{L^2(\Omega\setminus\Sigma_{8\varepsilon})}
+ C\Big\{\|u_0\|_{L^2(\Omega\setminus\Sigma_{8\varepsilon};\delta)}
+ \varepsilon\|\nabla u_0\|_{L^2(\Sigma_{4\varepsilon};\delta)}\Big\}\|\phi_0\|_{L^2(\Sigma_{3\varepsilon};\delta^{-1})}\\
&+ C\Big\{\|u_0\|_{L^2(\Omega\setminus\Sigma_{8\varepsilon})}
+\varepsilon\|\nabla u_0\|_{L^2(\Sigma_{4\varepsilon})}\Big\}\Big\{\|\xi_\varepsilon\|_{L^2(\Omega)}
+ \varepsilon\|\phi_0\|_{H^1(\Omega)}\Big\}.
\end{aligned}
\end{equation}
The rest thing is to estimate the third term in the right-hand side of $\eqref{f:3.66}$,
which is similar to the proof in
\textbf{Parts 4} and \textbf{5}. The core idea is still that
$\phi_\varepsilon$ is replaced by $\xi_\varepsilon$ in some proper place,
and we have practiced this argument serval times before.
Due to the expressions of $\mathcal{M}$ and $\mathcal{N}$ in $\eqref{eq:2.11}$, we have
\begin{equation}\label{f:3.69}
\begin{aligned}
\int_{\Omega}\big|(\mathcal{M}
&+\mathcal{N})\phi_\varepsilon\big| dx
\leq \int_\Omega\Big|\varpi_\varepsilon\nabla\big[S_\varepsilon(\psi_{4\varepsilon}u_0)+S_\varepsilon(\psi_{4\varepsilon}\nabla u_0)\Big]
\phi_\varepsilon\Big| dx
+  \int_\Omega\Big|\varpi_\varepsilon\big[S_\varepsilon(\psi_{4\varepsilon}u_0)+S_\varepsilon(\psi_{4\varepsilon}\nabla u_0)\Big]
\phi_\varepsilon\Big| dx \\
&\leq \varepsilon^{-1}\Big\{\big\|\varpi_\varepsilon S_\varepsilon(\nabla\psi_{4\varepsilon}u_0)\big\|_{L^2(\mathbb{R}^d)}
+\big\|\varpi_\varepsilon S_\varepsilon(\nabla\psi_{4\varepsilon}\nabla u_0)\big\|_{L^2(\mathbb{R}^d)}\Big\}
\|\phi_\varepsilon\|_{L^2(\Omega\setminus\Sigma_{8\varepsilon})}
+ \|\varpi_\varepsilon S_\varepsilon(\psi_{4\varepsilon}\nabla u_0)\|_{L^2(\mathbb{R}^d)}\|\phi_\varepsilon\|_{L^2(\Sigma_{3\varepsilon})} \\
& + \big\|\varpi_\varepsilon S_\varepsilon(\psi_{4\varepsilon}\nabla^2 u_0)\big\|_{L^2(\mathbb{R}^d)}
\Big\{\|\xi_\varepsilon\|_{L^2(\Omega)}+\varepsilon\|\phi_0\|_{H^1(\Omega)}\Big\}
+ \big\|\varpi_\varepsilon S_\varepsilon(\psi_{4\varepsilon}\nabla^2 u_0)\big\|_{L^2(\Sigma_{3\varepsilon;\delta})}
\|\phi_0\|_{L^2(\Sigma_{3\varepsilon};\delta^{-1})} \\
&+ \Big\{\big\|\varpi_\varepsilon S_\varepsilon(\psi_{4\varepsilon} u_0)\big\|_{L^2(\mathbb{R}^d)}
+ \big\|\varpi_\varepsilon S_\varepsilon(\psi_{4\varepsilon}\nabla u_0)\big\|_{L^2(\mathbb{R}^d)}\Big\}
\|\phi_\varepsilon\|_{L^2(\Sigma_{3\varepsilon})}\\
& \leq C\varepsilon^{-1}\|u_0\|_{H^1(\Omega\setminus\Sigma_{8\varepsilon})}\|\phi_\varepsilon\|_{H^1(\Omega\setminus\Sigma_{8\varepsilon})}
+ C\|u_0\|_{H^1(\Omega)}\|\phi_\varepsilon\|_{L^2(\Omega)}
+ C\|\nabla^2 u_0\|_{L^2(\Sigma_{3\varepsilon};\delta)}\|\phi_0\|_{L^2(\Sigma_{3\varepsilon};\delta^{-1})}\\
&+ C\|\nabla^2 u_0\|_{L^2(\Sigma_{4\varepsilon})}\Big\{\|\xi_\varepsilon\|_{L^2(\Omega)}
+\varepsilon\|\phi_0\|_{H^1(\Omega)}\Big\}
\end{aligned}
\end{equation}
where $\varpi_\varepsilon$ depends on
the coefficients $B,c,\lambda$ and auxiliary functions $\nabla\vartheta_k$
with $k=0,\cdots,d$.

Collecting $\eqref{f:3.66}-\eqref{f:3.69}$, we arrive at
\begin{equation}\label{f:3.70}
\begin{aligned}
\Big|\int_\Omega\tilde{F}\phi_\varepsilon dx\Big|
&\leq \|u_0\|_{H^1(\Omega\setminus\Sigma_{8\varepsilon})}\|\phi_\varepsilon\|_{L^2(\Omega\setminus\Sigma_{8\varepsilon})}
+ C\varepsilon\|u_0\|_{H^1(\Omega)}\|\phi_\varepsilon\|_{L^2(\Omega)}\\
&+ C\Big\{\|u_0\|_{H^1(\Omega\setminus\Sigma_{8\varepsilon};\delta)}
+ \varepsilon\|\nabla u_0\|_{L^2(\Sigma_{4\varepsilon};\delta)}
+ \varepsilon\|\nabla^2 u_0\|_{L^2(\Sigma_{4\varepsilon};\delta)}\Big\}\|\phi_0\|_{L^2(\Sigma_{3\varepsilon};\delta^{-1})}\\
&+ C\Big\{\|u_0\|_{H^1(\Omega\setminus\Sigma_{8\varepsilon})}
+\varepsilon\|\nabla u_0\|_{L^2(\Sigma_{4\varepsilon})}
+\varepsilon\|\nabla^2 u_0\|_{L^2(\Sigma_{4\varepsilon})}\Big\}
\Big\{\|\xi_\varepsilon\|_{L^2(\Omega)}
+ \varepsilon\|\phi_0\|_{H^1(\Omega)}\Big\}.
\end{aligned}
\end{equation}

Consequently, plugging $\eqref{f:3.65}$ and $\eqref{f:3.70}$ back into $\eqref{eq:2.13}$, we obtain
\begin{equation}\label{f:3.71}
\begin{aligned}
\Big|\int_\Omega w_\varepsilon\Phi dx\Big|
&\leq C\|u_0\|_{H^1(\Omega\setminus\Sigma_{9\varepsilon})}\|\phi_\varepsilon\|_{H^1(\Omega\setminus\Sigma_{9\varepsilon})}
 + C\varepsilon\|u_0\|_{H^1(\Omega)}\|\phi_\varepsilon\|_{H^1(\Omega)} \\
& + C\Big\{\| u_0\|_{H^1(\Omega\setminus\Sigma_{8\varepsilon};\delta)}
+\varepsilon \|\nabla u_0\|_{L^2(\Sigma_{4\varepsilon};\delta)}
+\varepsilon \|\nabla^2 u_0\|_{L^2(\Sigma_{4\varepsilon};\delta)}\Big\}\|\phi_0\|_{H^1(\Sigma_{4\varepsilon};\delta^{-1})} \\
& \qquad\qquad+ C\Big\{\| u_0\|_{H^1(\Omega\setminus\Sigma_{8\varepsilon})}
+\varepsilon \|\nabla u_0\|_{L^2(\Sigma_{4\varepsilon})}
+\varepsilon \|\nabla^2 u_0\|_{L^2(\Sigma_{4\varepsilon})}\Big\}\\
&\qquad\qquad\qquad\qquad\cdot\Big\{\|\xi_\varepsilon\|_{H^1(\Omega)} + \|\phi_0\|_{H^1(\Omega\setminus\Sigma_{20\varepsilon})}
+ \varepsilon\|\phi_0\|_{H^1(\Omega)}+ \varepsilon\|\nabla^2\phi_0\|_{L^2(\Sigma_{10\varepsilon})}\Big\}.
\end{aligned}
\end{equation}

The rest task is to handle the term of $\|\phi_\varepsilon\|_{H^1(\Omega\setminus\Sigma_{9\varepsilon})}$.
Hence, in view of $\eqref{eq:2.14}$, we instead $\phi_\varepsilon$ by the first order corrector $\xi_\varepsilon$.
Observing that $S_\varepsilon(\psi_{10\varepsilon}\phi_0)$ and $S_\varepsilon(\psi_{10\varepsilon}\nabla_j \phi_0)$ are supported in
$\Sigma_{9\varepsilon}$, we arrive at
\begin{equation}\label{f:3.72}
\|\phi_\varepsilon\|_{H^1(\Omega\setminus\Sigma_{9\varepsilon})}
\leq \|\xi_\varepsilon\|_{H^1(\Omega\setminus\Sigma_{9\varepsilon})}
+ \|\phi_0\|_{H^1(\Omega\setminus\Sigma_{9\varepsilon})}.
\end{equation}
Finally, inserting $\eqref{f:3.72}$ into $\eqref{f:3.71}$ we have the desired estimate $\eqref{pri:2.18}$,
and the proof is completed.
\qed
\end{pf}

\begin{lemma}[Duality lemma II]\label{lemma:2.11}
Assume $u_0\in H^2(\Omega)$ satisfy $\eqref{pde:1.2}$ or $\eqref{pde:1.3}$.
Let $w_\varepsilon$ be given in $\eqref{eq:2.10}$ by choosing $\varphi_0 = S_\varepsilon(\psi_{4\varepsilon}\tilde{u}_0)$ and
$\varphi_k = S_\varepsilon(\psi_{4\varepsilon}\nabla_k \tilde{u}_0)$,
where the weak solutions $\tilde {u}_0$ is the extension of $u_0$ such that $\tilde{u}_0 = u$ on $\Omega$ and
$\|\tilde{u}_0\|_{H^2(\mathbb{R}^d)}\leq C\|u_0\|_{H^2(\Omega)}$.
For any $\Phi\in L^q(\Omega;\mathbb{R}^m)$ with $q=2d/(d+1)$,
we let $\phi_\varepsilon$ and $\phi_0$ be the weak solutions to the corresponding adjoint problems $\eqref{pde:2.6}$ or $\eqref{pde:2.7}$.
Then the equation $\eqref{eq:2.13}$ also follows. Moreover, we have
\begin{equation}\label{pri:2.20}
\Big|\int_\Omega w_\varepsilon\Phi dx\Big|
\leq C\Big\{\|u_0\|_{H^1(\Omega\setminus\Sigma_{8\varepsilon})}\|\phi_\varepsilon\|_{H^1(\Omega\setminus\Sigma_{9\varepsilon})}
+ \varepsilon\|u_0\|_{H^2(\Omega)}\|\phi_\varepsilon\|_{H^1(\Omega)}\Big\},
\end{equation}
where $C$ depends on $\mu,\kappa,m,d$ and $\Omega$.
\end{lemma}

\begin{pf}
The proof of this lemma originally appeared in \cite[Lemma 5.5]{QXS1} for a Neumann problem.
We provide the reader with a proof for the sake of completeness.

In view of $\eqref{eq:2.13}$, $\eqref{eq:2.10}$ and $\eqref{eq:2.11}$, we have
\begin{equation}\label{f:3.83}
\begin{aligned}
 \Big|\int_\Omega w_\varepsilon \Phi dx \Big|
 & \leq \Big|\int_\Omega \tilde{f}\cdot\nabla \phi_\varepsilon dx\Big| +  \Big|\int_\Omega \tilde{F} \phi_\varepsilon dx \big| \\
 & \leq \Big|\int_\Omega \mathcal{K}\cdot\nabla\phi_\varepsilon dx\Big|
 + C\int_\Omega\Big(\big|\nabla\myu{u}-S_\varepsilon(\psi_{4\varepsilon}\nabla\myu{u})\big|
 + \big|\myu{u}-S_\varepsilon(\psi_{4\varepsilon}\myu{u})\big|\Big)
 \big(|\nabla\phi_\varepsilon|+|\phi_\varepsilon|\big) dx   \\
 & + \varepsilon\int_\Omega \big|\mathcal{I}+\mathcal{J}\big||\nabla\phi_\varepsilon|dx
 + \varepsilon\int_\Omega \big|\mathcal{M}+\mathcal{N}\big||\phi_\varepsilon|dx    =:I_1 + I_2 + I_3 + I_4.          
\end{aligned}
\end{equation}

Below we do some calculations in more details. We first estimate $I_1$. In view of $\eqref{f:3.3}$, we have
\begin{equation*}
\int_\Omega \mathcal{K}\cdot \nabla \phi_\varepsilon dx = R(\phi_\varepsilon) = R_1(\phi_\varepsilon) - R_2(\phi_\varepsilon).
\end{equation*}
It follows from the estimate $\eqref{f:3.4}$ that $R_1(\phi_\varepsilon)=0$,
because $S_\varepsilon(\psi_{4\varepsilon}\myu{u})$ and
$S_\varepsilon(\psi_{4\varepsilon}\nabla\myu{u})$ are supported in $\Sigma_{3\varepsilon}$.
The rest thing is to estimate $R_2(\phi_\varepsilon)$. By noting that $\nabla_0\myu{u}$ means $\myu{u}$, we have
\begin{equation*}
\begin{aligned}
R_2(\phi_\varepsilon)
& = \varepsilon \int_\Omega E_{jik,\varepsilon}^{\alpha\gamma}\nabla_j\big\{
S_\varepsilon(\psi_{4\varepsilon}\nabla_k\myu{u}^\gamma)\big\}\nabla_i\phi_\varepsilon^\alpha dx
+  \varepsilon \int_\Omega E_{jik,\varepsilon}^{\alpha\gamma}\nabla_j\big\{
S_\varepsilon(\psi_{4\varepsilon}\nabla_k\myu{u}^\gamma)\big\}\nabla_i\phi_\varepsilon^\alpha dx  \\
&\leq \int_\Omega\Big(\big|E_{jik,\varepsilon}^{\alpha\gamma}
S_\varepsilon(\nabla_j\psi_{4\varepsilon}\nabla_k\myu{u}^\gamma)\big|
+ \varepsilon\big|E_{jik,\varepsilon}^{\alpha\gamma}
S_\varepsilon(\psi_{4\varepsilon}\nabla_{kj}^2\myu{u}^\gamma)\big|\Big)|\nabla_i\phi_\varepsilon^\alpha| dx \\
&\leq \|\varpi_\varepsilon S_\varepsilon(\nabla\psi_{4\varepsilon}(\nabla\myu{u}+\myu{u}))\|_{L^2(\mathbb{R}^d)}
\|\nabla\phi_\varepsilon\|_{L^2(\Omega\setminus\Sigma_{9\varepsilon})}
+\varepsilon\|\varpi_\varepsilon S_\varepsilon(\psi_{4\varepsilon}(\nabla^2\myu{u}+\nabla\myu{u}))\|_{L^2(\mathbb{R}^d)}
\|\nabla\phi_\varepsilon\|_{L^2(\Omega)} \\
&\leq \|\nabla\psi_{4\varepsilon}(\nabla\myu{u}+\myu{u})\|_{L^2(\mathbb{R}^d)}
\|\nabla\phi_\varepsilon\|_{L^2(\Omega\setminus\Sigma_{9\varepsilon})}
+ C\varepsilon \|\psi_{4\varepsilon}(\nabla^2\myu{u}+\nabla\myu{u})\|_{L^2(\mathbb{R}^d)}\|\nabla\phi_\varepsilon\|_{L^2(\Omega)} \\
&\leq C\|u_0\|_{H^1(\Omega\setminus\Sigma_{9\varepsilon})}\|\nabla\phi_\varepsilon\|_{L^2(\Omega\setminus\Sigma_{9\varepsilon})}
+ C\varepsilon\|u_0\|_{H^2(\Omega)}\|\nabla\phi_\varepsilon\|_{L^2(\Omega)},
\end{aligned}
\end{equation*}
where $\varpi$ depending on $E$ is a periodic function,
and the first term in the first inequality is supported in $\Omega\setminus\Sigma_{9\varepsilon}$.
Also, we employ the estimate $\eqref{pri:2.5}$ in the second inequality. The above inequality gives
\begin{equation}\label{f:3.84}
 I_1 \leq C\|u_0\|_{H^1(\Omega\setminus\Sigma_{9\varepsilon})}\|\nabla\phi_\varepsilon\|_{L^2(\Omega\setminus\Sigma_{9\varepsilon})}
+ C\varepsilon\|u_0\|_{H^2(\Omega)}\|\nabla\phi_\varepsilon\|_{L^2(\Omega)}.
\end{equation}

By the fact that $\myu{u}$ is the extension of $u_0$,
we have $(1-\psi_{4\varepsilon})\nabla\myu{u} = (1-\psi_{4\varepsilon})\nabla u_0$ and
 $(1-\psi_{4\varepsilon})\myu{u} = (1-\psi_{4\varepsilon}) u_0$ on $\Omega$.
Furthermore, we have
\begin{equation*}
\begin{aligned}
\nabla u_0 - S_\varepsilon(\psi_{4\varepsilon}\nabla\myu{u})
&= \big[\nabla\myu{u} - S_\varepsilon(\nabla\myu{u})\big]
+ \big[S_\varepsilon\big((1-\psi_{4\varepsilon})\nabla\myu{u}\big)
- (1-\psi_{4\varepsilon})\nabla \myu{u}\big]
+ (1-\psi_{4\varepsilon})\nabla u_0
\quad &\text{in}~ \Omega, \\
u_0-S_\varepsilon(\psi_{4\varepsilon}\myu{u})
& =\big[\myu{u}-S_\varepsilon(\myu{u})\big]
+ \big[S_\varepsilon\big((1-\psi_{4\varepsilon})\myu{u}\big)
-(1-\psi_{4\varepsilon})\myu{u}\big]
+(1-\psi_{4\varepsilon})u_0, \quad &\text{in}~ \Omega,
\end{aligned}
\end{equation*}
and then
\begin{equation}\label{f:3.85}
\begin{aligned}
I_2
&\leq C\int_\Omega \Big\{\big|\nabla\myu{u} -S_\varepsilon(\nabla\myu{u})\big|
+ \big|\myu{u} -S_\varepsilon(\myu{u})\big|
+ \big|(1-\psi_{4\varepsilon})\nabla\myu{u} -S_\varepsilon((1-\psi_{4\varepsilon})\nabla\myu{u})\big| \\
& \qquad\qquad\qquad\qquad\qquad\quad+ \big|(1-\psi_{4\varepsilon})\myu{u} -S_\varepsilon((1-\psi_{4\varepsilon})\myu{u})\big|
+ (1-\psi_{4\varepsilon})(|\nabla\myu{u}|+|\myu{u}|)\Big\}
\big(|\nabla\phi_\varepsilon|+|\phi_\varepsilon|\big) dx \\
&\leq C\Big\{\big\|\nabla\myu{u} -S_\varepsilon(\nabla\myu{u})\big\|_{L^2(\mathbb{R}^d)}
+ \big\|\myu{u} -S_\varepsilon(\myu{u})\big\|_{L^2(\mathbb{R}^d)}\Big\}\|\phi_\varepsilon\|_{H^1(\Omega)}
+ C\|\myu{u}\|_{H^1(\Omega\setminus\Sigma_{8\varepsilon})}\|\phi_\varepsilon\|_{H^1(\Omega\setminus\Sigma_{8\varepsilon})}\\
& + C\Big\{\big\|(1-\psi_{4\varepsilon})\nabla\myu{u} -S_\varepsilon((1-\psi_{4\varepsilon})\nabla\myu{u})\big\|_{L^2(\mathbb{R}^d)}
+ \big\|(1-\psi_{4\varepsilon})\myu{u} -S_\varepsilon((1-\psi_{4\varepsilon})\myu{u})\big\|_{L^2(\mathbb{R}^d)}\Big\}
\|\phi_\varepsilon\|_{H^1(\Omega\setminus\Sigma_{9\varepsilon})} \\
&\leq C\varepsilon\Big\{\|\nabla^2\myu{u}\|_{L^2(\mathbb{R}^d)}+\|\nabla^2\myu{u}\|_{L^2(\mathbb{R}^d)}\Big\}
\|\phi_\varepsilon\|_{H^1(\Omega)}
+ C\|u_0\|_{H^1(\Omega\setminus\Sigma_{8\varepsilon})}\|\phi_\varepsilon\|_{H^1(\Omega\setminus\Sigma_{8\varepsilon})}\\
& + C\varepsilon\Big\{\|\nabla[(1-\psi_{4\varepsilon})\nabla\myu{u}]\|_{L^2(\mathbb{R}^d)}
+\|\nabla[(1-\psi_{4\varepsilon})\myu{u}]\|_{L^2(\mathbb{R}^d)} \Big\}
\|\phi_\varepsilon\|_{H^1(\Omega\setminus\Sigma_{8\varepsilon})}\\
&\leq C\varepsilon\|u_0\|_{H^2(\Omega)}\|\phi_\varepsilon\|_{H^1(\Omega)}
+ C\|u_0\|_{H^1(\Omega\setminus\Sigma_{8\varepsilon})}\|\phi_\varepsilon\|_{H^1(\Omega\setminus\Sigma_{8\varepsilon})}
\end{aligned}
\end{equation}
where we use Cauchy's inequality to derive the second inequality
and the observation that
$S_\varepsilon\big((1-\psi_{4\varepsilon})\nabla\myu{u}\big)$ restricted to $\Omega$ is supported in $\Omega\setminus\Sigma_{9\varepsilon}$.
In the third one, we employ the estimate $\eqref{pri:2.6}$,
and the last one follows from
\begin{equation*}
\begin{aligned}
\|\nabla[(1-\psi_{4\varepsilon})\nabla\myu{u}]\|_{L^2(\mathbb{R}^d)}
& +\|\nabla[(1-\psi_{4\varepsilon})\myu{u}]\|_{L^2(\mathbb{R}^d)} \\
&\leq C\varepsilon^{-1}\Big\{\|\nabla \myu{u}\|_{L^2(\Omega\setminus\Sigma_{8\varepsilon})}
+ \|\myu{u}\|_{L^2(\Omega\setminus\Sigma_{8\varepsilon})}\Big\}
+ \|\nabla^2\myu{u}\|_{L^2(\mathbb{R}^d)}
+ \|\nabla\myu{u}\|_{L^2(\mathbb{R}^d)} \\
& \leq C\varepsilon^{-1}\|\nabla u_0\|_{H^1(\Omega\setminus\Sigma_{9\varepsilon})}
+ C\|u_0\|_{H^2(\Omega)}.
\end{aligned}
\end{equation*}

Finally, we handle the terms of $I_3,I_4$ in $\eqref{f:3.83}$.
Note that compared with $\mathcal{J}$, $\mathcal{M}$ and $\mathcal{N}$, the structure of $\mathcal{I}$ is most complicated.
Hence, we only write down the whole proof for $\mathcal{I}$. In view of $\eqref{f:3.59}$ and
each of the first inequalities in $\eqref{f:3.60}$, $\eqref{f:3.61}$ and $\eqref{f:3.62}$, we obtain
\begin{equation}\label{f:3.86}
\begin{aligned}
\int_\Omega \big|\mathcal{I}\cdot\nabla\phi_\varepsilon\big| dx
&\leq \varepsilon^{-1}\int_\Omega\Big\{
\big|\varpi_\varepsilon S_\varepsilon(\nabla\psi_{4\varepsilon}\myu{u})\big|
+\big|\varpi_\varepsilon S_\varepsilon(\nabla\psi_{4\varepsilon}\nabla \myu{u})\big|\Big\} |\nabla\phi_\varepsilon| dx \\
& \qquad\qquad\qquad+ \int_\Omega\Big\{\big|\varpi_\varepsilon S_\varepsilon(\psi_{4\varepsilon}\myu{u})\big|
+2\big|\varpi_\varepsilon S_\varepsilon(\psi_{4\varepsilon}\nabla \myu{u})\big|
+ \big|\varpi_\varepsilon S_\varepsilon(\psi_{4\varepsilon}\nabla^2 \myu{u})\big|\Big\} |\nabla\phi_\varepsilon| dx \\
&\leq C\varepsilon^{-1}\|u_0\|_{H^1(\Omega\setminus\Sigma_{8\varepsilon})}\|\nabla\phi_\varepsilon\|_{L^2(\Omega\setminus\Sigma_{9\varepsilon})}
+ C\|u_0\|_{H^2(\Omega)}\|\nabla\phi_\varepsilon\|_{L^2(\Omega)},
\end{aligned}
\end{equation}
where we use the estimate $\eqref{pri:2.5}$ in the last inequality, and the fact that
$S_\varepsilon(\nabla\psi_{4\varepsilon}\nabla\myu{u})$ and $S_\varepsilon(\nabla\psi_{4\varepsilon}\myu{u})$ are supported in
$\Omega\setminus\Sigma_{9\varepsilon}$.
By the same token, we can show the following estimates without real difficulties,
\begin{equation}\label{f:3.87}
\begin{aligned}
&\Big|\int_\Omega\mathcal{J}\cdot\nabla\phi_\varepsilon dx\Big|
\leq C\|u_0\|_{H^1(\Omega)}\|\nabla\phi_\varepsilon\|_{L^2(\Omega)},\\
&\Big|\int_\Omega \mathcal{M}\phi_\varepsilon dx\Big|
\leq C\|u_0\|_{H^2(\Omega)}\|\phi_\varepsilon\|_{L^2(\Omega)}
+ C\varepsilon^{-1}\|u_0\|_{H^1(\Omega\setminus\Sigma_{8\varepsilon})}\|\phi_\varepsilon\|_{L^2(\Omega\setminus\Sigma_{9\varepsilon})},\\
&\Big|\int_\Omega \mathcal{N}\phi_\varepsilon dx\Big|
\leq C\|u_0\|_{H^1(\Omega)}\|\phi_\varepsilon\|_{L^2(\Omega)}.
\end{aligned}
\end{equation}
Combining $\eqref{f:3.86}$ and $\eqref{f:3.87}$ leads to
\begin{equation}\label{f:3.88}
 I_3+I_4 \leq C\|u_0\|_{H^1(\Omega\setminus\Sigma_{8\varepsilon})}
 \|\phi_\varepsilon\|_{H^1(\Omega\setminus\Sigma_{9\varepsilon})}
+ C\varepsilon\|u_0\|_{H^2(\Omega)}\|\phi_\varepsilon\|_{H^1(\Omega)}.
\end{equation}

Inserting the estimates $\eqref{f:3.84}$, $\eqref{f:3.85}$ and $\eqref{f:3.88}$ into
$\eqref{f:3.83}$, we derive
\begin{equation*}
\begin{aligned}
\Big|\int_{\Omega}w_\varepsilon\Phi dx\Big|
\leq C\Big\{\|u_0\|_{H^1(\Omega\setminus\Sigma_{8\varepsilon})}\|\phi_\varepsilon\|_{H^1(\Omega\setminus\Sigma_{9\varepsilon})}
+ \varepsilon\|u_0\|_{H^2(\Omega)}\|\phi_\varepsilon\|_{H^1(\Omega)}\Big\},
\end{aligned}
\end{equation*}
where $C$ depends on $\mu,\kappa,m,d$ and $\Omega$.
We have completed the proof.
\qed
\end{pf}

\section{Dirichlet problem}

\begin{thm}\label{thm:3.1}
 Suppose that the coefficients of $\mathcal{L}_\varepsilon$ satisfy $\eqref{a:1}$, $\eqref{a:2}$ and $\eqref{a:3}$.
 Let $u_\varepsilon$ and $u_0$ be the solutions to $\eqref{pde:1.2}$ with $F\in L^2(\Omega;\mathbb{R}^m)$
 and $g\in H^1(\partial\Omega;\mathbb{R}^m)$. If we additionally assume $A=A^*$, then we have
 \begin{equation}\label{pri:3.5}
 \big\|u_\varepsilon -u_0 - \varepsilon\chi_{0,\varepsilon}S_\varepsilon(\psi_{4\varepsilon}u_0)
 -\varepsilon\chi_{k,\varepsilon}S_\varepsilon(\psi_{4\varepsilon}\nabla_k u_0)\big\|_{H^1_0(\Omega)}
 \leq C\varepsilon^{\frac{1}{2}}\big\{\|F\|_{L^2(\Omega)} + \|g\|_{H^1(\partial\Omega)}\big\},
 \end{equation}
 where $C$ depends only on $\mu,\kappa,m,d$ and $\Omega$.
\end{thm}

\begin{lemma}\label{lemma:3.5}
Suppose that $\mathcal{L}_0$ is the homogenized operator of $\mathcal{L}_\varepsilon$ under the same conditions as in Theorem $\ref{thm:3.1}$.
Let $u_0$ be the solution to $(\mathbf{DH})_0$ in $\eqref{pde:1.2}$ with $F\in L^{2}(\Omega;\mathbb{R}^m)$
and $g\in H^1(\partial\Omega;\mathbb{R}^m)$. Then we have
\begin{equation}\label{pri:3.6}
 \|u_0\|_{H^1(\Omega\setminus\Sigma_{p_1\varepsilon})} \leq C\varepsilon^{\frac{1}{2}}\big\{\|F\|_{L^{2}(\Omega)}
 + \|g\|_{H^1(\partial\Omega)}\big\},
\end{equation}
and
\begin{equation}\label{pri:3.7}
 \|\nabla^2 u_0\|_{L^2(\Sigma_{p_2\varepsilon})} \leq C\varepsilon^{-\frac{1}{2}}\big\{\|F\|_{L^2(\Omega)}
 + \|g\|_{H^1(\partial\Omega)}\big\},
\end{equation}
where $p_1,p_2>0$ are fixed real number, and $C$ depends on $\mu,\kappa,m,d,p_1,p_2$ and $\Omega$.
\end{lemma}

\begin{remark}\label{re:3.1}
\emph{The results of $\eqref{pri:3.6}$ and $\eqref{pri:3.7}$ were originally
established by Z. Shen in \cite[Theorem 2.6]{SZW2} for $L_\varepsilon$ with the Dirichlet boundary condition.
Here we employ the radial maximal function to extend his results to our cases,
and it is convenient to let $p_1,p_2$ be two proper integers.}
\end{remark}
\begin{remark}\label{re:3.2}
\emph{For ease of presentation, we call the estimate of $\|u_0\|_{H^1(\Omega\setminus\Sigma_{p_1\varepsilon})}$ the ``layer type estimate'',
while the estimate of $\|\nabla^2 u_0\|_{L^2(\Sigma_{p_2\varepsilon})}$ is regarded as ``co-layer type estimate'',
where ``co-layer'' means of the complementary layer for short. Obviously,
the name of the two terms are based on the different parts of the domain.}
\end{remark}

\begin{pf}
By setting $L_0 = -\text{div}(\widehat{A}\nabla)$, we can rewrite $(\mathbf{DH})_0$ as
\begin{equation*}
 L_0(u_0) = (\widehat{V} -\widehat{B})\nabla u_0 - (\widehat{c}+\lambda I)u_0 + F \quad \text{in}~\Omega,
 \qquad u_0 = g \quad \text{on}~\partial\Omega.
\end{equation*}
Then we consider $u_0 = v+w$, and they satisfy
\begin{equation}\label{pde:3.3}
 (1)~ L_0(v) = \breve{F} \quad \text{in}~\mathbb{R}^d,
 \qquad\quad
 (2)~ \left\{\begin{aligned}
 L_0(w) &= 0  &\quad&\text{in}~~\Omega,\\
     w  &= g-v &\quad&\text{on}~\partial\Omega,
 \end{aligned}\right.
\end{equation}
where $\breve{F} = (\widehat{V} -\widehat{B})\nabla u_0 - (\widehat{c}+\lambda I)u_0 + F $ in $\Omega$ and $\breve{F} = 0$
on $\mathbb{R}^d\setminus\Omega$.

We first study $(1)$.
Let $\Gamma_0$ denote the fundamental solution of $L_0$, then we have
$v = \Gamma_0\ast\tilde{F}$ in $\mathbb{R}^d$. Moreover, it follows from the Calder\'on-Zygmund theorem (see \cite[Theorem 7.22]{MGLM}) that
\begin{equation}\label{f:3.27}
\begin{aligned}
\|\nabla^2 v\|_{L^{p}(\mathbb{R}^d)}
&\leq C\|\breve{F}\|_{L^{p}(\mathbb{R}^d)} \\
&\leq C\big\{\|F\|_{L^{p}(\Omega)}+\|u_0\|_{H^1(\Omega)}\big\}
\leq C\big\{\|F\|_{L^{p}(\Omega)}+\|g\|_{H^{1/2}(\partial\Omega)}\big\},
\end{aligned}
\end{equation}
where $p\in [\frac{2d}{d+2},2]$, and we use the estimate $\eqref{pri:2.3}$ and H\"older's inequality in the last inequality.
In view of
\begin{equation*}
|\nabla v(x)|\leq C\int_{\mathbb{R}^d}\frac{|\breve{F}(y)|}{|x-y|^{d-1}}dy,
\end{equation*}
and the Hardy-Littlewood-Sobolev inequality on fractional integration (see \cite[Theorem 7.25]{MGLM}), we have
\begin{equation}\label{f:3.28}
\begin{aligned}
\|\nabla v\|_{L^2(\mathbb{R}^d)}
&\leq C\|\breve{F}\|_{L^{\frac{2d}{d+2}}(\mathbb{R}^d)} \\
& \leq C\big\{\|F\|_{L^{\frac{2d}{d+2}}(\Omega)}+\|u_0\|_{H^1(\Omega)}\big\}
\leq C\big\{\|F\|_{L^{\frac{2d}{d+2}}(\Omega)}+\|g\|_{H^{1/2}(\partial\Omega)}\big\}.
\end{aligned}
\end{equation}
Then by Sobolev's inequality, we have
\begin{equation}\label{f:3.29}
\begin{aligned}
 &\|v\|_{L^{\frac{2d}{d-2}}(\mathbb{R}^d)}\leq C\|\nabla v\|_{L^2(\mathbb{R}^d)}
 \leq C \big\{\|F\|_{L^{\frac{2d}{d+2}}(\Omega)} + \|g\|_{H^{1/2}(\partial\Omega)} \big\}, \\
 &\|\nabla v\|_{L^{\frac{2d}{d-1}}(\mathbb{R}^d)}
 \leq C \|\nabla^2 v\|_{L^{\frac{2d}{d+1}}(\mathbb{R}^d)}
 \leq C \big\{\|F\|_{L^{\frac{2d}{d+1}}(\Omega)}+\|g\|_{H^{1/2}(\partial\Omega)}\big\}.
\end{aligned}
\end{equation}
Set $\varrho = (\varrho_1,\cdots,\varrho_d)\in C_0^1(\mathbb{R}^d;\mathbb{R}^d)$ be a
vector field such that $(\varrho,n) \geq c >0$ on $\partial\Omega$, where $n$ denotes the outward unit normal vector to $\partial\Omega$.
Since the divergence theorem, we have
\begin{equation}\label{f:3.94}
\begin{aligned}
c\int_{\partial\Omega} \big(|\nabla v|^2 + |v|^2\big) dS
&\leq \int_{\partial\Omega}<\varrho,n> (|\nabla v|^2 + |v|^2)dS\\
&= \int_\Omega \text{div}(\varrho)\big(|\nabla v|^2+|v|^2\big) dx
+ 2\int_\Omega \varrho\big(\nabla^2 v \nabla v + \nabla v v\big)dx \\
&\leq C\big\{\|v\|_{L^{\frac{2d}{d-2}}(\Omega)}^2 +\|\nabla v\|_{L^2(\Omega)}^2
+ \|\nabla v\|_{L^{\frac{2d}{d-1}}(\Omega)}\|\nabla^2 v\|_{L^{\frac{2d}{d+1}}(\Omega)}\big\} \\
&\leq C\big\{\|F\|_{L^{\frac{2d}{d+1}}(\Omega)}^2+\|g\|_{H^{1/2}(\partial\Omega)}^2\big\}.
\end{aligned}
\end{equation}
where we use H\"older's inequality and Young's inequality in the second inequality,
and the estimates $\eqref{f:3.27}$, $\eqref{f:3.28}$ and $\eqref{f:3.29}$ in the last one.
Furthermore, it is not hard to see that $(\varrho,n) \geq c/2 >0$ on $S_t$ for any $t\in[0,p_1\varepsilon]$,
and therefore the constant $C$ in the right-hand side of
\begin{equation}\label{f:3.30}
\begin{aligned}
\int_{S_t} \big(|\nabla v|^2 + |v|^2\big) dS
\leq C\big\{\|F\|_{L^{\frac{2d}{d+1}}(\Omega)}^2+\|g\|_{H^{1/2}(\partial\Omega)}^2\big\}
\end{aligned}
\end{equation}
is independent of $t$. It follows from $\eqref{pri:2.8}$ and $\eqref{f:3.30}$ that
\begin{equation}\label{f:3.31}
\begin{aligned}
\|v\|_{H^1(\Omega\setminus\Sigma_{p_1\varepsilon})}
=\Big(\int_0^{p_1\varepsilon}\int_{S_t}\big(|\nabla v|^2 + |v|^2\big) dS_tdt\Big)^{\frac{1}{2}}
\leq C\varepsilon^{\frac{1}{2}}\big\{\|F\|_{L^{\frac{2d}{d+1}}(\Omega)}+\|g\|_{H^{1}(\partial\Omega)}\big\}
\end{aligned}
\end{equation}
where we use the fact of $H^1(\partial\Omega;\mathbb{R}^m)\subset H^{1/2}(\partial\Omega;\mathbb{R}^m)$
in the last inequality.

The next thing is to estimate $\|w\|_{H^1(\Omega\setminus\Sigma_{p_1\varepsilon})}$.
According to (2) in $\eqref{pde:3.3}$,
it follows from the $L^2$ regularity problem (see \cite[Theorem 1.4]{SZW11}) that
\begin{equation}\label{f:3.32}
 \|(\nabla w)^*\|_{L^2(\partial\Omega)}
 \leq C\big\{\|g\|_{H^1(\partial\Omega)} + \|v\|_{H^1(\partial\Omega)}\big\}
 \leq C\big\{\|F\|_{L^{\frac{2d}{d+1}}(\Omega)}+\|g\|_{H^1(\partial\Omega)} \big\},
\end{equation}
where we use the estimate $\eqref{f:3.30}$ (for $t=0$) in the last inequality. In view of Lemma $\ref{lemma:2.5}$, we have
\begin{equation}\label{f:3.33}
\begin{aligned}
 \|\mathcal{M}(w)\|
 &\leq C\|w\|_{H^1(\Omega\setminus\Sigma_{c_0})} \\
 &\leq C\|w\|_{H^1(\Omega)}
  \leq C\big\{\|g\|_{H^{1/2}(\partial\Omega)} + \|v\|_{H^{1/2}(\partial\Omega)}\big\} \\
 & \leq C\big\{\|g\|_{H^{1/2}(\partial\Omega)} + \|v\|_{H^{1}(\Omega)}\big\}
 \leq C\big\{ \|F\|_{L^{\frac{2d}{d+2}}(\Omega)}+\|g\|_{H^{1/2}(\partial\Omega)}\big\},
\end{aligned}
\end{equation}
where we use the trace theorem in fourth inequality, and the estimate $\eqref{pri:2.3}$
is employed in the third and the last inequalities.
Similarly, it follows from $\eqref{pri:2.8}$, $\eqref{f:3.32}$ and $\eqref{f:3.33}$ that
\begin{equation}\label{f:3.34}
\begin{aligned}
\|w\|_{H^1(\Omega\setminus\Sigma_{p_1\varepsilon})}
\leq C\varepsilon^{\frac{1}{2}}
\big\{\|(\nabla w)^*\|_{L^2(\partial\Omega)} + \|\mathcal{M}(w)\|_{L^2(\partial\Omega)}\big\}
\leq C\varepsilon^{\frac{1}{2}}\big\{\|F\|_{L^{\frac{2d}{d+1}}(\Omega)}+\|g\|_{H^{1}(\partial\Omega)}\big\}.
\end{aligned}
\end{equation}

Thus combining $\eqref{f:3.31}$ and $\eqref{f:3.34}$,  we have
\begin{equation}\label{pri:3.11}
 \|u_0\|_{H^1(\Omega\setminus\Sigma_{p_1\varepsilon})} \leq C\varepsilon^{\frac{1}{2}}\big\{\|F\|_{L^{\frac{2d}{d+1}}(\Omega)}
 + \|g\|_{H^1(\partial\Omega)}\big\},
\end{equation}
and the desired estimate $\eqref{pri:3.6}$ is established by using H\"older's inequality.

We now turn to prove the estimate $\eqref{pri:3.7}$. Recalling $u_0 = v+w$, it is clear to see that
\begin{equation}\label{f:3.35}
\begin{aligned}
 \|\nabla^2 u_0\|_{L^2(\Sigma_{p_2\varepsilon})}
 &\leq \|\nabla^2 v\|_{L^2(\mathbb{R}^d)}  + \|\nabla^2 w\|_{L^2(\Sigma_{p_2\varepsilon})} \\
 &\leq C\big\{\|F\|_{L^2(\Omega)} + \|g\|_{H^{1/2}(\partial\Omega)}\big\} + \|\nabla^2 w\|_{L^2(\Sigma_{p_2\varepsilon})},
\end{aligned}
\end{equation}
where we use the estimate $\eqref{f:3.27}$ (for $p=2$) in the second inequality. The remaining thing is to estimate
the term of $\|\nabla^2 w\|_{L^2(\Sigma_{p_2\varepsilon})}$.
Noting (2) in $\eqref{pde:3.3}$ again, we obtain the interior estimate
 \begin{equation}\label{f:3.36}
  |\nabla^2 w(x)|\leq \frac{C}{\delta(x)}\Big(\dashint_{B(x,\delta(x)/8)}|\nabla w|^2 dy\Big)^{1/2},
 \end{equation}
 where $C$ depends on $\mu,m,d$. We show the explanation below.
 For any $B(P,r)\subset4B\subset\Omega$, we may assume $P=0$ and $r=1$ from the translation and rescaling arguments. Then due to
 the interior $H^k$ regularity theory (see \cite[Theorem 4.11]{MGLM}), we have $\|w\|_{H^k(B)}\leq C(\mu,m,d)\|u\|_{L^2(2B)}$, where
 $H^k(\Omega;\mathbb{R}^m) = W^{k,2}(\Omega;\mathbb{R}^m)$. By the Sobolev embedding theorem (for $2k>d$), we arrive at
 \begin{equation*}
   |w(0)|\leq \|w\|_{L^\infty(B)} \leq C\|u\|_{H^k(B)} \leq C\|u\|_{L^2(2B)}.
 \end{equation*}
 Since $v=\nabla^2_{ij}w$ also satisfies $L_0(v) = 0$ in $\Omega$, it is clear to see $|\nabla^2w(0)|\leq C\|\nabla^2 w\|_{L^2(2B)}$. Thus
 in view of Cacciopolli's inequality (see \cite[Theorem 4.1]{MGLM}), we obtain $|\nabla^2 u(0)|\leq C\|\nabla u\|_{L^2(4B)}$. Let
 $\tilde{w}(x) = w(rx)$ be the scaled function, a routine computation gives
 $|\nabla^2 w(0)|\leq \frac{C}{r}\big(\dashint_{B(0,r)}|\nabla w|^2 dy\big)^{1/2}$, and this implies the estimate $\eqref{f:3.36}$.

 Moreover, in view of $\eqref{f:3.36}$, we have
 \begin{equation*}
  \big|\nabla^2 w(x)\big|^2 \big[\delta(x)\big]^{d+2} \leq C \int_\Omega |\nabla w|^2 dy.
 \end{equation*}
 Integrating by parts with respect to $x$ on $\Sigma_{c_0}$ (where $c_0$ is layer constant), we obtain
 \begin{equation*}
  c_0^{d+2}\int_{\Sigma_{c_0}}\big|\nabla^2 w(x)\big|^2 dx
  \leq \int_{\Sigma_{c_0}}\big|\nabla^2 w(x)\big|^2 \big[\delta(x)\big]^{d+2} dx
  \leq C|\Omega|\int_{\Omega}|\nabla w|^2 dy.
 \end{equation*}
 As a result, we derive the following (interior) estimate
 \begin{equation}\label{f:3.76}
  \|\nabla^2 w\|_{L^2(\Sigma_{c_0})} \leq C(\mu,m,d,c_0,\Omega) \|\nabla w\|_{L^2(\Omega)}.
 \end{equation}
 Clearly, the above estimate gives
 \begin{equation}\label{f:3.77}
 \int_{\Sigma_{p_2\varepsilon}}|\nabla^2w|^2 dx
  = \int_{\Sigma_{p_2\varepsilon}\setminus\Sigma_{c_0}}|\nabla^2w|^2 dx
  + \int_{\Sigma_{c_0}}|\nabla^2w|^2 dx
  \leq \int_{\Sigma_{p_2\varepsilon}\setminus\Sigma_{c_0}}|\nabla^2w|^2 dx + C\|\nabla w\|_{L^2(\Omega)}^2.
 \end{equation}

 Hence, the remaining thing is to estimate the first term in the right-hand side of $\eqref{f:3.77}$.
 It follows from the estimates $\eqref{f:3.36}$ that
 \begin{equation}\label{f:3.73}
 \begin{aligned}
\int_{\Sigma_{p_2\varepsilon}\setminus\Sigma_{c_0}}|\nabla^2w|^2 dx
 &\leq C\int_{\Sigma_{p_2\varepsilon}\setminus\Sigma_{c_0}}\dashint_{B(x,\delta(x)/8)} \frac{|\nabla w(y)|^2}{[\delta(x)]^2} dydx \\
 & \leq C\int_{p_2\varepsilon}^{c_0}\int_{S_t}\frac{|(\nabla w)^*(x^\prime)|^2}{t^2} dS_t(\Lambda_t(x^\prime))dt
  ~~~\Bigg(= C\int_{p_2\varepsilon}^{c_0}\int_{\partial\Omega}\frac{|(\nabla w)^*(x^\prime)|^2}{t^2} |\Lambda_t(x^\prime)| dS(x^\prime)dt\Bigg) \\
 &\leq C\int_{\partial\Omega}\big|(\nabla w)^*(x^\prime)\big|^2dS(x^\prime)\int_{p_2\varepsilon}^{\infty}\frac{dt}{t^2} \\
 &\leq C\varepsilon^{-1}\|(\nabla w)^*\|_{L^2(\partial\Omega)}^2,
  \end{aligned}
 \end{equation}
 where we use co-area formula $\eqref{eq:2.9}$ in the second inequality.
 Note that $x^\prime\in \partial\Omega$ such that $\delta(x) = \text{dist}(x,x^\prime) = t$ and
 $|\nabla w(y)|\leq (\nabla w)^*(x^\prime)$ for any $y\in B(x,\delta(x)/8)$.

 Hence, inserting $\eqref{f:3.73}$ into $\eqref{f:3.77}$ we have
 \begin{equation}\label{f:3.37}
 \begin{aligned}
   \|\nabla^2 w\|_{L^2(\Sigma_{p_2\varepsilon})}
  & \leq C\varepsilon^{-\frac{1}{2}}\|(\nabla w)^*\|_{L^2(\partial\Omega)} + C\|\nabla w\|_{L^2(\Omega)} \\
  & \leq C\varepsilon^{-\frac{1}{2}}\Big\{\|F\|_{L^{\frac{2d}{d+1}}(\partial\Omega)}+\|g\|_{H^1(\partial\Omega)}\Big\},
 \end{aligned}
 \end{equation}
 where we use $\eqref{pri:2.3}$ and $\eqref{f:3.32}$ in the second inequality. Plugging $\eqref{f:3.37}$ back into $\eqref{f:3.35}$, we arrive at
 \begin{equation*}
 \|\nabla^2 u_0\|_{L^2(\Sigma_{p_2\varepsilon})}
 \leq C\varepsilon^{-\frac{1}{2}}\big\{ \|F\|_{L^{2}(\Omega)}+\|g\|_{H^1(\partial\Omega)}\big\},
 \end{equation*}
 where $C$ depends on $\mu,\kappa,m,d,p_2$ and $\Omega$, and the proof is complete.
\qed
\end{pf}

\begin{flushleft}
\textbf{Proof of Theorem \ref{thm:3.1}}\textbf{.}
Let $\varphi_0 = \mys{S}(\psi_{4\varepsilon}u_0)$ and $\varphi_k = \mys{S}(\psi_{4\varepsilon}\nabla_k u_0)$ in $\eqref{eq:2.12}$. Then we have
\begin{equation*}
 w_\varepsilon = u_\varepsilon - u_0
-\varepsilon\chi_{0,\varepsilon}S_\varepsilon(\psi_{4\varepsilon}u_0)
-\varepsilon\chi_{k,\varepsilon}S_\varepsilon(\psi_{4\varepsilon}\nabla_k u_0).
\end{equation*}
It follows from Lemma $\ref{lemma:3.1}$ that
\begin{equation}\label{f:3.38}
\begin{aligned}
 \|w_\varepsilon\|_{H^1_0(\Omega)}
 &\leq C \bigg\{ \| u_0 - S_\varepsilon(\psi_{4\varepsilon} u_0)\|_{L^2(\Omega)}
 + \|\nabla u_0 - S_\varepsilon(\psi_{4\varepsilon}\nabla u_0)\|_{L^2(\Omega)}
 + \varepsilon\|\varpi_\varepsilon S_\varepsilon(\psi_{4\varepsilon}u_0)\|_{L^2(\Omega)}  \\
 & + \varepsilon\|\varpi_\varepsilon S_\varepsilon(\psi_{4\varepsilon}\nabla u_0)\|_{L^2(\Omega)}
   + \varepsilon\|\varpi_\varepsilon \nabla S_\varepsilon(\psi_{4\varepsilon}u_0)\|_{L^2(\Omega)}
   + \varepsilon\|\varpi_\varepsilon \nabla S_\varepsilon(\psi_{4\varepsilon}\nabla u_0)\|_{L^2(\Omega)}
 \bigg\}.
\end{aligned}
\end{equation}
We note that $S_\varepsilon(\psi_{4\varepsilon} u_0)$ and $S_\varepsilon(\psi_{4\varepsilon} \nabla u_0)$
is supported in $\Sigma_{3\varepsilon}$. To complete the proof, we need the following estimates.
\end{flushleft}

Due to $\eqref{pri:2.6}$, we have
\begin{equation}\label{f:3.39}
\begin{aligned}
 \| u_0 - S_\varepsilon(\psi_{4\varepsilon} u_0)\|_{L^2(\Omega)}
 &\leq \| (1-\psi_{4\varepsilon})u_0\|_{L^2(\Omega)}
  + \| \psi_{4\varepsilon}u_0 - S_\varepsilon(\psi_{4\varepsilon} u_0)\|_{L^2(\Omega)} \\
 & \leq \|u_0\|_{L^2(\Omega\setminus\Sigma_{8\varepsilon})}
 + C\varepsilon\|\nabla(\psi_{4\varepsilon}u_0)\|_{L^2(\Omega)} \\
 & \leq C\|u_0\|_{L^2(\Omega\setminus\Sigma_{8\varepsilon})} + C\varepsilon\|\nabla u_0\|_{L^2(\Sigma_{4\varepsilon})}
\end{aligned}
\end{equation}
and
\begin{equation}\label{f:3.40}
 \| \nabla u_0 - S_\varepsilon(\psi_{4\varepsilon} \nabla u_0)\|_{L^2(\Omega)}
  \leq C\|\nabla u_0\|_{L^2(\Omega\setminus\Sigma_{8\varepsilon})} + C\varepsilon\|\nabla^2 u_0\|_{L^2(\Sigma_{4\varepsilon})}.
\end{equation}
From $\eqref{pri:2.5}$, it follows that
\begin{equation}\label{f:3.41}
\begin{aligned}
\|\varpi_\varepsilon S_\varepsilon(\psi_{4\varepsilon}u_0)\|_{L^2(\Omega)}
&+ \|\varpi_\varepsilon \nabla S_\varepsilon(\psi_{4\varepsilon}u_0)\|_{L^2(\Omega)} \\
& \leq C\|\psi_{4\varepsilon}u_0\|_{L^2(\Omega)} + C\|\nabla(\psi_{4\varepsilon}u_0)\|_{L^2(\Omega)}  \\
& \leq C\|u_0\|_{L^2(\Sigma_{4\varepsilon})} + C\varepsilon^{-1}\|u_0\|_{L^2(\Omega\setminus\Sigma_{8\varepsilon})}
+ C\|\nabla u_0\|_{L^2(\Sigma_{4\varepsilon})}
\end{aligned}
\end{equation}
and
\begin{equation}\label{f:3.42}
\begin{aligned}
\|\varpi_\varepsilon S_\varepsilon(\psi_{4\varepsilon}\nabla u_0)\|_{L^2(\Omega)}
& +\|\varpi_\varepsilon \nabla S_\varepsilon(\psi_{4\varepsilon}\nabla u_0)\|_{L^2(\Omega)} \\
& \leq C\|\nabla u_0\|_{L^2(\Sigma_{4\varepsilon})} + C\varepsilon^{-1}\|\nabla u_0\|_{L^2(\Omega\setminus\Sigma_{8\varepsilon})}
+ C\|\nabla^2 u_0\|_{L^2(\Sigma_{4\varepsilon})}.
\end{aligned}
\end{equation}
Plugging the estimates $\eqref{f:3.39}$, $\eqref{f:3.40}$, $\eqref{f:3.41}$ and $\eqref{f:3.42}$ into $\eqref{f:3.38}$, we obtain
\begin{equation*}
\|w_\varepsilon\|_{H^1_0(\Omega)}
\leq C\|u_0\|_{H^1(\Omega\setminus\Sigma_{8\varepsilon})} + C\varepsilon\|u_0\|_{H^2(\Sigma_{4\varepsilon})}.
\end{equation*}
Moreover, this together with the estimates $\eqref{pri:3.6}$ and $\eqref{pri:3.7}$ implies
\begin{equation}
 \|w_\varepsilon\|_{H^1_0(\Omega)}
 \leq C\varepsilon^{\frac{1}{2}}\big\{\|F\|_{L^2(\Omega)}+\|g\|_{H^1(\partial\Omega)}\big\},
\end{equation}
where $C$ depends on $\mu,\kappa,m,d$ and $\Omega$, and
we have completed the proof.
\qed

\begin{lemma}[Improved lemma]\label{lemma:3.6}
Assume the same conditions as in Theorem $\ref{thm:3.1}$.
Let $u_0$ be the solution to $(\mathbf{DH})_0$ in $\eqref{pde:1.2}$ with $F\in L^{2}(\Omega;\mathbb{R}^m)$
and $g\in H^1(\partial\Omega;\mathbb{R}^m)$. Then we have
\begin{equation}\label{pri:3.9}
\|u_0\|_{H^1(\Omega\setminus\Sigma_{p_1\varepsilon};\delta)} \leq C\varepsilon\big\{\|F\|_{L^2(\Omega)}
 + \|g\|_{H^1(\partial\Omega)}\big\},
\end{equation}
and
\begin{equation}\label{pri:3.8}
\max\big\{\|\nabla^2 u_0\|_{L^2(\Sigma_{p_2\varepsilon};\delta)},
~\|u_0\|_{H^1(\Sigma_{p_2\varepsilon};\delta^{-1})}\big\} \leq C\big[\ln(c_0/\varepsilon)\big]^{1/2}\big\{\|F\|_{L^2(\Omega)}
 + \|g\|_{H^1(\partial\Omega)}\big\},
\end{equation}
where $p_1,p_2>0$ are fixed real numbers and $c_0$ is the layer constant, and $C$ depends on $\mu,\kappa,m,d,p_1$ and $\Omega$.
\end{lemma}

\begin{remark}
\emph{Recall that the layer constant $c_0$ is defined in Subsection $\ref{sec:1.1}$.
Compared with the results of Lemma $\ref{lemma:3.5}$,
we can see that the weighted-type norms can notably improve the $\varepsilon$' power both
in the layer type estimate and in the co-layer type estimate. Here the weighted function $\delta$
plays a key role, briefly speaking, which can produce a good factor $\varepsilon^{\frac{1}{2}}$.
So we call Lemma $\ref{lemma:3.6}$ the improved lemma.}
\end{remark}

\begin{pf}
The proof of the estimate $\eqref{pri:3.9}$ is straightforward. In view of $\eqref{pri:2.19}$ and $\eqref{pri:3.6}$, we have
\begin{equation*}
 \|u_0\|_{H^1(\Omega\setminus\Sigma_{p_1\varepsilon};\delta)}
 \leq C\varepsilon^{\frac{1}{2}}\|u_0\|_{H^1(\Omega\setminus\Sigma_{p_1\varepsilon})}
 \leq C\varepsilon\big\{\|F\|_{L^2(\Omega)}
 + \|g\|_{H^1(\partial\Omega)}\big\}.
\end{equation*}

We now prove the estimate $\eqref{pri:3.8}$.  Proceeding as in the proof of $\eqref{pri:3.7}$ in Lemma $\ref{lemma:3.5}$, we first have
 \begin{equation}\label{f:3.74}
 \begin{aligned}
 \|\nabla^2 u_0\|_{L^2(\Sigma_{p_2\varepsilon};\delta)}
 &\leq C\|\nabla^2 v\|_{L^2(\mathbb{R}^d)} + \|\nabla^2 w\|_{L^2(\Sigma_{p_2\varepsilon};\delta)} \\
 &\leq C\big\{\|F\|_{L^2(\Omega)} + \|g\|_{H^1(\partial\Omega)}\big\} + \|\nabla^2 w\|_{L^2(\Sigma_{p_2\varepsilon};\delta)},
 \end{aligned}
 \end{equation}
 where we use the hypothesis of $\delta(x) = 0$ when $x\in\mathbb{R}^d\setminus\Omega$
 (see Subsection $\ref{sec:1.1}$) in the first inequality,
 and $\eqref{f:3.27}$ (for $p=2$) in the second one. Due to the estimate $\eqref{f:3.76}$, we have
 \begin{equation}\label{f:3.78}
 \begin{aligned}
   \|\nabla^2 w\|_{L^2(\Sigma_{p_2\varepsilon};\delta)}
  &\leq \|\nabla^2 w\|_{L^2(\Sigma_{p_2\varepsilon}\setminus\Sigma_{c_0};\delta)}
  + C\|\nabla w\|_{L^2(\Omega)}\\
  &\leq \|\nabla^2 w\|_{L^2(\Sigma_{p_2\varepsilon}\setminus\Sigma_{c_0};\delta)}
  + C\big\{\|F\|_{L^2(\Omega)} + \|g\|_{H^{1}(\partial\Omega)}\big\},
 \end{aligned}
 \end{equation}
 where we use $\eqref{pri:2.3}$ in the second inequality. Hence,
 it is sufficient to study the term of $\|\nabla^2 w\|_{L^2(\Sigma_{p_2\varepsilon}\setminus\Sigma_{c_0};\delta)}$.
 By the estimate $\eqref{f:3.36}$,
 \begin{equation}\label{f:3.75}
 \begin{aligned}
 \|\nabla^2 w\|_{L^2(\Sigma_{p_2\varepsilon}\setminus\Sigma_{c_0};\delta)}^2
 &= \int_{\Sigma_{p_2\varepsilon}\setminus\Sigma_{c_0}} |\nabla^2 w|^2 \delta(x)dx
 \leq C\int_{\Sigma_{p_2\varepsilon}\setminus\Sigma_{c_0}}\dashint_{B(x,\delta(x)/8)} \frac{|\nabla w(y)|^2}{\delta(x)} dydx \\
 &\leq C\int_{p_2\varepsilon}^{c_0}\int_{S_t}\frac{|(\nabla w)^*(x^\prime)|^2}{t} dS_t(\Lambda_t(x^\prime))dt
 \leq C\ln(c_0/\varepsilon)\|(\nabla w)^*\|_{L^2(\partial\Omega)} \\
 &\leq C\ln(c_0/\varepsilon)\big\{\|F\|_{L^2(\Omega)}^2 + \|g\|_{H^1(\partial\Omega)}^2\big\},
 \end{aligned}
 \end{equation}
 where the second inequality follows from the same observation as we did in $\eqref{f:3.73}$,
 and we use $\eqref{f:3.32}$ in the last inequality. Then combining $\eqref{f:3.74}$ and $\eqref{f:3.75}$,
 we partially derive the estimate $\eqref{pri:3.8}$.

 The rest task is to estimate $\|u_0\|_{H^1(\Sigma_{p_2\varepsilon};\delta^{-1})}$. By the same idea used above, we have
 \begin{equation}\label{f:3.95}
 \begin{aligned}
  \|u_0\|_{H^1(\Sigma_{p_2\varepsilon};\delta^{-1})}^2
 &\leq \int_{\Sigma_{p_2\varepsilon}\setminus\Sigma_{c_0}}\big(|u_0|^2+|\nabla u_0|^2\big)\delta^{-1}(x) dx
 + c_0^{-1}\int_{\Sigma_{c_0}}\big(|u_0|^2+|\nabla u_0|^2\big)dx \\
 &\leq \int_{p_2\varepsilon}^{c_0}\int_{\partial\Omega}\big(|u_0(\Lambda_t(z))|^2 + |\nabla u_0(\Lambda_t(z))|^2\big)
 |\nabla\Lambda_t|dS(z)t^{-1}dt + C\|u_0\|_{H^1(\Omega)}^2 \\
 &\leq C\int_{\partial\Omega}\big(|\mathcal{M}(u_0)|^2 + |\mathcal{M}(\nabla u_0)|^2\big)dS\int_{p_2\varepsilon}^{c_0}\frac{1}{t}dt
 + C\|u_0\|_{H^1(\Omega)}^2 \\
 &\leq C\ln(c_0/\varepsilon)\Big\{\|\mathcal{M}(u_0)\|_{L^2(\partial\Omega)}^2 + \|\mathcal{M}(\nabla u_0)\|_{L^2(\partial\Omega)}^2
 + \|u_0\|_{H^1(\Omega)}^2\Big\},
 \end{aligned}
 \end{equation}
 where we use co-area formula $\eqref{eq:2.9}$ in the second inequality,
 and the estimate $\eqref{pri:2.8}$ in the third one,
 and the estimate $\eqref{pri:2.9}$ in the last one.
 This implies
 \begin{equation}
   \|u_0\|_{H^1(\Sigma_{p_2\varepsilon};\delta^{-1})}
   \leq C[\ln(c_0/\varepsilon)]^{1/2}\big\{\|F\|_{L^2(\Omega)}+\|g\|_{H^1(\partial\Omega)}\big\},
 \end{equation}
 where we use the estimates $\eqref{pri:2.3}$ $\eqref{f:3.32}$ and $\eqref{f:3.33}$. The proof is complete.
 \qed

\end{pf}

\begin{thm}\label{thm:3.2}
Assume the same conditions as in Theorem $\ref{thm:3.1}$.
Let $u_\varepsilon,u_0$ be the solutions to $\eqref{pde:1.2}$ with $F\in L^2(\Omega;\mathbb{R}^m)$ and $g\in H^1(\partial\Omega;\mathbb{R}^m)$.
Then we have
\begin{equation}\label{pri:3.10}
 \big\|u_\varepsilon -u_0 - \varepsilon\chi_{0,\varepsilon}S_\varepsilon(\psi_{4\varepsilon}u_0)
 -\varepsilon\chi_{k,\varepsilon}S_\varepsilon(\psi_{4\varepsilon}\nabla_k u_0)\big\|_{L^2(\Omega)}
 \leq C\varepsilon\ln(c_0/\varepsilon)\big\{\|F\|_{L^2(\Omega)} + \|g\|_{H^1(\partial\Omega)}\big\},
\end{equation}
where $c_0$ is the layer constant, and $C$ depends on $\mu,\kappa,m,d,c_0$ and $\Omega$.
\end{thm}

\begin{pf}
We prove this theorem by a duality argument.
For any $\Phi\in L^2(\Omega;\mathbb{R}^m)$,
we find two weak solutions $\phi_\varepsilon$ and $\phi_0$ solving $\eqref{pde:2.6}$.
By recalling
\begin{equation*}
\xi_\varepsilon = \phi_\varepsilon - \phi_0 - \varepsilon\chi_{0,\varepsilon}^*S_\varepsilon(\psi_{10\varepsilon}\phi_0)
+ \varepsilon\chi_{k,\varepsilon}^*S_\varepsilon(\psi_{10\varepsilon}\nabla_k\phi_0)
\end{equation*}
in $\eqref{eq:2.14}$, it straightforward follows from Theorem $\ref{thm:3.1}$ that
\begin{equation}
 \|\xi_\varepsilon\|_{H^1_0(\Omega)} \leq C\varepsilon^{\frac{1}{2}}\|\Phi\|_{L^2(\Omega)}.
\end{equation}

Due to the linearity of $\mathcal{L}_\varepsilon$,
it is convenient to assume $\|F\|_{L^2(\Omega)} + \|g\|_{H^1(\partial\Omega)} = 1$,
otherwise we replace $u_\varepsilon$ and $u_0$ by $u_\varepsilon/(\|F\|_{L^2(\Omega)} + \|g\|_{H^1(\partial\Omega)})$ and
$u_0/(\|F\|_{L^2(\Omega)} + \|g\|_{H^1(\partial\Omega)})$, respectively. Hence, by setting
\begin{equation*}
  w_\varepsilon = u_\varepsilon -u_0 - \varepsilon\chi_{0,\varepsilon}S_\varepsilon(\psi_{4\varepsilon}u_0)
 -\varepsilon\chi_{k,\varepsilon}S_\varepsilon(\psi_{4\varepsilon}\nabla_k u_0),
\end{equation*}
it is equivalent to proving $\|w_\varepsilon\|_{L^2(\Omega)}\leq C\varepsilon\ln(r_0/\varepsilon)$. To do so,
in view of Lemma $\ref{lemma:2.10}$, we have
\begin{equation}\label{f:3.82}
\begin{aligned}
\Big|\int_\Omega w_\varepsilon\Phi dx\Big|
&\leq C\|u_0\|_{H^1(\Omega\setminus\Sigma_{9\varepsilon})}\|\phi_0\|_{H^1(\Omega\setminus\Sigma_{9\varepsilon})}
 + C\varepsilon\|u_0\|_{H^1(\Omega)}\|\phi_\varepsilon\|_{H^1(\Omega)} \\
& \qquad + C\Big\{\| u_0\|_{H^1(\Omega\setminus\Sigma_{8\varepsilon};\delta)}
+\varepsilon \|\nabla u_0\|_{L^2(\Sigma_{4\varepsilon};\delta)}
+\varepsilon \|\nabla^2 u_0\|_{L^2(\Sigma_{4\varepsilon};\delta)}\Big\}\|\phi_0\|_{H^1(\Sigma_{4\varepsilon};\delta^{-1})} \\
& \qquad\qquad\qquad+ C\Big\{\| u_0\|_{H^1(\Omega\setminus\Sigma_{8\varepsilon})}
+\varepsilon \|\nabla u_0\|_{L^2(\Sigma_{4\varepsilon})}
+\varepsilon \|\nabla^2 u_0\|_{L^2(\Sigma_{4\varepsilon})}\Big\}\\
&\qquad\qquad\qquad\qquad\qquad\cdot\Big\{\|\xi_\varepsilon\|_{H^1(\Omega)} + \|\phi_0\|_{H^1(\Omega\setminus\Sigma_{20\varepsilon})}
+ \varepsilon\|\phi_0\|_{H^1(\Omega)}+ \varepsilon\|\nabla^2\phi_0\|_{L^2(\Sigma_{10\varepsilon})}\Big\}.
\end{aligned}
\end{equation}
To complete the proof, we need the following computations:
\begin{equation}\label{f:3.79}
\begin{aligned}
\|u_0\|_{H^1(\Omega\setminus\Sigma_{9\varepsilon})}\|\phi_0\|_{H^1(\Omega\setminus\Sigma_{9\varepsilon})}
 &+ C\varepsilon\|u_0\|_{H^1(\Omega)}\|\phi_\varepsilon\|_{H^1(\Omega)}\\
 &\leq C\varepsilon^{\frac{1}{2}}\cdot C\varepsilon^{\frac{1}{2}}\|\Phi\|_{L^2(\Omega)} + C\varepsilon\|\Phi\|_{L^2(\Omega)}
 \leq C\varepsilon\|\Phi\|_{L^2(\Omega)},
\end{aligned}
\end{equation}
where we use the estimates $\eqref{pri:3.6}$ and $\eqref{pri:2.3}$ in the first inequality.
\begin{equation}\label{f:3.80}
\begin{aligned}
\Big\{\| u_0\|_{H^1(\Omega\setminus\Sigma_{8\varepsilon};\delta)}
+\varepsilon \|\nabla u_0\|_{L^2(\Sigma_{4\varepsilon};\delta)}
&+\varepsilon \|\nabla^2 u_0\|_{L^2(\Sigma_{4\varepsilon};\delta)}\Big\}\|\phi_0\|_{H^1(\Sigma_{4\varepsilon};\delta^{-1})} \\
&\leq C\Big\{\varepsilon
+ \varepsilon\big[\ln(c_0/\varepsilon)\big]^{\frac{1}{2}}\Big\}\big[\ln(c_0/\varepsilon)\big]^{\frac{1}{2}}\|\Phi\|_{L^2(\Omega)}
\leq C\varepsilon\ln(c_0/\varepsilon)\|\Phi\|_{L^2(\Omega)},
\end{aligned}
\end{equation}
where we use the estimates $\eqref{pri:2.3}$, $\eqref{pri:3.9}$ and $\eqref{pri:3.8}$ in the first inequality.
\begin{equation}\label{f:3.81}
\begin{aligned}
\Big\{\| u_0\|_{H^1(\Omega\setminus\Sigma_{8\varepsilon})}
&+\varepsilon \|\nabla u_0\|_{L^2(\Sigma_{4\varepsilon})}
+\varepsilon \|\nabla^2 u_0\|_{L^2(\Sigma_{4\varepsilon})}\Big\}\\
&\cdot\Big\{\|\xi_\varepsilon\|_{H^1(\Omega)} + \|\phi_0\|_{H^1(\Omega\setminus\Sigma_{20\varepsilon})}
+ \varepsilon\|\phi_0\|_{H^1(\Omega)}+ \varepsilon\|\nabla^2\phi_0\|_{L^2(\Sigma_{10\varepsilon})}\Big\}\\
& \qquad\qquad\qquad\qquad\qquad \leq C\Big\{\varepsilon + \varepsilon + \varepsilon^{1/2}\Big\}
\cdot\Big\{\varepsilon^{1/2}+\varepsilon + \varepsilon + \varepsilon^{1/2}\Big\}\|\Phi\|_{L^2(\Omega)}
\leq C\varepsilon \|\Phi\|_{L^2(\Omega)},
\end{aligned}
\end{equation}
where we use the estimates $\eqref{pri:2.3}$, $\eqref{pri:3.6}$ and $\eqref{pri:3.7}$ in the first inequality.

 Consequently, plugging the estimates $\eqref{f:3.79}-\eqref{f:3.81}$ back into $\eqref{f:3.82}$, we have
 \begin{equation*}
 \Big|\int_{\Omega}w_\varepsilon\Phi dx\Big|
 \leq C\varepsilon\ln(c_0/\varepsilon)\|\Phi\|_{L^2(\Omega)},
 \end{equation*}
 and this implies the desired estimate $\eqref{pri:3.10}$, and we have completed the proof.
\qed

\end{pf}

\begin{thm}\label{thm:3.3}
Suppose that the coefficients of $\myl{L}{\varepsilon}$ satisfy $\eqref{a:1}$, $\eqref{a:2}$ and $\eqref{a:3}$,
and $A$ additionally satisfies $A=A^*$. Let $u_\varepsilon,u_0$ be the weak solutions to $\eqref{pde:1.2}$
with $F\in L^p(\Omega;\mathbb{R}^m)$ and $g\in H^1(\partial\Omega;\mathbb{R}^m)$, where $q=\frac{2d}{d+1}$. Then
\begin{equation}\label{pri:3.2}
\big\|u_\varepsilon - u_0 -\varepsilon\chi_{0,\varepsilon}S^2_\varepsilon(\psi_{2\varepsilon}u_0)
-\varepsilon\chi_{k,\varepsilon}S^2_\varepsilon(\psi_{2\varepsilon}\nabla_k u_0)\big\|_{H^1_0(\Omega)}
\leq C\varepsilon^{\frac{1}{2}}\big\{\|F\|_{L^q(\Omega)}+\|g\|_{H^1(\partial\Omega)}\big\},
\end{equation}
where $C$ depends only on $\mu,\kappa,m,d$ and $\Omega$.
\end{thm}

\begin{pf}
We note that $\|w_\varepsilon\|_{H^1(\Omega)}$  is exactly the left-hand side of $\eqref{pri:3.2}$
by setting $\varphi_0 = \mys{S}^2(\psi_{2\varepsilon}u_0)$ and $\varphi_k = \mys{S}^2(\psi_{2\varepsilon}\nabla_k u_0)$ in $\eqref{eq:2.12}$.
Then it follows from $\eqref{pri:3.1}$ that
\begin{equation}\label{f:3.8}
\begin{aligned}
\|w_{\varepsilon}\|_{H^1(\Omega)}
&\leq C\Big\{
\|\nabla u_0 - \mys{S}^2(\psi_{2\varepsilon}\nabla u_0)\|_{L^2(\Omega)}
+\|u_0 - \mys{S}^2(\psi_{2\varepsilon}u_0)\|_{L^2(\Omega)}
+\varepsilon\|\varpi_{\varepsilon}\nabla\mys{S}^2(\psi_{2\varepsilon} u_0)\|_{L^2(\Omega)}\\
&+ \varepsilon \|\varpi_{\varepsilon}\nabla\mys{S}^2(\psi_{2\varepsilon}\nabla u_0)\|_{L^2(\Omega)}
+ \varepsilon \|\varpi_{\varepsilon}\mys{S}^2(\psi_{2\varepsilon}u_0)\|_{L^2(\Omega)}
+\varepsilon \|\varpi_{\varepsilon}\mys{S}^2(\psi_{2\varepsilon}\nabla u_0)\|_{L^2(\Omega)}\Big\}.
\end{aligned}
\end{equation}
Before proceeding further, let us do some calculations:
\begin{equation}\label{f:3.9}
\begin{aligned}
\|\nabla u_0 - \mys{S}^2(\psi_{2\varepsilon}\nabla u_0)\|_{L^2(\Omega)}
&\leq\|(1-\psi_{2\varepsilon})\nabla u_0\|_{L^2(\Omega)}
+ \|\psi_{2\varepsilon}\nabla u_0 - \mys{S}(\psi_{2\varepsilon}\nabla u_0)\|_{L^2(\Omega)} \\
& +\big\|\mys{S}\big[\psi_{2\varepsilon}\nabla u_0 - \mys{S}(\psi_{2\varepsilon}\nabla u_0)\big]\|_{L^2(\Omega)} \\
&\leq \|\nabla u_0\|_{L^2(\Omega\setminus\Sigma_{4\varepsilon})}
+ C\|\psi_{2\varepsilon}\nabla u_0 - \mys{S}(\psi_{2\varepsilon}\nabla u_0)\|_{L^2(\Omega)},
\end{aligned}
\end{equation}
\begin{equation}\label{f:3.10}
\begin{aligned}
\|\varpi_{\varepsilon}\nabla\mys{S}^2(\psi_{2\varepsilon}\nabla u_0)\|_{L^2(\Omega)}
&\leq C\|\nabla\mys{S}(\psi_{2\varepsilon}\nabla u_0)\|_{L^2(\Omega)} \\
&\leq C\big\{\varepsilon^{-1}\|\nabla u_0\|_{L^2(\Omega\setminus\Sigma_{4\varepsilon})}
+ \|\mys{S}(\psi_{2\varepsilon}\nabla^2 u_0)\|_{L^2(\Omega)}\big\},
\end{aligned}
\end{equation}
and
\begin{equation}\label{f:3.11}
\|\varpi_{\varepsilon}\mys{S}^2(\psi_{2\varepsilon}\nabla u_0)\|_{L^2(\Omega)}
\leq C\|\mys{S}(\psi_{2\varepsilon}\nabla u_0)\|_{L^2(\Omega)}
\end{equation}
where we mainly use the estimate $\eqref{pri:2.5}$ in the second inequality of $\eqref{f:3.9}$, and in the first inequality
of $\eqref{f:3.10}$,
as well as in $\eqref{f:3.11}$. After a similar computation, we have
\begin{equation}\label{f:3.12}
\begin{aligned}
&\|\varpi_{\varepsilon}\mys{S}^2(\psi_{2\varepsilon} u_0)\|_{L^2(\Omega)}
\leq C\|\mys{S}(\psi_{2\varepsilon}u_0)\|_{L^2(\Omega)}, \\
& \| u_0 - \mys{S}^2(\psi_{2\varepsilon} u_0)\|_{L^2(\Omega)}
\leq \|u_0\|_{L^2(\Omega\setminus\Sigma_{4\varepsilon})}
+ C\|\psi_{2\varepsilon} u_0 - \mys{S}(\psi_{2\varepsilon} u_0)\|_{L^2(\Omega)}, \\
&\|\varpi_{\varepsilon}\nabla\mys{S}^2(\psi_{2\varepsilon} u_0)\|_{L^2(\Omega)}
\leq C\Big\{\frac{1}{\varepsilon}\|u_0\|_{L^2(\Omega\setminus\Sigma_{4\varepsilon})}
+ \|\mys{S}(\psi_{2\varepsilon}\nabla u_0)\|_{L^2(\Omega)}\Big\}. \\
\end{aligned}
\end{equation}
By substituting $\eqref{f:3.9}$, $\eqref{f:3.10}$, $\eqref{f:3.11}$ and $\eqref{f:3.12}$ into $\eqref{f:3.8}$, we find
\begin{equation}\label{f:3.22}
\begin{aligned}
\|w_\varepsilon\|_{H^1(\Omega)}& \leq C\Big\{
\|u_0\|_{H^1(\Omega\setminus\Sigma_{4\varepsilon})}
+\|\psi_{2\varepsilon}\nabla u_0 - \mys{S}(\psi_{2\varepsilon}\nabla u_0)\|_{L^2(\Omega)}
+\|\psi_{2\varepsilon}u_0 - \mys{S}(\psi_{2\varepsilon}u_0)\|_{L^2(\Omega)}\Big\} \\
& + C\varepsilon\Big\{\|\mys{S}(\psi_{2\varepsilon}\nabla^2 u_0)\|_{L^2(\Omega)}
+ \|\mys{S}(\psi_{2\varepsilon}\nabla u_0)\|_{L^2(\Omega)}
+\|\mys{S}(\psi_{2\varepsilon}u_0)\|_{L^2(\Omega)} \Big\}.
\end{aligned}
\end{equation}

We now estimate $\|\psi_{2\varepsilon}\nabla u_0 - \mys{S}(\psi_{2\varepsilon}\nabla u_0)\|_{L^2(\Omega)}$.
Noting that $\psi_{2\varepsilon}\nabla u_0 - \mys{S}(\psi_{2\varepsilon}\nabla u_0)$ is supported in $\Sigma_\varepsilon$, it is
equivalent to estimating $\|\psi_{2\varepsilon}\nabla u_0 - \mys{S}(\psi_{2\varepsilon}\nabla u_0)\|_{L^2(\Sigma_\varepsilon)}$. Hence,
\begin{equation*}
\begin{aligned}
&\|\psi_{2\varepsilon}\nabla u_0 - \mys{S}(\psi_{2\varepsilon}\nabla u_0)\|_{L^2(\Sigma_\varepsilon)} \\
&\leq \|\nabla v - \mys{S}(\nabla v)\|_{L^2(\Sigma_\varepsilon)}
+\|(\psi_{2\varepsilon}-1)\nabla v\|_{L^2(\Sigma_\varepsilon)}
 +\|S_\varepsilon\big((\psi_{2\varepsilon}-1)\nabla v\big)\|_{L^2(\Sigma_\varepsilon)}
+ \|\psi_{2\varepsilon}\nabla w - \mys{S}(\psi_{2\varepsilon}\nabla w)\|_{L^2(\Sigma_\varepsilon)} \\
&\leq \|\nabla v - \mys{S}(\nabla v)\|_{L^2(\mathbb{R}^d)}
+2\|(\psi_{2\varepsilon}-1)\nabla v\|_{L^2(\Omega)}
+ \|\psi_{2\varepsilon}\nabla w - \mys{S}(\psi_{2\varepsilon}\nabla w)\|_{L^2(\Omega)}
\end{aligned}
\end{equation*}
where we use the fact of
$\|S_\varepsilon\big((\psi_{2\varepsilon}-1)\nabla v\big)\|_{L^2(\Sigma_\varepsilon)}
\leq C\|(\psi_{2\varepsilon}-1)\nabla v\|_{L^2(\Omega)}$ in the second inequality. This implies
\begin{equation}\label{f:3.23}
\begin{aligned}
&\|\psi_{2\varepsilon}\nabla u_0 - \mys{S}(\psi_{2\varepsilon}\nabla u_0)\|_{L^2(\Sigma_\varepsilon)} \\
& \leq  C\varepsilon^{\frac{1}{2}}\|\nabla^2 v\|_{L^{\frac{2d}{d+1}}(\mathbb{R}^d)}
+ C\|\nabla v\|_{L^2(\Omega\setminus\Sigma_{4\varepsilon})}
+ C\varepsilon\|\nabla(\nabla w\psi_{2\varepsilon})\|_{L^2(\Omega)} \\
& \leq C\varepsilon^{\frac{1}{2}}\|\nabla^2 v\|_{L^{\frac{2d}{d+1}}(\mathbb{R}^d)}
+ C\|\nabla v\|_{L^2(\Omega\setminus\Sigma_{4\varepsilon})}
+ C\big\{\varepsilon\|\nabla^2 w\|_{L^2(\Sigma_{2\varepsilon})}
+ \|\nabla w\|_{L^2(\Omega\setminus\Sigma_{4\varepsilon})}\big\}\\
&\leq C\varepsilon^{\frac{1}{2}}\big\{\|F\|_{L^{\frac{2d}{d+1}}(\Omega)} + \|g\|_{H^1(\partial\Omega)}\big\},
\end{aligned}
\end{equation}
where we use the estimates $\eqref{pri:2.6}$ and $\eqref{pri:2.7}$ in the first inequality,
and the estimates $\eqref{f:3.27}$ (for $p=\frac{2d}{d+1}$), $\eqref{f:3.31}$ (for $t=0$),
$\eqref{f:3.37}$ and $\eqref{f:3.34}$ are employed in the last
inequality.
Also, we have
\begin{equation}\label{f:3.24}
\begin{aligned}
\|\psi_{2\varepsilon} u_0 - \mys{S}(\psi_{2\varepsilon} u_0)\|_{L^2(\Omega)}
& \leq C\varepsilon\|\nabla(\psi_{2\varepsilon}u_0)\|_{L^2(\Omega)} \\
&\leq C\|u_0\|_{L^2(\Omega\setminus\Sigma_{4\varepsilon})} + C\varepsilon\|\nabla u_0\|_{L^2(\Omega)} \\
&\leq  C\varepsilon^{\frac{1}{2}}\big\{\|F\|_{L^{\frac{2d}{d+1}}(\Omega)}+\|g\|_{H^1(\partial\Omega)}\big\},
\end{aligned}
\end{equation}
where we use the estimates $\eqref{pri:3.11}$ and $\eqref{pri:2.3}$ in the last inequality.

We still need to estimate
\begin{equation}\label{f:3.25}
\begin{aligned}
\|S_\varepsilon(\psi_{2\varepsilon}\nabla^2u_0)\|_{L^2(\Omega)}
&\leq \|S_\varepsilon(\psi_{2\varepsilon}\nabla^2v)\|_{L^2(\Omega)}
+ \|S_\varepsilon(\psi_{2\varepsilon}\nabla^2 w)\|_{L^2(\Omega)} \\
& \leq C\varepsilon^{-\frac{1}{2}}\|\psi_{2\varepsilon}\nabla^2 v\|_{L^{\frac{2d}{d+1}}(\mathbb{R}^d)}
+ C\|\nabla^2 w\|_{L^2(\Sigma_{2\varepsilon})}\\
&\leq C\varepsilon^{-\frac{1}{2}}\big\{\|g\|_{H^1(\partial\Omega)}+\|F\|_{L^{\frac{2d}{d+1}}(\Omega)}\big\},
\end{aligned}
\end{equation}
where we use $\eqref{pri:2.7}$ in the second inequality, and $\eqref{f:3.27}$, $\eqref{f:3.37}$ in the last one. Also,
\begin{equation}\label{f:3.26}
\begin{aligned}
\|S_\varepsilon(\psi_{2\varepsilon}\nabla u_0)\|_{L^2(\Omega)}
+ \|S_\varepsilon(\psi_{2\varepsilon}u_0)\|_{L^2(\Omega)}
\leq C\|u_0\|_{H^1(\Omega)}\leq C\big\{\|g\|_{H^{1/2}(\partial\Omega)}+\|F\|_{L^{\frac{2d}{d+2}}(\Omega)}\big\}.
\end{aligned}
\end{equation}

Consequently, collecting $\eqref{f:3.22}-\eqref{f:3.26}$ and $\eqref{pri:3.11}$, we obtain
\begin{equation*}
 \|w_\varepsilon\|_{H^1(\Omega)} \leq C\varepsilon^{\frac{1}{2}}\big\{\|F\|_{L^{\frac{2d}{d+1}}(\Omega)} + \|g\|_{H^1(\partial\Omega)} \big\},
\end{equation*}
and the proof is finished.
\qed
\end{pf}

\begin{flushleft}
\textbf{Proof of Theorem \ref{thm:1.1}}\textbf{.}
This theorem includes the estimates $\eqref{pri:1.1}$ and $\eqref{pri:1.2}$, and we thus divide the proof into two parts.
\end{flushleft}

\textbf{Part I}, we prove the estimate $\eqref{pri:1.1}$. Let
\begin{equation}
 w_\varepsilon = u_\varepsilon - u_0 - \varepsilon\chi_{0,\varepsilon}S_\varepsilon(\psi_{4\varepsilon}u_0)
 + \varepsilon\chi_{k,\varepsilon}S_\varepsilon(\psi_{4\varepsilon}\nabla_ku_0),
\end{equation}
where $u_\varepsilon$ and $u_0$ satisfy $\eqref{pde:1.2}$, and $\psi_{4\varepsilon}$ is cut-off function
defined in Subsection $\ref{sec:1.1}$.
Then we have
\begin{equation*}
\begin{aligned}
 \|u_\varepsilon - u_0\|_{L^2(\Omega)}
 &\leq \|w_\varepsilon\|_{L^2(\Omega)} + \varepsilon\big\|\chi_{0,\varepsilon}S_\varepsilon(\psi_{4\varepsilon}u_0)\big\|_{L^2(\Omega)}
 + \varepsilon\big\|\chi_{k,\varepsilon}S_\varepsilon(\psi_{4\varepsilon}\nabla_ku_0)\big\|_{L^2(\Omega)}\\
 &\leq C\varepsilon\ln(c_0/\varepsilon)\Big\{\|F\|_{L^2(\Omega)} + \|g\|_{H^1(\partial\Omega)}\Big\} + C\varepsilon\|u_0\|_{H^1(\Omega)} \\
 &\leq C\varepsilon\ln(c_0/\varepsilon)\Big\{\|F\|_{L^2(\Omega)} + \|g\|_{H^1(\partial\Omega)}\Big\},
\end{aligned}
\end{equation*}
where we employ Theorem $\ref{thm:3.2}$ and the estimate $\eqref{pri:2.5}$ in the second inequality, and the estimate $\eqref{pri:2.3}$ is
used in the last one. We have already proved the estimate $\eqref{pri:1.1}$.

\textbf{Part II}, we now proceed to prove the estimate $\eqref{pri:1.2}$. Let $\tilde{w}_\varepsilon$ be given by
\begin{equation*}
 \tilde{w}_\varepsilon = u_\varepsilon - u_0 - \varepsilon\chi_{0,\varepsilon}S_\varepsilon(\psi_{4\varepsilon}\myu{u})
 + \varepsilon\chi_{k,\varepsilon}S_\varepsilon(\psi_{4\varepsilon}\nabla_k\myu{u}),
\end{equation*}
where $\myu{u}$ is the extension of $u_0$. For any $\Phi\in L^{q}(\Omega;\mathbb{R}^m)$ with $q=\frac{2d}{d+1}$,
there exist $\phi_\varepsilon,\phi_0$ solving $\eqref{pde:2.6}$ (due to Theorem $\ref{thm:2.1}$),
and satisfying $H^1$ estimate (see (\eqref{pri:2.3}))
\begin{equation}\label{f:3.89}
   \max\big\{\|\phi_\varepsilon\|_{H^1(\Omega)},\|\phi_0\|_{H^1(\Omega)}\big\}
   \leq C\|\Phi\|_{L^{\frac{2d}{d+2}}(\Omega)} \leq C\|\Phi\|_{L^q(\Omega)},
\end{equation}
where we use H\"older's inequality in the last inequality. Set
\begin{equation}\label{eq:3.1}
\Theta_\varepsilon = \phi_\varepsilon - \phi_0 - \varepsilon \chi_{0,\varepsilon}^* S_\varepsilon^2(\psi_{20\varepsilon}\phi_0)
- \varepsilon \chi_{k,\varepsilon}^* S_\varepsilon^2(\psi_{20\varepsilon}\nabla_k\phi_0),
\end{equation}
where $\chi_i^*$ with $i=0,\cdots,d$ are corresponding correctors of $\mathcal{L}_\varepsilon^*$.
Then it follows from the estimate $\eqref{pri:3.2}$ that $\|\Theta_\varepsilon\|_{H^1(\Omega)}\leq C\varepsilon^{1/2}\|\Phi\|_{L^q(\Omega)}$.
Moreover, in view of $\eqref{eq:3.1}$ again, we have
\begin{equation}\label{f:3.90}
 \|\phi_\varepsilon\|_{H^1(\Omega\setminus\Sigma_{9\varepsilon})}
 \leq \|\Theta_\varepsilon\|_{H^1(\Omega\setminus\Sigma_{9\varepsilon})} + \|\phi_0\|_{H^1(\Omega\setminus\Sigma_{9\varepsilon})}
 \leq \|\Theta_\varepsilon\|_{H^1(\Omega)} + \|\phi_0\|_{H^1(\Omega\setminus\Sigma_{9\varepsilon})}
 \leq C\varepsilon^{\frac{1}{2}}\|\Phi\|_{L^q(\Omega)},
\end{equation}
where we note that $\chi_{0,\varepsilon}^* S_\varepsilon^2(\psi_{20\varepsilon}\phi_0)$ and
$\chi_{k,\varepsilon}^* S_\varepsilon^2(\psi_{20\varepsilon}\nabla_k\phi_0)$ are supported in $\Sigma_{18\varepsilon}$,
and we use the estimate $\eqref{pri:3.11}$ in the last inequality.

In view of Lemma $\ref{lemma:2.11}$, we have
\begin{equation*}
\Big|\int_\Omega \tilde{w}_\varepsilon\Phi dx\Big|
\leq C\Big\{\|u_0\|_{H^1(\Omega\setminus\Sigma_{8\varepsilon})}\|\phi_\varepsilon\|_{H^1(\Omega\setminus\Sigma_{9\varepsilon})}
+ \varepsilon\|u_0\|_{H^2(\Omega)}\|\phi_\varepsilon\|_{H^1(\Omega)}\Big\},
\end{equation*}
and this together with $\eqref{f:3.89}$ and $\eqref{f:3.90}$ leads to
\begin{equation}\label{f:3.91}
\Big|\int_\Omega \tilde{w}_\varepsilon\Phi dx\Big| \leq C\varepsilon\|u_0\|_{H^2(\Omega)}\|\Phi\|_{L^q(\Omega)},
\end{equation}
where we use the following fact (due to the assumption of $u_0\in H^2(\Omega;\mathbb{R}^m)$),
\begin{equation}
\begin{aligned}
 \|u_0\|_{H^1(\Omega\setminus\Sigma_{8\varepsilon})}
 & \leq C\varepsilon^{1/2}\big\{\|\mathcal{M}(u_0)\|_{L^2(\partial\Omega)}
 + \|\mathcal{M}(\nabla u_0)\|_{L^2(\partial\Omega)}\big\} \\
 &\leq C\varepsilon^{1/2}\big\{\|u_0\|_{H^1(\Omega\setminus\Sigma_{c_0})}
 + \|u_0\|_{H^2(\Omega\setminus\Sigma_{c_0})}\big\}
 \leq C\varepsilon^{1/2}\|u_0\|_{H^2(\Omega)}.
\end{aligned}
\end{equation}
We mention that the estimate $\eqref{pri:2.8}$ is used in the first inequality,
and we apply the estimate $\eqref{pri:2.9}$ to the second one.

We turn back to the estimate $\eqref{f:3.91}$, and this implies
\begin{equation}
 \|\tilde{w}_\varepsilon\|_{L^{p}(\Omega)} \leq C\varepsilon\|u_0\|_{H^2(\Omega)},
\end{equation}
with $p=\frac{2d}{d-1}$.

The reminder of the argument is analogous to that shown in Part I. Before approaching the estimate $\eqref{pri:1.2}$,
we need to pose some simple estimates based on the Sobolev embedding theorem and interpolation inequality.
Let $2^* = \frac{2d}{d-2}$, and for all $k=0,\cdots,d$, we first have
\begin{equation}\label{f:3.92}
\|\chi_k\|_{L^p(Y)}\leq C\|\chi_k\|_{L^{2^*}(Y)}\leq C\|\chi_k\|_{H^1(Y)} \leq C(\mu,\kappa,m,d),
\end{equation}
where we use H\"older's inequality in the first inequality, and the Sobolev embedding theorem in second one. Besides, we have
\begin{equation}\label{f:3.93}
\begin{aligned}
& \left\{\begin{aligned}
&\|\myu{u}\|_{L^p(\mathbb{R}^d)}\leq \|\myu{u}\|_{L^2(\mathbb{R}^d)}^{\frac{1}{2}}\|\myu{u}\|_{L^{2^*}(\mathbb{R}^d)}^{\frac{1}{2}}
\leq C\|\myu{u}\|_{L^2(\mathbb{R}^d)}^{\frac{1}{2}}\|\nabla\myu{u}\|_{L^{2}(\mathbb{R}^d)}^{\frac{1}{2}}
\leq C\|\myu{u}\|_{H^1(\mathbb{R}^d)}, \\
&\|\nabla\myu{u}\|_{L^p(\mathbb{R}^d)}\leq \|\nabla\myu{u}\|_{L^2(\mathbb{R}^d)}^{\frac{1}{2}}\|\nabla\myu{u}\|_{L^{2^*}(\mathbb{R}^d)}^{\frac{1}{2}}
\leq C\|\nabla\myu{u}\|_{L^2(\mathbb{R}^d)}^{\frac{1}{2}}\|\nabla^2\myu{u}\|_{L^{2}(\mathbb{R}^d)}^{\frac{1}{2}}
\leq C\|\myu{u}\|_{H^2(\mathbb{R}^d)}
\end{aligned}\right. \\
&\Longrightarrow \qquad \|\myu{u}\|_{W^{1,p}(\mathbb{R}^d)}\leq C\|\myu{u}\|_{H^{2}(\mathbb{R}^d)}\leq C\|u_0\|_{H^2(\Omega)},
\end{aligned}
\end{equation}
where we use the interpolation inequality and Sobolev's inequality.

It is clear to see that
\begin{equation*}
\begin{aligned}
 \|u_\varepsilon - u_0\|_{L^p(\Omega)}
 &\leq \|\tilde{w}_\varepsilon\|_{L^p(\Omega)}
 + \varepsilon\big\|\chi_{0,\varepsilon}S_\varepsilon(\psi_{4\varepsilon}\myu{u})\big\|_{L^p(\Omega)}
 + \varepsilon\big\|\chi_{k,\varepsilon}S_\varepsilon(\psi_{4\varepsilon}\nabla_k\myu{u})\big\|_{L^p(\Omega)}\\
 &\leq C\varepsilon\|u_0\|_{H^2(\Omega)}
 + C\varepsilon\Big\{\|\chi_0\|_{L^p(Y)}\|\myu{u}\|_{L^p(\mathbb{R}^d)} + \|\chi_k\|_{L^p(Y)}\|\nabla_k\myu{u}\|_{L^p(\mathbb{R}^d)}\Big\}
 \leq C\varepsilon\|u_0\|_{H^2(\Omega)},
\end{aligned}
\end{equation*}
where we use the estimates $\eqref{f:3.92}$ and $\eqref{f:3.93}$ in the last inequality, and we have completed the proof.
\qed

\section{Neumann problem}

\begin{lemma}\label{lemma:4.1}
Suppose that the coefficients of $\mathcal{L}_\varepsilon$ satisfy $\eqref{a:1}-\eqref{a:3}$.
Let $u_\varepsilon$ and $u_0$ be the weak solutions to $\eqref{pde:1.3}$, and $w_\varepsilon$ is defined in $\eqref{eq:2.12}$.
Then we have
\begin{equation}\label{pri:4.1}
\begin{aligned}
\|w_{\varepsilon}\|_{H^1(\Omega)}
&\leq C\Big\{\|\varpi_{\varepsilon}\vec{\phi}\|_{L^2(\Omega\setminus\Sigma_{2\varepsilon})}
+ \|\nabla u_0 - \vec{\varphi}\|_{L^2(\Omega)}
+\|u_0 - \varphi_0\|_{L^2(\Omega)} \\
&+\varepsilon\|\varpi_{\varepsilon}\nabla\vec{\phi}\|_{L^2(\Omega)}
+ \varepsilon \|\varpi_{\varepsilon}\vec{\phi}\|_{L^2(\Omega)}\Big\}
+C\varepsilon\|h_\varepsilon\vec{\phi}\|_{L^2(\partial\Omega)},
\end{aligned}
\end{equation}
where $\vec{\varphi}=(\varphi_1,\cdots,\varphi_d)$, $\vec{\phi} = (\varphi_0,\vec{\varphi})$,
and $\varpi$ represents the periodic function depending on the periodic functions
such as the coefficients of $\myl{L}{\varepsilon}$,
the correctors $\{\chi_{k}\}_{k=0}^{d}$, and auxiliary functions
$\{b_{ik}, E_{jik},\nabla\vartheta_k\}_{k=0}^d$.
\end{lemma}

\begin{pf}
See \cite[Lemma 5.2]{QXS1}.
\qed
\end{pf}

\begin{lemma}\label{lemma:4.2}
Suppose that $\mathcal{L}_0$ is the homogenized operator of $\mathcal{L}_\varepsilon$ under the same conditions as in Lemma $\ref{lemma:4.1}$,
and we additionally assume $A=A^*$.
Let $u_0$ be the solution to $(\mathbf{NH})_0$ in $\eqref{pde:1.3}$ with $F\in L^{2}(\Omega;\mathbb{R}^m)$
and $h\in L^2(\partial\Omega;\mathbb{R}^m)$. Then we have
\begin{equation}\label{pri:4.2}
 \|u_0\|_{H^1(\Omega\setminus\Sigma_{p_1\varepsilon})} \leq C\varepsilon^{\frac{1}{2}}\big\{\|F\|_{L^{2}(\Omega)}
 + \|h\|_{L^2(\partial\Omega)}\big\},
\end{equation}
and
\begin{equation}\label{pri:4.3}
 \|\nabla^2 u_0\|_{L^2(\Sigma_{p_2\varepsilon})} \leq C\varepsilon^{-\frac{1}{2}}\big\{\|F\|_{L^2(\Omega)}
 + \|h\|_{L^2(\partial\Omega)}\big\},
\end{equation}
where $p_1,p_2>0$ are fixed real numbers, and $C$ depends on $\mu,\kappa,m,d,p_1,p_2$ and $\Omega$.
\end{lemma}

\begin{remark}\label{re:3.1}
\emph{The results of $\eqref{pri:4.2}$ and $\eqref{pri:4.3}$ were originally
established by Z. Shen in \cite[Theorem 2.7]{SZW2} for $L_\varepsilon$ with the Neumann boundary condition.
Here we employ the radial maximal function to extend his results to our cases.}
\end{remark}

\begin{pf}
We first mention that this lemma is the counterpart of Lemma $\ref{lemma:3.5}$ in Neumann problems,
which was actually proved in \cite[Lemma 5.6]{QXS1}. We provide the reader with a proof for the sake of the completeness.
So, proceeding as in the proof Lemma $\ref{lemma:3.5}$, we rewrite $(\mathbf{NH})_0$ as
\begin{equation*}
 L_0(u_0) = (\widehat{V} -\widehat{B})\nabla u_0 - (\widehat{c}+\lambda I)u_0 + F \quad \text{in}~\Omega,
 \qquad \partial u_0/\partial\nu_0 = h - n\cdot\widehat{V}u_0 \quad \text{on}~\partial\Omega,
\end{equation*}
where $L_0 = \text{div}(\widehat{A}\nabla)$, and then consider $u_0 = v + w$ such that
\begin{equation}\label{pde:4.1}
 (1)~ L_0(v) = \breve{F} \quad \text{in}~\mathbb{R}^d,
 \qquad\quad
 (2)~ \left\{\begin{aligned}
 L_0(w) &= 0  &\quad&\text{in}~~\Omega,\\
     \partial w/\partial\nu_0  &= h-n\cdot\widehat{V}u_0 - \partial v/\partial\nu_0 &\quad&\text{on}~\partial\Omega,
 \end{aligned}\right.
\end{equation}
where $\partial/\partial\nu_0 = n\cdot\widehat{A}\nabla$, and $\breve{F}$ is the same thing in $\eqref{pde:3.3}$.
For (1), we apply the estimate $\eqref{pri:2.4}$ to the term of $\|u_0\|_{H^1(\Omega)}$ in $\eqref{f:3.27}$ and $\eqref{f:3.28}$,
and then obtain
\begin{equation}\label{f:4.1}
\|\nabla^2 v\|_{L^{p}(\mathbb{R}^d)}
\leq C\big\{\|F\|_{L^p(\Omega)}+\|h\|_{H^{-1/2}(\partial\Omega)}\big\},
\qquad\|\nabla v\|_{L^2(\mathbb{R}^d)}
\leq C\big\{\|F\|_{L^{\frac{2d}{d+2}}(\Omega)} + \|h\|_{H^{-1/2}(\partial\Omega)}\big\}
\end{equation}
where $p\in[\frac{2d}{d+2},2]$. Thus by Sobolev's inequality, we have
\begin{equation}\label{f:4.2}
 \|v\|_{L^{\frac{2d}{d-2}}(\mathbb{R}^d)}
 \leq C \big\{\|F\|_{L^{\frac{2d}{d+2}}(\Omega)} + \|h\|_{H^{-1/2}(\partial\Omega)} \big\},
 \qquad
 \|\nabla v\|_{L^{\frac{2d}{d-1}}(\mathbb{R}^d)}
 \leq C \big\{\|F\|_{L^{\frac{2d}{d+1}}(\Omega)}+\|h\|_{H^{-1/2}(\partial\Omega)}\big\}.
\end{equation}
Then plugging $\eqref{f:4.1}$ and $\eqref{f:4.2}$ back into the third line of $\eqref{f:3.94}$ leads to
\begin{equation}\label{f:4.3}
\int_{S_t} \big(|\nabla v|^2 + |v|^2\big) dS
\leq C\big\{\|F\|_{L^{\frac{2d}{d+1}}(\Omega)}^2+\|h\|_{H^{-1/2}(\partial\Omega)}^2\big\}
\end{equation}
for any $t\in[0,p_1\varepsilon]$. and this together with the co-area formula $\eqref{eq:2.9}$ shows
\begin{equation}\label{f:4.4}
\|v\|_{H^1(\Omega\setminus\Sigma_{p_1\varepsilon})}
\leq C\varepsilon^{\frac{1}{2}}\big\{\|F\|_{L^{\frac{2d}{d+1}}(\Omega)}+\|h\|_{L^{2}(\partial\Omega)}\big\},
\end{equation}
where we use the fact of
$L^2(\partial\Omega;\mathbb{R}^m)\subset H^{-1/2}(\partial\Omega;\mathbb{R}^m)$.

The next thing is to estimate $\|w\|_{H^1(\Omega\setminus\Sigma_{p_1\varepsilon})}$. For (2), in view of \cite[Remark 5.7]{QXS1}, we have
\begin{equation}\label{f:4.5}
 \|(\nabla w)^*\|_{L^2(\partial\Omega)}
 \leq C\big\{\|h\|_{L^2(\partial\Omega)}+\|u_0\|_{L^2(\partial\Omega)} + \|\nabla v\|_{L^2(\partial\Omega)}\big\}
 \leq C\big\{\|F\|_{L^{\frac{2d}{d+1}}(\Omega)} + \|h\|_{L^2(\partial\Omega)}\big\}
\end{equation}
where we employ the estimate $\eqref{f:4.3}$ (for $t=0$) and the trace theorem coupled with the estimate $\eqref{pri:2.4}$
in the last inequality. On the other hand, we apply Lemma $\ref{lemma:2.5}$ and the estimate $\eqref{pri:2.4}$  to derive
\begin{equation}\label{f:4.6}
 \|\mathcal{M}(w)\|_{L^2(\partial\Omega)}
 \leq C\big\{\|h\|_{L^{2}(\partial\Omega)} + \|u_0\|_{L^{2}(\partial\Omega)} + \|\nabla v\|_{L^2(\partial\Omega)}\big\}
 \leq C\big\{\|F\|_{L^{\frac{2d}{d+1}}(\Omega)}+\|h\|_{L^{2}(\partial\Omega)}\big\},
\end{equation}
Hence, it follows from $\eqref{pri:2.8}$, $\eqref{f:4.5}$ and $\eqref{f:4.6}$ that
\begin{equation}\label{f:4.7}
\|w\|_{H^1(\Omega\setminus\Sigma_{p_1\varepsilon})}
\leq C\varepsilon^{\frac{1}{2}}
\big\{\|(\nabla w)^*\|_{L^2(\partial\Omega)} + \|\mathcal{M}(w)\|_{L^2(\partial\Omega)}\big\}
\leq C\varepsilon^{\frac{1}{2}}\big\{\|F\|_{L^{\frac{2d}{d+1}}(\Omega)}+\|h\|_{L^{2}(\partial\Omega)}\big\}.
\end{equation}
Combining the estimates $\eqref{f:4.4}$ and $\eqref{f:4.7}$, we have
\begin{equation}\label{pri:4.10}
\|u_0\|_{H^1(\Omega\setminus\Sigma_{p_1\varepsilon})}
\leq C\varepsilon^{\frac{1}{2}}\big\{\|F\|_{L^{\frac{2d}{d+1}}(\Omega)}+\|h\|_{L^{2}(\partial\Omega)}\big\},
\end{equation}
and this coupled with H\"older's inequality leads to the estimate $\eqref{pri:4.2}$.

We are now in the position to prove the estimate $\eqref{pri:4.3}$. Thanks to the first line of the estimate $\eqref{f:3.37}$,
we actually have
 \begin{equation}
   \|\nabla^2 w\|_{L^2(\Sigma_{p_2\varepsilon})}
   \leq C\varepsilon^{-\frac{1}{2}}\|(\nabla w)^*\|_{L^2(\partial\Omega)} + C\|\nabla w\|_{L^2(\Omega)}
   \leq C\varepsilon^{-\frac{1}{2}}\big\{\|F\|_{L^{\frac{2d}{d+1}}(\partial\Omega)}+\|h\|_{L^2(\partial\Omega)}\big\},
 \end{equation}
 where we use the estimates $\eqref{f:4.5}$ and $\eqref{pri:2.4}$ in the second inequality.
 This together with the estimate $\eqref{f:4.1}$ (with $p=2$) gives the desired estimate $\eqref{pri:4.3}$,
 and we have completed the proof.
\qed
\end{pf}

\begin{lemma}[Improved lemma]\label{lemma:4.3}
Assume the same conditions as in Lemma $\ref{lemma:4.2}$.
Let $u_0$ be the solution to $(\mathbf{NH})_0$ in $\eqref{pde:1.3}$ with $F\in L^{2}(\Omega;\mathbb{R}^m)$
and $h\in L^2(\partial\Omega;\mathbb{R}^m)$. Then we have
\begin{equation}\label{pri:4.4}
\|u_0\|_{H^1(\Omega\setminus\Sigma_{p_1\varepsilon};\delta)} \leq C\varepsilon\big\{\|F\|_{L^2(\Omega)}
 + \|h\|_{L^2(\partial\Omega)}\big\},
\end{equation}
and
\begin{equation}\label{pri:4.5}
\max\big\{\|\nabla^2 u_0\|_{L^2(\Sigma_{p_2\varepsilon};\delta)},
~\|\nabla u_0\|_{L^2(\Sigma_{p_2\varepsilon};\delta^{-1})}\big\}
\leq C\big[\ln(c_0/\varepsilon)\big]^{1/2}\big\{\|F\|_{L^2(\Omega)}
 + \|h\|_{L^2(\partial\Omega)}\big\},
\end{equation}
where $p_1,p_2>0$ are fixed real numbers and $c_0$ is the layer constant, and $C$ depends on $\mu,\kappa,m,d,p_1$ and $\Omega$.
\end{lemma}

\begin{pf}
The proof is quite similar to that shown in Lemma $\ref{lemma:3.6}$,
so we can straightforward use some key inequalities obtained there to prove this one.
The estimate $\eqref{pri:4.4}$ follows from $\eqref{pri:2.19}$ and $\eqref{pri:4.2}$ that
\begin{equation*}
 \|u_0\|_{H^1(\Omega\setminus\Sigma_{p_1\varepsilon};\delta)}
 \leq C\varepsilon^{1/2}\|u_0\|_{H^1(\Omega\setminus\Sigma_{p_1\varepsilon})}
 \leq C\varepsilon\big\{\|F\|_{L^2(\Omega)}+\|h\|_{L^2(\partial\Omega)}\big\}.
\end{equation*}

We now study the term of $\|\nabla^2 u_0\|_{L^2(\Sigma_{p_2\varepsilon};\delta)}$.
Thanks to the first line of $\eqref{f:3.74}$ and the second line of $\eqref{f:3.75}$, we have
\begin{equation}\label{f:4.8}
\begin{aligned}
 \|\nabla^2 u_0\|_{L^2(\Sigma_{p_2\varepsilon};\delta)}
 &\leq C\|\nabla^2 v\|_{L^2(\mathbb{R}^d)} + \|\nabla^2 w\|_{L^2(\Sigma_{p_2\varepsilon};\delta)} \\
 &\leq C\|\nabla^2 v\|_{L^2(\mathbb{R}^d)} + C\big[\ln(c_0/\varepsilon)\big]^{\frac{1}{2}}\|(\nabla w)^*\|_{L^2(\partial\Omega)}
 \leq C\big[\ln(c_0/\varepsilon)\big]^{\frac{1}{2}}\big\{\|F\|_{L^2(\Omega)}+\|h\|_{L^2(\partial\Omega)}\big\},
\end{aligned}
\end{equation}
where we use the estimate $\eqref{f:4.1}$ (for $p=2$) and the estimate $\eqref{f:4.5}$
in the last inequality.

The rest thing is to estimate the term of $\|u_0\|_{H^1(\Sigma_{p_2\varepsilon};\delta^{-1})}$. Due to the estimate $\eqref{f:3.95}$, we have
\begin{equation}\label{f:4.9}
\begin{aligned}
\|u_0\|_{H^1(\Sigma_{p_2\varepsilon};\delta^{-1})}
&\leq C\big[\ln(c_0/\varepsilon)\big]^{1/2}\Big\{\|\mathcal{M}(u_0)\|_{L^2(\partial\Omega)} +
\|\mathcal{M}(\nabla u_0)\|_{L^2(\partial\Omega)}
 + \|u_0\|_{H^1(\Omega)}\Big\} \\
&\leq C\big[\ln(c_0/\varepsilon)\big]^{1/2}\big\{\|F\|_{L^2(\Omega)}+\|h\|_{L^2(\partial\Omega)}\big\},
\end{aligned}
\end{equation}
where we use the estimates $\eqref{pri:2.4}$, $\eqref{f:4.5}$ and $\eqref{f:4.6}$ in the last inequality.

Combining $\eqref{f:4.8}$ and $\eqref{f:4.9}$ leads to the desired estimate $\eqref{pri:4.5}$.
We have completed the proof.
\qed
\end{pf}

\begin{thm}\label{thm:4.1}
Suppose that the coefficients of $\mathcal{L}_\varepsilon$ satisfy $\eqref{a:1}-\eqref{a:3}$,
and we additionally assume $A=A^*$.
Let $u_\varepsilon,u_0$ be the weak solutions to $\eqref{pde:1.2}$
with $F\in L^q(\Omega;\mathbb{R}^m)$ and $h\in L^2(\partial\Omega;\mathbb{R}^m)$, where $q=\frac{2d}{d+1}$. Then we have
\begin{equation}\label{pri:4.6}
\big\|u_\varepsilon - u_0 -\varepsilon\chi_{0,\varepsilon}S^2_\varepsilon(\psi_{2\varepsilon}u_0)
-\varepsilon\chi_{k,\varepsilon}S^2_\varepsilon(\psi_{2\varepsilon}\nabla_k u_0)\big\|_{H^1(\Omega)}
\leq C\varepsilon^{\frac{1}{2}}\big\{\|F\|_{L^q(\Omega)}+\|h\|_{L^2(\partial\Omega)}\big\},
\end{equation}
where $C$ depends only on $\mu,\kappa,m,d$ and $\Omega$.
\end{thm}

\begin{pf}
By suitable modification to the proof of Theorem $\ref{thm:3.3}$,
it is not hard to prove this one, which in fact was shown in \cite[Lemma 5.6]{QXS1}, so we omit this proof.
\qed
\end{pf}

\begin{cor}\label{cor:4.1}
Assume the same conditions as in Theorem $\ref{thm:4.1}$. Let $u_\varepsilon$ and $u_0$ be the
weak solutions to $\eqref{pde:1.3}$ with $F\in L^2(\Omega;\mathbb{R}^m)$ and $h\in L^2(\partial\Omega;\mathbb{R}^m)$.
Then we have
\begin{equation}\label{pri:4.7}
\big\|u_\varepsilon - u_0 -\varepsilon\chi_{0,\varepsilon}S_\varepsilon(\psi_{4\varepsilon}u_0)
+\varepsilon\chi_{k,\varepsilon}S_\varepsilon(\psi_{4\varepsilon}\nabla_k u_0)\big\|_{H^1(\Omega)}
\leq C\varepsilon^{\frac{1}{2}}\big\{\|F\|_{L^2(\Omega)} + \|h\|_{L^2(\partial\Omega)}\big\},
\end{equation}
where $C$ depends only on $\mu,\kappa,m,d$ and $\Omega$.
\end{cor}

\begin{pf}
By the triangle inequality, we have
\begin{equation}\label{f:4.10}
\begin{aligned}
\big\|u_\varepsilon - u_0
&-\varepsilon\chi_{0,\varepsilon}S_\varepsilon(\psi_{4\varepsilon}u_0)
+\varepsilon\chi_{k,\varepsilon}S_\varepsilon(\psi_{4\varepsilon}\nabla_k u_0)\big\|_{H^1(\Omega)}
\leq \varepsilon\big\|\chi_{k,\varepsilon}S_\varepsilon\big(S_\varepsilon(\psi_{4\varepsilon}\nabla_ku_0)
- \psi_{4\varepsilon}\nabla_ku_0\big)\big\|_{H^1(\Omega)} \\
&  + \varepsilon\big\|\chi_{0,\varepsilon}S_\varepsilon\big(S_\varepsilon(\psi_{4\varepsilon}u_0)
 - \psi_{4\varepsilon}u_0\big)\big\|_{H^1(\Omega)}
+ \big\|u_\varepsilon - u_0 -\varepsilon\chi_{0,\varepsilon}S^2_\varepsilon(\psi_{4\varepsilon}u_0)
+\varepsilon\chi_{k,\varepsilon}S^2_\varepsilon(\psi_{4\varepsilon}\nabla_k u_0)\big\|_{H^1(\Omega)}.
\end{aligned}
\end{equation}
The next thing is to handle the first two terms in the right-hand side of $\eqref{f:4.10}$. It is not hard to derive the following estimates:
\begin{equation}\label{f:4.11}
\begin{aligned}
\big\|\chi_{k,\varepsilon}S_\varepsilon\big(&S_\varepsilon(\psi_{4\varepsilon}\nabla_ku_0)
- \psi_{4\varepsilon}\nabla_ku_0\big)\big\|_{H^1(\Omega)}\\
&\leq \big\|\chi_{k,\varepsilon}S_\varepsilon\big(S_\varepsilon(\psi_{4\varepsilon}\nabla_ku_0)
- \psi_{4\varepsilon}\nabla_ku_0\big)\big\|_{L^2(\Omega)}
+ \varepsilon^{-1}\big\|(\nabla\chi_{k})_\varepsilon S_\varepsilon\big(S_\varepsilon(\psi_{4\varepsilon}\nabla_ku_0)
- \psi_{4\varepsilon}\nabla_ku_0\big)\big\|_{L^2(\Omega)}\\
&\qquad \qquad+ \big\|\chi_{k,\varepsilon}S_\varepsilon\big(\nabla(S_\varepsilon(\psi_{4\varepsilon}\nabla_ku_0))
- \nabla(\psi_{4\varepsilon}\nabla_ku_0)\big)\big\|_{L^2(\Omega)}\\
&\leq C\big\{1+\varepsilon^{-1}\big\}\big\|S_\varepsilon(\psi_{4\varepsilon}\nabla u_0)
- \psi_{4\varepsilon}\nabla u_0\big\|_{L^2(\mathbb{R}^d)}
+ C\big\{\|\nabla S_\varepsilon(\psi_{4\varepsilon}\nabla u_0)\|_{L^2(\mathbb{R}^d)}
+ \|\nabla(\psi_{4\varepsilon}\nabla u_0)\|_{L^2(\mathbb{R}^d)}\big\} \\
&\leq C\|\nabla(\psi_{4\varepsilon}\nabla u_0)\|_{L^2(\mathbb{R}^d)}
\leq C\big\{\varepsilon^{-1}\|\nabla u_0\|_{L^2(\Omega\setminus\Sigma_{8\varepsilon})}
+ \|\nabla^2 u_0\|_{L^2(\Sigma_{4\varepsilon})}\big\}\\
& \leq C\varepsilon^{-1/2}\big\{\|F\|_{L^2(\Omega)} + \|h\|_{L^2(\partial\Omega)}\big\},
\end{aligned}
\end{equation}
where we use the estimate $\eqref{pri:2.5}$ in the second inequality,
and the estimates $\eqref{pri:2.5}$ and $\eqref{pri:2.6}$ in the third one. For the last one,
we employ the estimates $\eqref{pri:4.2}$ and $\eqref{pri:4.3}$. By the same token, we can derive
\begin{equation}\label{f:4.12}
  \big\|\chi_{0,\varepsilon}S_\varepsilon\big(S_\varepsilon(\psi_{4\varepsilon}u_0)
 - \psi_{4\varepsilon}u_0\big)\big\|_{H^1(\Omega)}
 \leq C\|\nabla(\psi_{4\varepsilon}u_0)\|_{L^2(\mathbb{R}^d)}
 \leq C\varepsilon^{-1/2}\big\{\|F\|_{L^2(\Omega)} + \|h\|_{L^2(\partial\Omega)}\big\}.
\end{equation}

Collecting the estimates $\eqref{pri:4.6}$, $\eqref{f:4.10}$, $\eqref{f:4.11}$ and $\eqref{f:4.12}$,
we have the desired estimate $\eqref{pri:4.7}$.
In the end, we mention that the proof also follows from Lemmas $\ref{lemma:4.1}$ and $\ref{lemma:4.2}$,
and we have completed the proof.
\qed

\end{pf}

\begin{thm}\label{thm:4.2}
Suppose that the coefficients of $\mathcal{L}_\varepsilon$ satisfy the same conditions as in Theorem $\ref{thm:4.1}$.
Let $u_\varepsilon,u_0$ be the solutions to $\eqref{pde:1.3}$.
\begin{itemize}
  \item [\emph{(i)}] If $F\in L^2(\Omega;\mathbb{R}^m)$ and $h\in L^2(\partial\Omega;\mathbb{R}^m)$, then we have
  \begin{equation}\label{pri:4.8}
  \big\|u_\varepsilon -u_0 - \varepsilon\chi_{0,\varepsilon}S_\varepsilon(\psi_{4\varepsilon}u_0)
  -\varepsilon\chi_{k,\varepsilon}S_\varepsilon(\psi_{4\varepsilon}\nabla_k u_0)\big\|_{L^2(\Omega)}
  \leq C\varepsilon\ln(c_0/\varepsilon)\big\{\|F\|_{L^2(\Omega)} + \|h\|_{L^2(\partial\Omega)}\big\},
  \end{equation}
  where $c_0$ is the layer constant.
  \item [\emph{(ii)}] If $u_0\in H^2(\Omega;\mathbb{R}^m)$, then for $p=\frac{2d}{d-1}$, we have
  \begin{equation}\label{pri:4.9}
  \big\|u_\varepsilon -u_0 - \varepsilon\chi_{0,\varepsilon}S_\varepsilon(\psi_{4\varepsilon}\myu{u})
  -\varepsilon\chi_{k,\varepsilon}S_\varepsilon(\psi_{4\varepsilon}\nabla_k \myu{u})\big\|_{L^p(\Omega)}
  \leq C\varepsilon\|u_0\|_{H^2(\Omega)},
 \end{equation}
 where $\myu{u}$ is the extension of $u_0$ \emph{(}see Lemma $\ref{lemma:2.11}$\emph{)}.
\end{itemize}
Note that $C$ depends on $\mu,\kappa,m,d,c_0$ and $\Omega$.
\end{thm}

\begin{pf}
Again, the duality argument is employed in this theorem. (i). For any $\Phi\in L^2(\Omega;\mathbb{R}^m)$,
there exist two weak solutions $\phi_\varepsilon$ and $\phi_0$ satisfying $\eqref{pde:2.7}$.
According to Lemma $\ref{lemma:2.10}$, we construct
\begin{equation*}
\xi_\varepsilon = u_\varepsilon - u_0 - \varepsilon\chi_{0,\varepsilon}^*S_\varepsilon(\psi_{10\varepsilon}\phi_0)
- \varepsilon\chi_{k,\varepsilon}^*S_\varepsilon(\psi_{10\varepsilon}\nabla_k\phi_0).
\end{equation*}
Then it follows from Corollary $\ref{cor:4.1}$ that
\begin{equation}\label{f:4.13}
 \|\xi_\varepsilon\|_{H^1(\Omega)}\leq C\varepsilon^{1/2}\|\Phi\|_{L^2(\Omega)}.
\end{equation}

Hence the rest ting is to estimate the right-hand side of $\eqref{pri:2.18}$ term by term.
We may assume $\|F\|_{L^2(\Omega)} + \|h\|_{L^2(\partial\Omega)} = 1$ for ease of computations.
In view of the estimates $\eqref{pri:4.2}$ and $\eqref{pri:2.4}$, we have
\begin{equation}\label{f:4.14}
\|u_0\|_{H^1(\Omega\setminus\Sigma_{9\varepsilon})}\|\phi_0\|_{H^1(\Omega\setminus\Sigma_{9\varepsilon})}
 + \varepsilon\|u_0\|_{H^1(\Omega)}\|\phi_\varepsilon\|_{H^1(\Omega)}
 \leq C\varepsilon^{\frac{1}{2}}\cdot C\varepsilon^{\frac{1}{2}}\|\Phi\|_{L^2(\Omega)} + C\varepsilon\|\Phi\|_{L^2(\Omega)}
 \leq C\varepsilon\|\Phi\|_{L^2(\Omega)}.
\end{equation}
On account of the estimates $\eqref{pri:2.4}$, $\eqref{pri:4.4}$ and $\eqref{pri:4.5}$, we derive
\begin{equation}\label{f:4.15}
\begin{aligned}
\Big\{\| u_0\|_{H^1(\Omega\setminus\Sigma_{8\varepsilon};\delta)}
+\varepsilon \|\nabla u_0\|_{L^2(\Sigma_{4\varepsilon};\delta)}
&+\varepsilon \|\nabla^2 u_0\|_{L^2(\Sigma_{4\varepsilon};\delta)}\Big\}\|\phi_0\|_{H^1(\Sigma_{4\varepsilon};\delta^{-1})} \\
&\leq C\Big\{\varepsilon
+ \varepsilon\big[\ln(c_0/\varepsilon)\big]^{\frac{1}{2}}\Big\}\big[\ln(c_0/\varepsilon)\big]^{\frac{1}{2}}\|\Phi\|_{L^2(\Omega)}
\leq C\varepsilon\ln(c_0/\varepsilon)\|\Phi\|_{L^2(\Omega)}.
\end{aligned}
\end{equation}
Due to the estimates $\eqref{pri:2.4}$, $\eqref{pri:4.2}$, $\eqref{pri:4.3}$ and $\eqref{f:4.13}$, we arrive at
\begin{equation}\label{f:4.16}
\begin{aligned}
\Big\{\| u_0\|_{H^1(\Omega\setminus\Sigma_{8\varepsilon})}
&+\varepsilon \|\nabla u_0\|_{L^2(\Sigma_{4\varepsilon})}
+\varepsilon \|\nabla^2 u_0\|_{L^2(\Sigma_{4\varepsilon})}\Big\}\\
&\cdot\Big\{\|\xi_\varepsilon\|_{H^1(\Omega)} + \|\phi_0\|_{H^1(\Omega\setminus\Sigma_{20\varepsilon})}
+ \varepsilon\|\phi_0\|_{H^1(\Omega)}+ \varepsilon\|\nabla^2\phi_0\|_{L^2(\Sigma_{10\varepsilon})}\Big\}\\
& \qquad\qquad\qquad\qquad\qquad \leq C\Big\{\varepsilon + \varepsilon + \varepsilon^{1/2}\Big\}
\cdot\Big\{\varepsilon^{1/2}+\varepsilon + \varepsilon + \varepsilon^{1/2}\Big\}\|\Phi\|_{L^2(\Omega)}
\leq C\varepsilon \|\Phi\|_{L^2(\Omega)}.
\end{aligned}
\end{equation}
Consequently, collecting $\eqref{pri:2.18}$, $\eqref{f:4.13}-\eqref{f:4.16}$ gives the desired estimate $\eqref{pri:4.8}$.

(ii). The proof is similar to that shown in Part II of Theorem $\ref{thm:1.1}$.
Let $u_\varepsilon$ and $u_0$ be the solutions of
$\eqref{pde:1.3}$, and the corresponding $w_\varepsilon$ in Lemma $\ref{lemma:2.11}$ is given by
\begin{equation}\label{f:4.19}
  w_\varepsilon = u_\varepsilon -u_0 - \varepsilon\chi_{0,\varepsilon}S_\varepsilon(\psi_{4\varepsilon}\myu{u})
  -\varepsilon\chi_{k,\varepsilon}S_\varepsilon(\psi_{4\varepsilon}\nabla_k \myu{u}).
\end{equation}
For any $\Psi\in L^q(\Omega;\mathbb{R}^m)$ with $q = \frac{2d}{d+1}$,
there exist $\phi_\varepsilon$ and $\phi_0$ satisfying $\eqref{pde:2.7}$ and the $H^1$ estimates (see $\eqref{pri:2.4}$)
\begin{equation}\label{f:4.17}
\max\big\{\|\phi_\varepsilon\|_{H^1(\Omega)},\|\varphi_{0}\|_{H^1(\Omega)}\big\}
\leq C\|\Psi\|_{L^{\frac{2d}{d+2}}(\Omega)}\leq C\|\Psi\|_{L^{q}(\Omega)},
\end{equation}
where we use H\"older's inequality in the second inequality. According to Lemma $\ref{lemma:2.11}$, construct
\begin{equation*}
 \Pi_\varepsilon = \phi_\varepsilon - \phi_0 -\varepsilon\chi_{0,\varepsilon}^*S_\varepsilon^2(\psi_{20\varepsilon}\phi_0)
 -\varepsilon\chi_{k,\varepsilon}^*S_\varepsilon^2(\psi_{20\varepsilon}\nabla_k\phi_0),
\end{equation*}
and then we have
\begin{equation}\label{f:4.18}
\|\phi_\varepsilon\|_{H^1(\Omega\setminus\Sigma_{9\varepsilon})}
\leq \|\Pi_\varepsilon\|_{H^1(\Omega)} + \|\phi_0\|_{H^1(\Omega\setminus\Sigma_{9\varepsilon})}
\leq C\varepsilon^{1/2}\|\Psi\|_{L^q(\Omega)}.
\end{equation}
Note that $S_\varepsilon^2(\psi_{20\varepsilon}\phi_0)$ and $S_\varepsilon^2(\psi_{20\varepsilon}\nabla_k\phi_0)$ are supported in
$\Sigma_{18\varepsilon}$, which guarantees that the first inequality is valid,
and we use Theorem $\ref{thm:4.1}$ and the estimate $\eqref{pri:4.10}$ in the last inequality.

In view of Lemma $\ref{lemma:2.11}$, we plug the estimate $\eqref{f:4.17}$ and $\eqref{f:4.18}$
into $\eqref{pri:2.20}$, and obtain
\begin{equation*}
\Big|\int_\Omega w_\varepsilon\Psi dx\Big|
\leq C\varepsilon \|u_0\|_{H^2(\Omega)}\|\Psi\|_{L^q(\Omega)}.
\end{equation*}
This implies the desired estimate $\eqref{pri:4.9}$. We have completed the proof.
\qed
\end{pf}

\begin{flushleft}
\textbf{Proof of Theorem \ref{thm:1.2}}\textbf{.}
Based on Theorem $\ref{thm:4.2}$, the proof is almost the same one shown for Theorem $\ref{thm:1.1}$.
We first investigate the estimate $\eqref{pri:1.3}$. Let
\begin{equation*}
 w_\varepsilon = u_\varepsilon - u_0 - \varepsilon\chi_{0,\varepsilon}S_\varepsilon(\psi_{4\varepsilon}u_0)
 + \varepsilon\chi_{k,\varepsilon}S_\varepsilon(\psi_{4\varepsilon}\nabla_ku_0),
\end{equation*}
where $u_\varepsilon$ and $u_0$ satisfy $\eqref{pde:1.3}$. Hence,
\end{flushleft}
\begin{equation*}
\begin{aligned}
 \|u_\varepsilon - u_0\|_{L^2(\Omega)}
 &\leq \|w_\varepsilon\|_{L^2(\Omega)} + \varepsilon\big\|\chi_{0,\varepsilon}S_\varepsilon(\psi_{4\varepsilon}u_0)\big\|_{L^2(\Omega)}
 + \varepsilon\big\|\chi_{k,\varepsilon}S_\varepsilon(\psi_{4\varepsilon}\nabla_ku_0)\big\|_{L^2(\Omega)}\\
 &\leq C\varepsilon\ln(c_0/\varepsilon)\Big\{\|F\|_{L^2(\Omega)} + \|h\|_{L^2(\partial\Omega)}\Big\} + C\varepsilon\|u_0\|_{H^1(\Omega)} \\
 &\leq C\varepsilon\ln(c_0/\varepsilon)\Big\{\|F\|_{L^2(\Omega)} + \|h\|_{L^2(\partial\Omega)}\Big\},
\end{aligned}
\end{equation*}
where we employ Theorem $\ref{thm:4.2}$ and the estimate $\eqref{pri:2.5}$ in the second inequality, and the estimate $\eqref{pri:2.4}$ is
used in the last one.

We now turn to the estimate $\eqref{pri:1.4}$. Redefine $w_\varepsilon$ to be expressed by the right-hand side of $\eqref{f:4.19}$.
Then it follows from the estimates $\eqref{pri:4.9}$, $\eqref{f:3.92}$ and $\eqref{f:3.93}$ that
\begin{equation*}
 \|u_\varepsilon - u_0\|_{L^p(\Omega)}
 \leq \|w_\varepsilon\|_{L^p(\Omega)}
 + \varepsilon\big\|\chi_{0,\varepsilon}S_\varepsilon(\psi_{4\varepsilon}\myu{u})\big\|_{L^p(\Omega)}
 + \varepsilon\big\|\chi_{k,\varepsilon}S_\varepsilon(\psi_{4\varepsilon}\nabla_k\myu{u})\big\|_{L^p(\Omega)}
 \leq C\varepsilon\|u_0\|_{H^2(\Omega)}.
\end{equation*}
We have completed the proof.
\qed

\begin{center}
\textbf{Acknowledgements}
\end{center}

The author wishes to express his sincere appreciation to Professor Zhongwei Shen for his constant and illuminating instruction.
The author is greatly indebted to the referees for many useful comments and helpful suggestions.
The author also wants to express his heartfelt gratitude to Professor Peihao Zhao, who led him into the world of mathematics.
This work was supported by the National Natural Science Foundation of China (Grant NO.11471147).

%
%